\documentclass[12pt]{article}
\usepackage{local}

\newcounter{opagecount}
\let\oripage\relax
\renewcommand\thefootnote{(\arabic{footnote})}
\newcounter{footA}
\long\def\footletter#1{\stepcounter{footA}\renewcommand\thefootnote{(\Alph{footA})}\protect\footnote{#1}\renewcommand\thefootnote{(\arabic{footnote})}\addtocounter{footnote}{-1}}

\begin{document}
\renewcommand\bibname{}

\setcounter{page}{35}%
\setcounter{opagecount}{35}%

\title{Finiteness Problems in Diophantine Geometry}

\author{Yu.\,G.\,Zarkhin (Zarhin), A.\,N.\,Parshin}
\date{}
\maketitle

This survey contains an exposition of the results of the German mathematician
Gerd Faltings, who proved some longstanding conjectures in algebraic geometry
formulated by Shafarevich, Tate, and Mordell. The different parts of the proof
are described with varying degrees of detail. Sometimes we shall limit
ourselves to a sketch of the proof together with a discussion of the basic idea. Compared
with the original article by Faltings~\cite{39}, we have introduced several
additions and modifications contributed by the participants in Szpiro's
seminar~\cite{71} (\'Ecole Normale Sup\'erieure, Paris, 1983--1984) and by the
authors. \S1 contains the general plan of proof, and also a history of what was
known before Faltings. Concerning other work related to these conjectures,
see~\cite{19}, \cite{20}, and the basic text of Lang's book. Sections 2--6
contain all of the necessary definitions, of which we make free use in \S1, and
the details of the proofs. We conclude the survey by giving some unsolved
problems.

\section{Introduction}

The point of departure for the circle of problems we are considering was two conjectures proposed
by Shafarevich at the 1962 International Congress of Mathematicians in Stockholm~\cite{30}. He
examined the question of classifying algebraic curves $X$ of fixed genus $g\geq1$ over an
algebraic number field $K$. Besides these invariants, a curve $X$ also has its set $S$ of places
of bad reduction (see \S6). The first conjecture of Shafarevich (the \textit{finiteness
conjecture}) asserts that, if $g>1$ (or if $g=1$ and $X$ has a rational point), then these data
determine the curve up to a finite number of possibilities. The second conjecture concerns the
case when $K=\bfQ$ and $S$ is the empty set: it
\renewcommand{\thefootnote}{}\footnote{\kern-1em1980 \textit{Mathematics Subject Classification
\r(\r{1985} Revision\r)}. Primary 11G35, 11G10, 14H25, 14K15, 11D41, 14D10; Secondary 11G25,
14K10, 14L15, 14G25.\par Translation of the Appendix to the Russian translation of Serge Lang,
\textit{Fundamentals of Diophantine geometry} (original, Springer-Verlag, 1983---hereinafter
referred to as ``Lang's book''), ``Mir'', Moscow, 1986, pp.\,369--438; MR \textbf{88a}: 11054.}
\renewcommand\thefootnote{(\arabic{footnote})}\addtocounter{footnote}{-1}
\oripage states that there are no curves of genus $g\geq1$ with those invariants. These
conjectures are analogs of two classical results in algebraic number theory. The
first---Hermite's theorem (1857)---establishes finiteness of the number of extensions $L/K$ of
$K$ having fixed degree and fixed ramification. The second---Minkowski's theorem (1891)---states
that there are no unramified extensions of $\bfQ$. Comparing these theorems with Shafarevich's
conjectures, we see that those conjectures are their analogs for field extensions of relative
dimension (transcendence degree) $1$, while the theorems themselves concern the case of finite
extensions (i.\,e., relative dimension $0$). Here the analogs of the places of ramification in
the conjectures are the places of bad reduction.\footnote{By analogy with the concepts of
different and discriminant, Shafarevich introduced corresponding concepts for curves $X$ over a
field $K$ (see~\cite{8}). Incidentally, those concepts still have not been studied much.}

Shafarevich proved the finiteness conjecture for hyperelliptic curves as a consequence of
Siegel's finiteness theorem for the integer points (see Chapter 8 of Lang's book). (Variants of
the proof were published later in~\cite{Par2}, \cite{34}, and~\cite{56}.) In his report,
Shafarevich also raised the question of whether the analogous conjecture holds in the case when
$K=k(B)$ is the field of functions on an algebraic curve $B$ (with algebraically closed or finite
field of constants $k$). In this situation one must exclude ``constant'' curves $X$ for which
there exist finite extensions $L/K$ and $k'/k$ and a curve $X'$ over $k'$ such that $X\otimes
L\cong X'\otimes k'$. An approach to proving this conjecture in the case $k=\bfC$ was suggested
by Parshin in 1968 \cite{Par1}. The set of algebraic curves under consideration may be regarded
as the set ${\gotM}_g{(K)}_S$ of rational points of the moduli space of curves of genus $g$ (this
is a quasiprojective variety of dimension $3g-3$ having a canonical compactification
${\gotM}_g\subset{\ol{\gotM}}_g$) which satisfy the following additional condition: in
${\ol{\gotM}}_g$ these points are $S$-\textit{integral} relative to the closed subvariety
${\ol{\gotM}}_g-{\gotM}_g$ (concerning this notion, see Chapter 8, \S1, of Lang's book for the
case of a projective curve). In~\cite{Par1}, in accordance with general ideas of Diophantine
geometry, the problem was divided into two independent parts: proving that the heights of points
in ${\gotM}_g{(K)}_S$ (in some projective imbedding ${\ol\gotM}_g$) are bounded, and describing
the structure of the set of points of bounded height (see Chapter 3 of Lang's book). Boundedness
of the height was proved in the general case; and in the case when $S$ is empty it was shown that
the points of bounded height form a finite set. The generalization of the second part of the
proof to the case of arbitrary $S$ required new ideas, namely, the theory of theta-functions.
This was done in 1971 by S.\,Yu.\,Arakelov~\cite{Ar1}.\footnote{It was Szpiro~\cite{68},
\cite{69} who carried over the constructions and results in these papers to the case of a field
of finite characteristic We note that, using the techniques developed there Raynaud obtained a
counterexample to the Kodaira vanishing theorem for the cohomology of negative fiber bundles
(\cite{69}, Expos\'e 4).}

At the same time, Parshin~\cite{Par1}, \cite{57} (see also \S6) showed that in both the number
field and function field cases the finiteness conjecture implies the conjecture, advanced by
Mordell in 1922, that there are finitely many \oripage rational points on a curve of genus $g>1$.
(For the statement of the function field analog of Mordell's conjecture, see the commentary in
Chapter 8; this conjecture was first proved by Manin in 1963 \cite{Man1}.) The report~\cite{57}
at the 1970 International Congress of Mathematicians in Nice sketched a possible plan of proof of
the finiteness conjecture in the number field case. First of all, a natural generalization was
proposed to the case of an abelian variety $A$ over a number field $K$ with fixed set $S$ of
places of bad reduction (see \S2 for the definition, and see below for the fundamental theorem;
in~\cite{57} there was also a bound on the degree of the polarization on $A$). Since abelian
varieties are groups, it follows that they have the equivalence relation ``existence of an
isogeny'' (an epimorphism with finite kernel). The set $S$ is preserved under isogeny (see \S2),
and thus the finiteness conjecture splits into two parts. The first is the claim that there are
only finitely many abelian varieties over $K$ which are isogenous to a fixed one (Conjecture T).
The second part---which, at first glance, seems independent of the first---is the claim that the
number of isogeny classes is finite. This second part of the conjecture is closely connected with
some conjectures on homomorphisms of abelian varieties over $K$ and semisimplicity of the Tate
module that were proposed by Tate at Woods Hole in 1964 (see~\cite{Ta4}, \cite{26}, and Theorem 3
below). In particular, these conjectures of Tate imply that an abelian variety $A/K$ is
determined up to isogeny by the semisimple representation $\rho$ of the Galois group
$G_K=\Gal(\ol K/K)$ of the algebraic closure $\ol K$ of $K$ in the Tate module $V_{\ell}(A)$ (see
below).

Thus, if we accept Tate's conjectures, we reduce the second part to the question of finiteness of
the number of characters $\chi_\rho$ of representations $\rho$ in
$V_{\ell}(A)$.
\footletter{\textit{Added in translation.} The point is that any semisimple
representation over a field of characteristic zero $(={\bfQ}_{\ell})$ is uniquely determined by
its character.} See~\cite{5}, \S107. The next step consists in considering the Frobenius
automorphisms $Fr_v$, $v\not\in S$, in the group $G_K$ (see \S6). The Riemann hypothesis for
abelian varieties over a finite field (proved by Weil) shows that for every $v$ the number of
values ${\chi}_\rho(Fr_v)$ is bounded by a constant which depends only on the field $K$, the
place $v$ and $\dim A$. According to the Chebotarev density theorem~\cite{21}, \cite{67}, the
character $\chi_\rho$ is determined by its values ${\chi}_\rho(Fr_v)$; hence, the problem reduces
to constructing a finite subset $Q$ in $\{Fr_v\}$, depending only on $S$ and $\dim A$, from
which the values of $\chi_\rho$ for all $Fr_v$ are determined.

But, of course, there still remains the question of proving Tate's conjectures. In 1966, Tate
\cite{27} showed that they hold if the ground field $K$ is finite. Here he noted that an
essential role in the proof is played by Conjecture $T$ above (for a finite field, when it is
trivial), and he suggested that this conjecture might be useful for other fields $K$ as well.

That this is valid was confirmed in~\cite{21}, \cite{57}, and~\cite{Par2} for abelian varieties
of dimension~$1$ over a number field.

\oripage

Thus, the fundamental finiteness conjecture for abelian varieties was reduced to a series of
different but closely connected problems. (See the diagram of the definitive proof, below.)
\cite{57} also contains a discussion of the question of deducing Shafarevich's original
finiteness conjecture from the fundamental finiteness conjecture. Although it seems extremely
natural to make such a reduction---one goes from the curve $X$ to its Jacobian variety (see the
proof of the theorem in \S6)---nevertheless, \cite{57} contained some unjustified statements
about the insufficiency of the Torelli theorem for this purpose.

How did further events develop? The question of the finiteness of the number of characters (or
construction of the set $Q$) turned out to be very difficult. As noted in~\cite{57}, this
question has a positive answer if one assumes (as yet unproved) conjectures about zeta-functions
of abelian varieties.\footnote{Incidentally, now, after Faltings finally resolved this question,
the functional equation and the Taniyama--Weil conjecture have turned out to be useful in
questions of effectiveness. See~\cite{71}, Expos\'e~IX.\par \textit{Editor's note}. With
reference to the designation ``Taniyama--Weil conjecture'', see Note 1 by Serge Lang appended at
the end of this paper. } In 1981, Serre published a proof of the finiteness theorem for
characters in the case of elliptic curves over $\bfQ$ \cite{67}, assuming the generalized Riemann
hypothesis.

The question of the connection between Conjecture $T$ and Tate's conjectures was resolved for
varieties of any dimension in 1975--1977 by Yu.\,G.\,Zarkhin~\cite{11}--\cite{13}, using a
``quaternion trick'' he introduced. It turned out that for \textit{any} field $K$ finiteness
Conjecture $T$ for isogenies implies Tate's conjectures on homomorphisms and semisimplicity of
the Tate module of abelian varieties over $K$. Thus, in this part the arithmetic of $K$ plays no
role at all.

The greatest effort was required for the proof of Conjecture $T$ itself. Work ing with abelian
varieties rather than curves had the advantage that one can study the corresponding set of
rational points on the moduli space ${\gotM}^{ab}_g$ of abelian varieties of dimension $g$ (just
as this was done above for curves). Also, the isogeny property for abelian varieties has a
remarkable geometric interpretation in the language of the moduli space ${\gotM}^{ab}_g$. Over
${\gotM}^{ab}_g$ there exists a tower of Hecke correspondences
\begin{equation}\label{eq1.1}
H_{g,m}\mathop{\rightrightarrows}\limits_{p_2}^{p_1}\,{\gotM}^{ab}_g,
\end{equation}
which are the graphs of the equivalence relation ``existence of an isogeny of degree $m$'' (see
\S2.2). The universal family $A$ of abelian varieties determines an invertible sheaf $\omega$ on
${\gotM}^{ab}_g$ (\;=\;one-dimensional fiber bundle equal to the determinant of the conormal
bundle to the identity section; see \S2.1). The sheaf $\omega$ has the invariance property
\begin{equation}\label{eq1.2}
p^*_1\omega\cong p^*_2\omega.
\end{equation}
\oripage We now consider the height $h_\omega$, (the Weil height, or class of heights) on
${\gotM}^{ab}_g(K)$ which is associated to the sheaf $\omega$ (see Chapters 3 and 4 of Lang's
book),\footnote{In the terminology of Chapter 4, \S5, this is the height associated with the
divisor class of the invertible sheaf $\omega$.} and we try to apply its functorial properties
(Chapter 4, Theorem 5.1) to the diagram \eqref{eq1.1}, taking into account the isomorphism
\eqref{eq1.2}. For $m={\ell}^i$ and for points $P$,$Q\in{\gotM}^{ab}_g(K)$, for which ${\calA}_P$
and ${\calA}_Q$ are connected by an isogeny of degree ${\ell}^i$, we have
\begin{equation}\label{eq1.3}
h_\omega(P)\leq h_\omega(Q)+i\cdot c,\quad c=\const.
\end{equation}

In the early 1960's, Tate proposed modifying the definition of the Weil height $h_L(P)$, $P\in
A(K)$, $L\in\Pic A$, on an abelian variety $A$ in such a way that $h_L(P)$ is a canonical
function of $P$ (and not a class of functions, as in Chapter 3) and the following
\textit{equality} holds if $f\colon A\rightarrow B$ is a homomorphism:
\begin{equation}\label{eq1.4}
{\hat h}_L(f(P))={\hat h}_{f^*L}(P),
\end{equation}
(rather than the \textit{bound} in Theorem 5.1). (This is the well-known Tate height; see Chapter
5 of Lang's book and the appendix to~\cite{Mum2}.) If property \eqref{eq1.4} were to hold for the
morphisms in \eqref{eq1.1} and the sheaf $\omega$, then in place of the bound \eqref{eq1.3} we
would obtain the equality
\begin{equation}\label{eq1.5}
{\hat h}_\omega(P)={\hat h}_\omega(Q)
\end{equation}
for an isogeny of any degree $m$. Since $\omega$ is an ample sheaf on ${\gotM}^{ab}_g$, it
follows that Conjecture $T$ would give the finiteness property of the height (Chapter 3, Theorem
2.6). This program and a method for constructing the height ${\hat h}_\omega$ on ${\gotM}^{ab}_g$
were proposed by A.\,N.\,Parshin in his report at the International Conference on Number Theory
in Moscow (1971) \cite{17}.

Here the point of departure was again the case of an abelian variety $A$ over a function field
$K=k(B)$. If $f\colon A\rightarrow B$ is the N\'eron model of the abelian variety $A$ over $K$
(see \S2.1), then we set
\begin{equation}\label{eq1.6}
d(A)=c_1(\omega_{A/B}).
\end{equation}
This invariant arose in~\cite{Par1}, where it was noted that, when $S(A)$ is the empty set, then
the number of families $A/B$ with given $d(A)$ depends only on a finite number of parameters
(over $k=\bfC$). This fact was then proved in~\cite{Ar1} for families with arbitrary $S$ and
$k=\bfC$, and in~\cite{18} for families of Jacobian varieties and $k$ of arbitrary
characteristic. Since $d(A)$ is obviously preserved under isogenies (of degree prime to char
$k$), it follows that the function $d(A)$ is completely suitable for the applications mentioned
above. In order to construct a number field analog, it is natural to consider the expression
\eqref{eq1.6} and add on archimedean components corresponding to places of the field $K$. The
representation of the canonical height ${\hat h}_\omega$ as a product was suggested by the
fundamental work of N\'eron, which introduced the de n of the Tate height into local components
(see Chapter 11 of \oripage Lang's book).${}^*$
\renewcommand{\thefootnote}{}\footnote{\kern-1em${}^*$\textit{Editor's note}. Concerning N\'eron's contribution, see Note 2 by Serge Lang appended at
the end of this paper.}\renewcommand\thefootnote{(\arabic{footnote})}\addtocounter{footnote}{-1}
In this decomposition the archimedean components have the form ${\hat
h}_{L,\infty}=\log|\theta_D(z)|$, if $L={\goto}_A(D)$ and $\theta_D$ is the theta-function
corresponding to the divisor $D$ (see [Ne3], Chapter 3, and Chapter 13 of Lang's book).
In~\cite{17} it was noted that the following Poisson equation holds for the archimedean
components:
\begin{equation}\label{eq1.7}
\Delta_A({\hat h}_{L,\infty})=\const.
\end{equation}
In~\cite{17}, Equation \eqref{eq1.7} was taken as the definition of the canonical archimedean
components on an arbitrary variety $X/K$ having fixed K\"ahler metrics on $X\otimes\bfC$ for all
imbeddings $K\subset\bfC$. In the case of ${\gotM}^{ab}_g$ one can give an explicit expression
for ${\hat h}_{\omega,\infty}$ (see the functions $h_{\mod,v}$ in \S2.3).

A finiteness theorem was obtained for the height ${\hat h}_{\omega}$, (for abelian varieties with
potentially good reduction; an exposition is given in \S2.3). However, mistakes in the
computations of the behavior of the height under isogenies prevented one from finding the
situations when \eqref{eq1.5} is fulfilled.

A completely different method for constructing the canonical height is order to
prove Conjecture T was proposed by Zarkhin in 1974 (\cite{9}, \cite{10}). It is
well known that Mumford's theory of symmetric theta-structures~\cite{Mum3},
\cite{Z-M}, \cite{14} enables one to construct a canonical basis for the space
of sections $\Gamma(A,L)$ for a completely symmetric abelian invertible sheaf
$L$ on an abelian variety $A$. One can then take the Weil height of the zero
point corresponding to this basis (Chapter 3) as the canonical height of the
variety $A$ over $K$. Here the finiteness property holds by definition. But the
invariance of the finiteness of this height under isogeny is a much more
difficult question. Zarkhin resolved this question in the case of a function
field of arbitrary transcendence degree over ${\bfF}_p$ ($p\neq2$). By the same
token, in this situation Conjecture T was obtained, and also, in view of the
results mentioned above, Tate's homomorphism and semisimplicity conjectures. In
the number field case the possibility of applying this height remained an open
question.\footnote{Concerning the connection between these two heights on the
moduli space of abelian varieties, see the book~\cite{52} by Moret-Bailly. This
book gives an exposition of a result of S.\,Mori, who was able to remove the
condition $p\neq2$ in the theorems of Zarkhin mentioned above.}

After an interruption of several years, this direction of research was continued by the German
mathematician Gerd Faltings from Wuppertal. In June, 1983, his preprint \textit{Einige S\"atze
zum Thema abelsche Variet\"aten \"uber Zahlk\"orpern} appeared, and then his paper~\cite{39},
which caused a great commotion throughout the world. These papers contained proofs of the
conjectures of Tate, Shafarevich, and Mordell, following the program described above. The author
was able to obtain unexpected and very clever solutions to all the remaining unsolved parts of
the program.

First of all, Faltings gave an amazingly simple proof of the finiteness the orem for characters
of semisimple representations (the existence of the finite \oripage set $Q$); this proof takes up
less than a page in his paper (see the exposition in \S6). The catalyst for this result was a
theorem of Deligne on finiteness of the number of characters of complex representations of
discrete groups (see below).

Thus, the center of gravity of the whole problem immediately shifted to proving Conjecture T.
Faltings introduced the canonical height of an abelian variety (see \S2.1) and proved a
finiteness theorem for it (Theorem 1 in~\cite{39}, p.\,356). Faltings' definition of height was
influenced by the work of Arakelov, which we shall discuss later. The proof in~\cite{39} of the
theorem on finiteness of the number of abelian varieties with given height contains gaps (in the
proof of Lemma 2 on p.\,352), and goes through completely in the case of varieties with
potentially good reduction (this gives a proof of all of the conjectures when this condition is
fulfilled, in particular, the Mordell conjecture for such curves, for example, for the Fermat
curve). In order to treat points of bad reduction new arguments were needed. This was done by
Deligne, Gabber, and Moret-Bailly in~\cite{71}, Exposes IV and V (see also
\cite{40}).
\footletter{\textit{Added in translation.} The compactification in~\cite{40} of
the moduli scheme of abelian varieties over $\bfZ$ enables one to obtain
another version of the proof. The book~\cite{6*} constructs a compactification
of the moduli scheme of abelian varieties over $\bfZ[1/2]$ using Mumford's
theory of theta-structures~\cite{Mum3}.}

The last unexpected result in~\cite{39} concerned the behavior of the height under isogeny.
Although the height is not preserved, in general, nevertheless there exist two important classes
of isogenies for which invariance of height holds (see Claims 1 and 2 below). Faltings was able
to show that this is sufficient to prove Conjecture T. An important role in this part of the
paper is the combined application of Raynaud's theorem~\cite{60} on the determinant of group
schemes of period $p$, the Riemann hypothesis (Weil's theorem), and the corollary of Minkowski's
theorem which states that the Galois group of any extension of $\bfQ$ is generated by the inertia
subgroups of the ramified places (see \S5.1). In the special case of elliptic curves over $\bfQ$,
this argument\footnote{In many cases it enables one to find rational subgroups of the group of
points of an abelian variety (there are numerous examples in~\cite{Ser3}). A general method in
the case of elliptic curves was recently proposed in~\cite{36}.} occurs (in another connection)
in Serre's 1972 paper~\cite{Ser3} (pp.\,307 and 313).

Raynaud's formula for the determinant (see \S4.3.4) arose in answer to a
question of Serre relating to~\cite{Ser3}. Thus, \cite{39} showed the need to
make essential modifications in the rather straightforward behavior of the
height defined by \eqref{eq1.5}.\footnote{We give the central result
of~\cite{39} (see the fundamental theorem below) in a strong form (without the
assumption concerning the polarization) due to Zarkhin~\cite{74}. In Faltings'
paper~\cite{Fa2}, which appeared at the same time as~\cite{39}, the same
problem was examined in the case of a function field (for abelian varieties of
dimension $\leq3$ and $k=\bfC$ it had been resolved by Arakelov in~\cite{Ar1}).
The situation here turned out to be more complicated than in the case described
above of curves in the first Shafarevich conjecture (boundedness of the height
is preserved, but finiteness is not; see also~\cite{52}).} \oripage

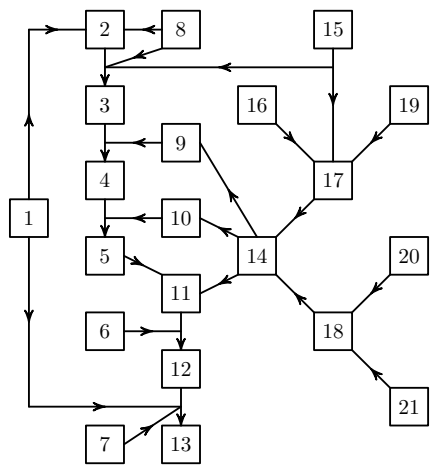
\bigskip

\centerline{\includegraphics{pic.1}}
\bigskip

\centerline{\bf Plan of Proof}

\medskip

{\footnotesize 1.~Hermite's theorem. 2.~Finiteness of the set $Q$. 3.~Finiteness of the number of
characters. 4.~Finiteness of the number of representations in $V_{\ell}(A)$. 5.~Finiteness of the
number of isogeny classes. 6.~The Torelli theorem. 7.~Construction of unramified coverings.
8.~The Chebotarev density theorem. 9.~Semisimplicity of the representations in $V_{\ell}(A)$.
10.~The Tate homomorphism conjecture. 11.~The finiteness conjecture for abelian varieties.
12.~Shafarevich's finiteness conjecture for curves. 13.~The Mordell conjecture. 14.~The
finiteness theorem for isogenies (Conjecture~T). 15.~The Riemann--Weil hypothesis for abelian
varieties. 16.~Raynaud's formula. 17.~The behavior of the canonical height under isogenies. 18.
The finiteness theorem for the canonical height. 19.~The arithmetic monodromy principle. (The
Galois group of any extension of $\bfQ$ is generated by the inertia groups of its ramified
places.) 20.~Projective imbedding of the moduli space. 21.~Gabber's lemma (Let $M$ be a dense
open subset of a normal integral excellent scheme $\ol M$, and let $f\colon A\rightarrow M$ be a
group object in the category of algebraic spaces proper over $M$ and whose fibers are abelian
varieties. There exists a proper surjective morphism $\pi\colon M'\rightarrow \ol M$ such that
the preimage of $A$ over $\pi^{-1}(M)$ extends to an algebraic space group object which is flat
over $M'$ and whose fibers are extensions of abelian varieties by tori.)}

\medskip

As we already noted, some unpublished work of Deligne in the early 1980's
played a role in the proof of Faltings' theorem on characters of
representations. In view of their importance, we shall give a statement of
Deligne's results.\footnote{They were communicated to one of the authors in a
letter from Deligne dated 23 November 1983 (see~\cite{8*}); they are published
here with his kind permission. For details, see the survey being prepared by
F.\,A.\,Bogomolov and A. N.\,Parshin entitled \textit{Differential-geometric
methods in algebraic geometry} (in the series ``Contemporary Problems of
Mathematics'', published by VINITI).
Added in December 2009: So far, this survey has not appeared.
} We first recall that to every family of Hodge structures of dimension $n$
over a variety $U$ there corresponds a local system of vector spaces of
dimension $n$ over $U$, or, equivalently, a representation $\rho\colon
{\pi}_1(U)\rightarrow\GL(n,\bfC)$.

\oripage

\begin{xxxi}[Theorem (Deligne)]Suppose we are given a nonsingular algebraic variety $U$ over
$\bfC$ and an integer $n$. There exist only finitely many representations $\rho\colon
\pi_1(U)\rightarrow\GL(n,\bfC)$, arising from families of polarized Hodge structures of dimension
$n$ over $U$.
\end{xxxi}

\begin{xxxr}[Sketch of the proof] The semisimplicity of the representations $\rho$ was proved earlier by
Deligne in~\cite{7}. Thus, it suffices to show that there are only finitely many possible
characters $\chi_\rho$ (see~\cite{5}, \S107). This follows from the following two independent
claims.
\end{xxxr}

\begin{enumerate}[\upshape1)]
\item \textit{There exists a constant $C_\gamma$ such that $|\chi_\rho(\gamma)|\leq C_\gamma$ for all
$\rho$.}
\item \textit{There exists a finite subset $Q\subset\pi_1(U)$ such that, for all $\rho_1$ and $\rho_2$, the equality
$\chi_{\rho_1}(\gamma)=\chi_{\rho_2}(\gamma)\ \forall\gamma\in Q$ implies
$\chi_{\rho_1}=\chi_{\rho_2}$.}
\end{enumerate}

The first estimate is obtained by mapping $U$ to the moduli space of Hodge structures and using
the hyperbolic Kobayashi metric (compare with the use of the Riemann hypothesis above). The
second result holds for all representations of fixed dimension $n$ of a given (discrete) group
$G$ ($=\pi_1(U)$) with a finite number $m$ of generators. It suffices to verify it for a free
group $G$. In that case it follows from the following result from invariant theory.

\begin{xxxi}[Theorem \cite{25}, \cite{58}]\footnote{See also the discussion in~\cite{53}, p.\,164. We are grateful to A.\,E.
Zalesskil for an informative letter on this question.} Let $X=\GL{(n,k)}^m$ be an affine
algebraic variety over a closed field $k$, on which the group $\GL(n)$ acts by (simultaneous)
conjugation. Then the ring of invariants $k{[X]}^{\GL(n)}$ is generated by the functions
\begin{equation}\label{eq1.8}
\Tr(g_{i_1}\dots g_{i_s}),
\end{equation}
where $g_i\in\GL(n,k)$ and $i_1,\dots,i_s$ run through all possible words of $m$ letters.
\end{xxxi}

According to Hilbert, the ring of invariants contains a finite number of generators of the form
\eqref{eq1.8}. Let $\gamma_1,\dots,\gamma_m$ be generators of the group $G$, and let $Q$ be the
set of elements of $G$ of the form $\gamma_{i_1}\dots\gamma_{i_s}$ where $i_1,\dots,i_s$ is run
through the words in $G$ which by \eqref{eq1.8} correspond to generators of the ring
$k{[X]}^{\GL(n)}$. The character $\chi_\rho(\gamma)$ as a function of $\rho$ determines a regular
invariant function $f_\gamma$ on $X$:
$$
f_\gamma(g_1,\dots,g_m)=\Tr(g_{i_1}\dots g_{i_s}),
$$
where $i_1\dots i_s$ is is the word giving the element $\gamma$ of the free group $G$. If we
represent $f_\gamma$ as a product of generators of the ring of invariants, we obtain the required
result.

We now return to the late 1960's and early 1970's. Parshin's paper~\cite{Par1}
also contained a proof of the Mordell conjecture in the function field case
carried out entirely within the framework of the theory of algebraic surfaces
(as presented, for example, in~\cite{32}), without bringing in the Shafarevich
conjecture, let alone the Tate conjectures. The basic ingredients were
intersection theory, the adjunction formula, the Riemann--Roch formula, and the
Hodge decomposition of cohomology. The proof of Mordell's conjecture \oripage
given by Grauert~\cite{Gra} was also to a large extent algebraic. This
stimulated attempts to carry over the concepts and constructions from the
geometric situation to the number field case. To every curve $X$ over a global
field $K$ one can associate its minimal model - two-dimensional scheme fibered
over the ring of integers of $K$ in the number field case (see \S6) or over a
curve in the geometric case (see~\cite{32} and~\cite{Sh3}). However, there is a
difference in principle between these two constructions: in the geometric case
the surfaces which arise are complete (projective), while in the number field
case one must add ``by hand'' the fibers over the archimedean places of $K$
(see the situation with the product formula in~\cite{4}, \cite{L5}, and Chapter
2 of Lang's book). This means that, in order to carry over some technique of
algebraic geometry to the number field situation, one must accomplish two
steps. First, one must consider arithmetic surfaces (in the framework of the
theory of schemes); then one must determine the ``archimedean analog'' of the
concept being considered. This program was realized by Shafarevich for
intersection theory and the theory of canonical classes (but without
archimedean components) in his Bombay lectures~\cite{Sh3} (see also~\cite{38},
Expos\'e X). At the end of those lectures, he posed the problem, still
unsolved, of carrying over the notion of tangent bundle to the number field
situation. We note that the ideas in this program go back to the deep analogy
between algebraic number fields and algebraic function fields that was first
noted in the nineteenth century by Kronecker and Hilbert. This analogy has
since been discussed and employed many times (see~\cite{31} and~\cite{6}).

Arakelov's papers, \cite{Ar2} and~\cite{Ar3}, heralded remarkable progress in this direction.
Starting from the definition of canonical height in~\cite{17} (see equation \eqref{eq1.7} above),
and adding to it the normalization condition
$$
\int_X\hat h_{L,\infty}=0,
$$
on arithmetic surfaces he defined such concepts as a divisor, the divisor of a function, linear
equivalence, the intersection index, and the canonical class $c_1$ and he proved an adjunction
formula. In~\cite{Ar3} he also formulated an analog of the Riemann--Roch theorem (unfortunately,
part of Arakelov's results remained unpublished, and as a consequence were re-proved
in~\cite{Fa1} and~\cite{46}). However, it was not then possible to construct analogs of such
concepts as the Hodge decomposition or the tangent bundle.

\begin{xxxr}[Remark added in translation] We note that even in the part of the theory of arithmetic surfaces
that has already been constructed it has not been possible to obtain a complete analogy with the
usual theory of algebraic surfaces. For example, the adjunction formula is true only for sections
(as noted by Szpiro). There is no projection formula for the intersection index under finite
morphisms. The entire theory is rather artificially divided into two parts: arithmetic surfaces
with generic fiber of genus $g\geq1$ and ``ruled'' arithmetic surfaces. In the theory we have,
the latter surfaces appear to be \oripage very trivial objects. Meanwhile, an analysis of
families of curves of genus 0 shows that the structure of the archimedean fibers $X_v/\bfC$
(\;=\;compact Riemann surfaces) is much more complicated. The tendency to regard the pair ($X_v$,
K\"ahler metric) as an analog of the fiber model $X_v/K_v$ (where $K_v$ is a nonarchimedean local
field) over Spec $O_v$ (where $O_v$ is its ring of integers) is clearly insufficient. From the
point of view of the analogy between nonarchimedean local fields and the field $\bfC$, the set of
models of the curve $X=P^1/\bfC$ can be described (or rather defined!) as a certain set of
probability measures in three-dimensional hyperbolic space for which the Riemann surface $X$ is
the absolute space. In order to pass to the case of curves $X$ of genus $g>0$ one must explicitly
examine the Schottky uniformization of the curves $X$ and the two-dimensional analog of the
Fourier--Tate transformation in local fields.\end{xxxr}

In this connection, Arakelov proposed another path to a proof of the Mordell conjecture in the
framework of the theory of arithmetic surfaces, using the $\zeta$-functions introduced in
\cite{Ar3} which are associated to divisors on a surface. Unfortunately, this program also
remained unrealized.

Since that time, in the theory of arithmetic surfaces one further fundamental result was obtained
which has a purely algebraic formulation---the Bogomolov--Miyaoka--Yao inequality (see~\cite{3}
and \cite{51}) for surfaces $V$ of general type over a field $k$ of characteristic 0:
\begin{equation}\label{eq1.9}
c_1{(V)}^2\leq3c_2(V),
\end{equation}
where $c_1(V)$ is the canonical class and $c_2(V)$ is the Euler characteristic. If one could
obtain inequality \eqref{eq1.9} for an arithmetic surface, then the constructions of Chapter III
of \cite{Par1} combined with Arakelov's intersection theory would enable one to give a new
\textit{effective} proof of the Mordell conjecture. A difficulty with this approach is that
\eqref{eq1.9} is false for fields $k$ of finite characteristic (see, for example, \cite{Par2}).
See also the Appendix to \S6, below.

\begin{xxxr}[Remark added in translation] The arithmetic analog of the Bogomolov--Miyaoki--Yao inequality for an
algebraic curve $X$ of genus $g\geq1$ defined over an algebraic number field
$K$ looks as follows :
$$
\omega\cdot\omega\leq3\delta+(2g-2)\log|D_{K/\bfQ}|.
$$
Here $\omega$ is the relative canonical class of the minimal model of $X$,
$\omega\cdot\omega$ is the Arakelov intersection index, $\delta$ is the
invariant introduced by Faltings~\cite{Fa1}, and $D_{K/\bfQ}$ is the
discriminant of the field $K$ over $\bfQ$. Here it is assumed that the curve
has stable reduction over $K$, and the metrics on the archimedean components of
the corresponding arithmetic surface are Arakelov metrics (i.\,e., they are
induced by a flat metric on the Jacobian variety). See~\cite{4*,5**,6**}. To
obtain the result mentioned above, and also an application to Fermat's theorem
(see p.\,435), it is enough to have the inequality with any absolute constant
$A$ in place of 3 in front of $\delta$. We further note that van de
Ven~\cite{5*} was the first to obtain an inequality of the type $c^2_1\leq
Ac_2$ in the geometric case. His proof uses Hodge theory for surfaces in an
essential way, and it gives $A=8$.
\end{xxxr}

\oripage

We note that Arakelov's papers, \cite{Ar2} and~\cite{Ar3}, contained one other approach to
defining the concepts of algebraic geometry mentioned above, using the technique of metrized
sheaves. This approach aroused great interest after the appearance of Faltings' papers~\cite{Fa1}
and~\cite{39}, and it led to a whole series of papers (\cite{71}, Expos\'es I--III, \cite{47},
\cite{16a}, \cite{45}, and~\cite{15b}).

\begin{xxxr}[Remark added in translation] It should be noted that the survey~\cite{16a} contains the following
inaccuracies. The comparison between the intersection indices in the archimedean and
nonarchimedean situations at the end of \S1.5 is wrong. The Green's function $G_v(D_v,x)$,
$K_v=\bfC$, $x\in(X-D_v)(\bfC)$, $D_v$ a divisor on the Riemann surface $X_v$ has as its analog
in the nonarchimedean situation the intersection index ${\langle D_v,x\rangle}_v$ on the regular
model $Y/O_v$ of the curve $X_v/K_v$ (where $K_v$ is a nonarchimedean local field and $O_v$ is
its ring of integers). The condition $x\not\in D_v$ means that $x$ and $D_v$ correspond to
subschemes of $Y$ which intersect only over the closed point of $\Spec O_v$. Even if the model
$Y$ is smooth over $\Spec O_v$, the function $f(x)=\langle D_v,x\rangle$, $x\in(X-D_v)(K_v)$, is
not constant, nor even locally constant. It is also clear that the presence of degenerations does
not (despite the claim in~\cite{16a}) have any influence on such properties of the function
$f(x)$. This function always has an infinite set of values.
\end{xxxr}

[The remark in \S3.1 concerning the special role of the K\"ahler--Einstein metric for arithmetic
varieties contradicts the concrete results in \S\S1 and 2---the adjunction formula and the
Riemann--Roch theorem, which are valid only for the Arakelov metric, which is not a
K\"ahler--Einstein metric. This was pointed out to one of the authors by M.\,Atiyah.]

At the beginning of this Introduction we alluded to Shafarevich's second conjecture in~\cite{30}.
It can be reformulated as the statement that there do not exist smooth proper schemes over
$\Spec\bfZ$ of relative dimension 1 and genus $\geq1$. The geometric analog of this conjecture
was thoroughly investigated in~\cite{Par1} and~\cite{69}. The generalization of this problem to
abelian varieties was examined in a series of papers by V.\,A.\,Abrashkin in 1976--1977 (see
\cite{20}). He found a general approach to this conjecture using finite group schemes, the
Riemann--Weil hypothesis, and bounds for the discriminants of algebraic number fields. This
problem was then solved for abelian varieties of dimension $\leq3$. Very recently, Fontaine
\cite{43} and independently Abrashkin~\cite{1}, \cite{1a} obtained a complete solution along
these lines.

\begin{xxxi}[Theorem \cite{43}, \cite{1}] There do not exist smooth abelian schemes over the rings of integers of the
algebraic number fields $\bfQ$, $\bfQ(\sqrt{-1})$, $\bfQ(\sqrt{-2})$, $\bfQ(\sqrt{-3})$,
$\bfQ(\sqrt{-7})$, $\bfQ(\sqrt{2})$, $\bfQ(\sqrt{5})$, $\bfQ(\sqrt[5]{1})$.
\end{xxxi}

We now proceed to describe the plan of proof of Faltings' results that will be used in this
survey.

Let $S$ be a finite set of nonarchimedean places of the field $K$, and let $g$ be a natural
number. We let $\III(K,g,S)$ denote the set of $g$-dimensional abelian varieties over $K$,
considered up to $K$-isomorphism, and having good \oripage reduction outside $S$. If
$A\in\III(K,g,S)$ and if $B$ is an abelian variety which is $K$-isogenous to $A$, then $B$ also
belongs to the set $\III(K,g,S)$. Let $S'$ be the set of all places of the field $L$ which lie
over the places of $K$ which belong to $S$. The set $S'$ is finite, and if $A\in\III(K,g,S)$,
then $A_L\in\III(L,g,S')$.

\begin{xxxi}[Fundamental theorem] The set $\III(K,g,S)$ is finite.
\end{xxxi}

\begin{xxxr}[Remark] To prove the fundamental theorem it suffices to show that the set $\III(L,g,S')$ is finite
(for any extension $L$ of $K$), because we have the following fact.
\end{xxxr}

\begin{xxxi}[Finiteness theorem for forms of an abelian variety]
Let $X$ be an abelian variety over a field $F$, and let $E$ be a finite
separable extension of $F$. The set of abelian varieties $Y$ over $F$ such that
$Y{\otimes}_FE\cong X{\otimes}_FE$ is finite \r(up to $F$-isomorphism\r).
\end{xxxi}

The fundamental theorem follows from the following two weaker assertions.

\begin{xxxi}[Finiteness theorem for isogeny classes] There exist a finite number of abelian varieties
$A_{(1)},\dots,A_{(t)}$) over $K$ which satisfy the following condition. If $A$ is a
$g$-dimensional abelian variety over $K$ which has good reduction outside $S$, then $A$ is
isogenous over $K$ to one of the abelian varieties $A_{(i)}$.
\end{xxxi}

\begin{xxxi}[Finiteness theorem for isogenies] The set of abelian varieties $B$ over $K$ \r(considered up to
$K$-isomorphism\r) for which there exists a $K$-isogeny $A\rightarrow B$ is finite.
\end{xxxi}

We shall derive these theorems from Theorems 1--4 below.

\begin{xxxi}[Theorem 1] Let $K$ be a number field, and let $p$ be a prime which is unramified in $K$. Fix a
natural number $g$. There exists a finite set of primes $M=M(K,p,g)$ which can be explicitly
determined from $K$, $p$, $g$ and which satisfies the following conditions.

Let $A$ be a $g$-dimensional abelian variety over $K$ which has good reduction at all places of
$K$ lying above $p$. Let s $(A)$ be the \r(finite\r) set of prime numbers which are the
characteristics of the residue fields of the places in $S(A)$. Then the set of abelian varieties
$B$ over $K$ \r(considered up to $K$-isomorphism\r) for which there exists a $K$-isogeny
$A\rightarrow B$ of degree not divisible by any prime number in $M\cup s(A)$ is finite.
\end{xxxi}

\begin{xxxi}[Corollary] Let $A$ be an abelian variety over $K$. There exists a finite set $M$ of prime numbers such
that the set of abelian varieties $B$ over $K$ \r(considered up to $K$-isomorphism\r) for which
there is a $K$-isogeny $A\rightarrow B$ of degree not divisible by any prime in $M$ is finite.
\end{xxxi}

We shall prove Theorem 1 only for $K=\bfQ$ (\S5.0.3). When $K$ is arbitrary we shall limit
ourselves to a proof of the corollary (\S5.0.4).

\oripage

\begin{xxxi}[Theorem 2] Let $A$ be an abelian variety over a number field $K$, and let $\ell$ be a prime number. The
set of abelian varieties $B$ for which there exists a $K$-isogeny $A\rightarrow B$ whose degree
is a power of $\ell$ is finite \r(up to isomorphism\r).
\end{xxxi}

In order to state Theorems 3 and 4 we shall need the concept of the Tate module of an abelian
variety (see~\cite{Mum2} and~\cite{24}).

Let $F$ be an arbitrary field, let $\ol F$ be the algebraic closure of $F$, and let $G$ be the
Galois group of $\ol F$. Let $X$ be an abelian variety defined over $F$, and let End $X$ be its
ring of $F$-endomorphisms. Let $n$ be a natural number which is prime to the characteristic of
$F$. We let $X_n$ denote the subgroup of the group $X(\ol F)$ of $\ol F$-points consisting of the
elements annihilated by multiplication by $n$. It is well known~\cite{Mum2} that $X_n$ is a
finite abelian group which is a free $\bfZ/n\bfZ$-module of rank 2 $\dim X$. The natural actions
of the group $G$ and the ring End $X$ on $X(\ol F)$ preserve the subgroup $X_n$, thereby making
it into a finite Galois module and determining an imbedding
$$
\End X\otimes\bfZ/n\bfZ\hookrightarrow{\End}_G X_n.
$$
If $m$ and $n$ are natural numbers prime to the characteristic of $F$, then $X_n$ is a Galois
submodule of $X_{nm}$, equal to the kernel of multiplication by $n$, and we have the exact
sequence of Galois modules
$$
0\rightarrow X_n\rightarrow X_{nm}\stackrel n\rightarrow X_m\rightarrow0,
$$
which is compatible with the action of the ring $\End X$.

Now suppose that $n={\ell}^i$ is a power of a prime $\ell$ not equal to the characteristic of
$F$. We define the \textit{Tate $\ell$-module \r(or $\bfZ_{\ell}$-module\r)} $T_{\ell}(X)$ of the
abelian variety $X$ to be the projective limit (with respect to $i$)
$$
T_{\ell}(X)=\lim\limits_\leftarrow X_{{\ell}^i}
$$
(here the maps are multiplication by $\ell$). This limit $T_{\ell}(X)$ is a free
$\bfZ_{\ell}$-module of rank 2 $\dim x$. We set
$$
V_{\ell}(X)=T_{\ell}(X)\otimes_{{\bfQ}_{\ell}}{\bfZ}_{\ell}.
$$
Clearly $V_{\ell}(X)$ is a vector space over the field ${\bfQ}_{\ell}$ of $\ell$-adic numbers of
dimension 2 $\dim X$; we shall call it the \textit{Tate ${\bfQ}_{\ell}$-module} of the abelian
variety $X$. The actions of the group $G$ and the algebra $\End X$ on the $X_{{\ell}^i}$ glue
together to give a continuous homomorphism \textit{\r($\ell$-adic representation\r)}
$$
\rho_{\ell,X}\colon G\rightarrow\Aut T_{\ell}(X)\subset\Aut V_{\ell}(X)
$$
and imbeddings
$$
\begin{matrix}
\End X\otimes{\bfZ}_{\ell}&\hookrightarrow&{\End}_GT_{\ell}(X),\\
\cap& &\cap\\
\End X\otimes{\bfQ}_{\ell}&\hookrightarrow&{\End}_GV_{\ell}(X),\\
\end{matrix}
$$
respectively. We let $\chi_\rho=\chi_{\rho_{\ell,X}}\colon G\rightarrow{\bfZ}_{\ell}$ denote the
function which is the character of the representation $\rho_{\ell,X}$. If $Y$ is an abelian
variety over $F$ which is $F$-isogenous to $X$, then the representations $\rho_{\ell,X}$ and
$\rho_{\ell,Y}$, i.\,e., the modules $V_{\ell}(X)$ and $V_{\ell}(X)$, are isomorphic; in
particular, the characters $\chi_\rho$ are the same for $X$ and $Y$.

\oripage

\begin{xxxi}[Theorem 3]Let $A$ be an abelian variety over a number field $K$, and let $\ell$ be a prime number.
\begin{enumerate}[\upshape1)]
\item \textit{If $W$ is a ${\bfQ}_{\ell}$-subspace of $V_{\ell}(A)$ which is a Galois submodule,
then there exists an endomorphism $u\in\End A\otimes{\bfQ}_{\ell}$, such that $u^2=u$ and
$V_{\ell}(A)=W$, i.\,e., the subspace $(1-u)V_{\ell}(A)$ is a complement of $W$ in
$V_{\ell}(A)$. In particular, the Galois module $V_{\ell}(A)$ is semisimple.}

\item The Tate homomorphism conjecture. \textit{The natural imbeddings}
$$
\End A\otimes{\bfZ}_{\ell}\hookrightarrow{\End}_GT_{\ell}(A)\quad\text{and}\quad
\End A\otimes{\bfQ}_{\ell}\hookrightarrow{\End}_GV_{\ell}(A)
$$
\textit{are bijections.}
\item\textit{If $B$ is an abelian variety over $K$, and if the Galois modules $V_{\ell}(A)$ and $V_{\ell}(B)$ are
isomorphic, then the abelian varieties $A$ and $B$ are isomorphic over $K$.}
\end{enumerate}
\end{xxxi}

\begin{xxxr}[Remark] Parts 2 and 3 follow from part 1 applied to the abelian varieties $A^2$ and $A\times B$,
respectively.\end{xxxr}

\begin{xxxr}[Remark added in translation] If one replaces the Tate module by the group of points of
sufficiently large prime order $\ell$, then one has the following analog of Theorem 3 (see
\cite{74}).
\end{xxxr}

\textit{For all except finitely many primes $\ell$ the Galois module $A_{\ell}$ is
semisimple, and the natural imbedding $\End A\otimes\bfZ/\ell\bfZ\rightarrow{\End}_GA_{\ell}$
is a bijection.}

\begin{xxxi}
[Serre's theorem \cite{7*,9**,8**}] Set $g=\dim A$ and let $G(\ell)$ denote the
image of the Galois group $G$ in $\Aut A_{\ell}\approx\GL(2g,\bfZ/\ell\bfZ)$.
For all primes $\ell$ except perhaps finitely many, there exists a connected
reductive algebraic subgroup $H_{\ell}$ in $\GL_{2g}$ which is defined over
$\bfZ/\ell\bfZ$ and has the property that, if we replace $K$ by a fixed finite
algebraic extension not depending on $\ell$, the group
$H_{\ell}(\bfZ/\ell\bfZ)$ of $\bfZ/\ell\bfZ$-rational points contains $G(\ell)$
as a subgroup of  index bounded by a constant independent of $\ell$.
\end{xxxi}

\begin{xxxi}[Theorem 4] Suppose that $K$ is a number field, $G$ is the Galois group of its algebraic closure,
$S$ is a finite set of nonarchimedean places, $g$ is an integer, and $\ell$ is a prime number.
Then the set of characters $\chi_\rho\colon G\rightarrow{\bfZ}_{\ell}$ corresponding to abelian
varieties $A$ in $\III(K,g,S)$ is finite.
\end{xxxi}

The finiteness theorem for isogeny classes is obtained from Theorems 3 and 4 using the reasoning
outlined above.

\begin{xxxr}[Proof of the finiteness theorem for isogeny classes] (modulo the corollary to Theorem 1 and
Theorem 2). Let $A$ be an abelian variety over $K$, and let $M$ be a finite set of prime numbers
which satisfies the conclusion of the corollary to Theorem 1. Let $n$ be the number of elements
in the finite set $M$, and let ${\ell}_1,\dots{\ell}_n$ be all of the primes in $M$. Any isogeny
from $A$ to $B$ factors into a composition of isogenies
\end{xxxr}
$$
A\rightarrow B_{(0)}\rightarrow B_{(1)}\rightarrow\dots\rightarrow B_{(n)}=B,
$$
\oripage which satisfy the following conditions. The degree of the isogeny $A\rightarrow B_0$ is
not divisible by any of the primes in $M$. The degree of the isogeny $B_{(i-1)}\rightarrow
B_{(i)}$ is a power of the prime ${\ell}_i$.

According to Theorem 1, there exist only finitely many possibilities for the abelian variety
$B_{(0)}$. According to Theorem 2, applied to the abelian variety $B_{(0)}$ and the prime
${\ell}_1$, there exist only finitely many possibilities for $B_{(1)}$. A simple induction (based
on applying Theorem 2 to $B_{(i-1)}$ and ${\ell}_i$) shows that there are only finitely many
possibilities for the abelian variety $B_{(n)}=B$.

We now discuss the proof of Theorems 1--4. Concerning Theorem 4, see \S6. The proofs of Theorems
1, 2, and 3 are much more complicated, and are based upon the notion of the \textit{canonical
height $h(A)$ of an abelian variety $A$} defined over a number field. This height is a certain
real number; for the precise definition, see \S2.1. We have the following

\begin{xxxi}[Finiteness theorem for the height (\S2.1)] Let $K$ be a number field, $C$ a real number, and $g$ an
integer. The set of abelian varieties $A$ over $K$ \r(considered up to $K$-isomorphism\r) for
which $\dim A=g$, $h(A)\leq C$, and $A$ is principally polarized is finite.
\end{xxxi}

This theorem and a special case of the next claim will be used in \S5 in order to prove the
corollary to Theorem 1.

\begin{xxxi}[Claim 1] If $A\rightarrow B$ is the isogeny in the hypothesis of Theorem 1, then $h(A)=h(B)$.
\end{xxxi}

In \S5 we prove this special case (Corollary 5.0.2), and then we use a Galois descent along with
the finiteness theorem for the height in order to obtain the corollary to Theorem 1.

\begin{xxxi}[Claim 2] Let $A$ be an abelian variety over a number field $K$. Suppose that we are given an infinite
sequence of $K$-isogenies of abelian varieties
$$
A\rightarrow B_{(1)}\rightarrow B_{(2)}\rightarrow\dots\rightarrow B_{(n)}\rightarrow\dots,
$$
such that the kernels $W_n=\Ker(A\rightarrow B_{(n)})$ form an $\ell$-divisible group \r(see
\S3\r). Then there exists a natural number $N$ such that $h(B_{(n)})=h(B_{(+1)})$, for all $n>N$,
i.\,e., the height does not depend on $n$.
\end{xxxi}

The special case of Claim 2 when $K=\bfQ$ (Lemma 5.3.1) and the finiteness theorem for the height
are used to prove the following fact.

\begin{xxxi}[Claim 3] (Corollary 5.4). Let $A$ be an abelian variety over a number field $K$. Suppose that we are
given an infinite sequence of $K$-isogenies of abelian varieties $A\rightarrow B_{(1)}\rightarrow
B_{(2)}\rightarrow\dots$ such that the kernels $\Ker(A\rightarrow B_{(n)})$ form an
$\ell$-divisible group. Then the sequence $B_{(1)}$, $B_{(2)}$ , $\dots$ contains only finitely
many pairwise nonisomorphic abelian varieties.
\end{xxxi}

Now Theorems 2 and 3 are derived from Claim 3. We note that this derivation is purely algebraic:
it in no way uses the arithmetic of the field $K$ (see \S3).

\oripage

The most complete exposition of Fallings' results (and Arakelov's theory of arithmetic surfaces)
is in~\cite{71}. For a more condensed presentation, see the Bourbaki seminar reports~\cite{37}
and~\cite{70}. In~\cite{41} one can find a very detailed exposition of the part connected with
isogenies, finite group schemes, and $\ell$-adic representations.\footletter{Added in
translation. After this was written, the book~\cite{1*} also appeared.}

\section{The canonical height of an abelian variety and the moduli space}

\ssect Let $K$ be an algebraic number field of finite degree $[K:\bfQ]$ over $\bfQ$, let
$\Lambda_K$ be its ring of integers, let $\bfS=\Spec\Lambda_k$, and let $P(K)$ be the set of
places of $K$. Then $P(K)=P{(K)}_f\cup P{(K)}_\infty$, where $P{(K)}_f$ is the infinite set of
nonarchimedean (finite) places (they can be identified with the closed points of the scheme $S$
or with the discrete valuations of the field $K$), and $P{(K)}_\infty$ is the finite set of
archimedean (infinite) places $v$ (they can be identified with imbeddings $v\colon
K\hookrightarrow\bfC)$; see~\cite{4} and~\cite{25}. It is well known that
$\#P{(K)}_\infty=[K:\bfQ]$.

If $v\in P{(K)}_f$, then we let ${\goto}_v$ denote the corresponding complete local ring (the
completion of the ring $\Lambda_K$ with respect to the discrete valuation corresponding to $v$).
We shall use the same letter $v$ to denote the place $v$, the corresponding discrete valuation
$v\colon K^*\rightarrow\bfZ$, and also the extension of the valuation to the field of fractions
$K_v$ of the ring ${\goto}_v$. For a prime element $\pi\in{\gotm}_v\subset{\goto}_v$ (${\gotm}_v$
is the maximal ideal of ${\goto}_v$) we have $v(\pi)=1$, and if $k(v)={\goto}_v/{\gotm}_v$ is the
residue field, then $\bfN v=\#k(v)$.

We consider an abelian variety $A$ defined over a local field $K_v$ (a finite extension of
${\bfQ}_p$). It corresponds to a smooth group scheme $\varphi\colon
\bfA\rightarrow\Spec{\goto}_v$ over the ring ${\goto}_v$, called the \textit{N\'eron
\r(minimal\r) model} of the abelian variety $A$; it has the following property.

Given any smooth scheme $\bfB$ over $\Spec{\goto}_v$ and any morphism $f\colon B=\bfB\otimes
K_v\rightarrow A$, there exists a unique morphism $\tilde f\colon \bfB\rightarrow\bfA$, which
coincides with $f$ on the generic fiber.

The scheme $\bfA$ is uniquely determined by this condition; see~\cite{46}, \cite{59}, Chapter 11,
\S5, and Lang's book, Chapter 12, \S3. In what follows we shall let ${\bfA}^0$ denote the
connected component of the identity section in $\bfA$.

We let $e: \Spec{\goto}_v\rightarrow\bfA$ denote the identity section of the scheme $\bfA$. The
closed fiber ${\bfA}_v=\bfA\otimes k(v)$ is a commutative algebraic group over the field $k(v)$,
and one says that $A$ has at $v$ (or over $K_v$)

\textit{good reduction}, if ${\bfA}_v$ is an abelian variety;

\textit{potentially good reduction}, if there exists a finite extension $L/K_v$ such that the
abelian variety $A\otimes L$ has good reduction; and

\textit{semistable reduction}, if ${\bfA}^0_v$ is an extension of an abelian variety by a torus.

The last condition is equivalent to saying that ${\bfA}_v$ has no unipotent component. For any
abelian variety $A$ over $K_v$ there exists a finite extension $L/K_v$ such that $A\otimes L$ has
semistable reduction over $L$ (the \textit{semistable reduction theorem}) \cite{45}. \oripage If
for some $n\geq3$ prime to the characteristic of the residue field all of the points of order $n$
on $A$ are rational over $K$, then $A$ always has semistable reduction (Raynaud's criterion). We
further note that all three properties above are preserved in going from $A$ to any abelian
variety $B$ which is isogenous to it (\textit{invariance under isogeny}) \cite{24}. If $A$ has
semistable reduction over $K_v$, then for any extension $L/K_v$ the N\'eron model of the abelian
variety $A\otimes L$ is $\bfA{{\otimes}_{\goto}}_v{\goto}_L$ (${\goto}_L$ is the ring of integers of
$L$).

Now suppose that $\varphi\colon \bfA\rightarrow\bfS$, $\bfS=\Spec{\goto}_v$ is an arbitrary
smooth group scheme over ${\goto}_v$, and $A=\bfA\otimes K_v$ is the generic fiber. We let
$\bfE=e(\bfS)\subset\bfA$ denote the image of the identity section $e$. This is a one-dimensional
subscheme which is regularly imbedded in $\bfA$ and is mapped by $\varphi$ isomorphically onto
$\bfS$. The sheaf ${\Omega}^1_{\bfA/\bfS}$ of relative differentials of degree 1 when restricted
to $\bfE$ gives a locally free sheaf of rank $g=\dim A$ (the conormal bundle). It is isomorphic
to $e^*{\Omega}^1_{\bfA/\bfS}$. We set
$$
\omega_{\bfA/\bfS}=\Lambda^g_{{\goto}_v}(e^*{\Omega}^1_{\bfA/\bfS}).
$$
This is a free ${\goto}_v$-module of rank 1. If $\omega\in{\Omega}^g_{A/K_v}$ is a differential
form which is nonzero at the identity point of the group $A$, then $\omega$ determines a rational
section of the sheaf $\omega_{\bfA/\bfS}$, or, equivalently, an element of the vector space
$\omega_{\bfA/\bfS}\otimes K_v\cong{\Omega}^g_{A/K_v}$ of dimension 1 over $K_v$. The valuation
$v$ can be extended to ${\Omega}^g_{A/K_v}$, using the lattice $\omega_{\bfA/\bfS}$. If $G$ and
$H$ are commensurable subgroups (equal up to a subgroup of finite index) of a group, we set
$[G:H]=\#G/(G\cap H){(\#H/(G\cap H))}^{-1}$. In particular,
$$
[G:H]=\begin{cases}\#(G/H),\ &\text{if}\ H\subseteq G,\\\#{(H/G)}^{-1},\ &\text{if}\ G\subseteq H.\end{cases}
$$

Then
$$
\bfN v^{v(\omega)}=[\omega_{\bfA/\bfS}\colon {\goto}_v\cdot\omega].
$$

\begin{xxxr}[Definition 1] Let $A$ be an abelian variety over the field $K_v$, $g=\dim A$, $\omega\in{\Omega}^g_{A/K_v}$.
Suppose that the form $\omega\neq0$ is invariant under translations. Then
$$
h_{K_v}(A,\omega)=v(\omega)\log\bfN v
$$
is called the \textit{local factor} (at the place $v$) of the canonical height of the abelian
variety $A$ over $K_v$.
\end{xxxr}
This definition makes sense, since the form $\omega$ must be nonzero at the identity point.

\begin{xxxi}[Proposition 1] The local factors $h_{K_v}$ have the following properties.
\begin{enumerate}[\upshape1)]
\item If $\lambda\in K^*_v$, then
$$
h_{K_v}(A,\lambda\omega)=h_{K_v}(A,\omega)+v(\lambda)\log\bfN v.
$$
\item If the abelian variety $A$ has semistable reduction over $K_v$ and if $L/K_v$ is a finite
extension, then
$$
[L:K_v]\cdot h_{K_v}(A,\omega)=h_L(A\otimes L,\omega\otimes L).
$$
\item If the abelian variety $A$ has semistable reduction over $K_v$ and if $f\colon A\rightarrow
B$ is an isogeny, then
$$
h_{K_v}(A,f^*\omega)=h_{K_v}(B,\omega)+\log(\#e^*\Omega^1_{\bfG/\bfC}),
$$
where $\bfG=\Ker\tilde f$, $\tilde f\colon \bfA\rightarrow\bfB$ is the homomorphism of
N\'eron models determined by $f$, and $e\colon \bfS\rightarrow\bfG$ is the identity section.
\end{enumerate}
\end{xxxi}

Property 1) is obvious; property 2) follows from the fact that if $\bfA$ is the N\'eron model of
$A$ and ${\goto}_L$ is the ring of integers of $L$, then $\bfA\otimes_{{\goto}_v}{\goto}_L$ is
the N\'eron model of the abelian variety $A\otimes L$ over $L$. If the fiber $A_v$ has a
unipotent component, then this stability property of the N\'eron model is not, in general,
fulfilled (nor is property 2) of the proposition). In order to obtain 3), we note that the
sequence of homomorphisms of group schemes
$$
\bfG\stackrel i{\rightarrow}\bfA\stackrel{\tilde f}{\rightarrow}\bfB\quad\text{($i$ is a closed imbedding)}
$$
which exists by the definition of the N\'eron model, has the following additional property in the
case of semistable reduction of $A$ (and hence of $B$): $\tilde f$ is a flat morphism and $\bfG$
is a quasifinite flat group scheme over ${\goto}_v$ (see~\cite{71}, p.\,201, or~\cite{41}, p.
129).

$\bfG$ is a finite group scheme, i.\,e., proper over ${\goto}_v$, only if $A$ has good reduction
at $v$. If the reduction is not semistable, then the example of multiplication by $p$ (where $p$
is the characteristic of the residue field $k(v)$) shows that $\dim{\bfG}_v$ may be greater than
zero (multiplication by $p$ annihilates the unipotent component of ${\bfA}_v$), at the same time
as $\dim(\bfG\times K_v)=0$.

Under these conditions it can be shown that the sequence of maps of modules of differentials
induced by $(*)$
$$
0\rightarrow e^*\Omega^1_{\bfB/\bfS}\stackrel{\tilde f^*}{\rightarrow}e^*\Omega^1_{\bfA/\bfS}
\stackrel{i^*}{\rightarrow}e^*\Omega^1_{\bfG/\bfS}\rightarrow0
$$
is exact [loc.~cit.]. Here the first two groups are free modules over ${\goto}_v$ of rank $d=\dim
A=\dim B$ (see above), and the last group is finite. Passing to exterior powers, we obtain an
imbedding of free rank 1 ${\goto}_v$-modules:
$\omega_{\bfB/\bfS}\hookrightarrow\omega_{\bfA/\bfS}$, whose cokernel, which we shall denote
$\det\bfG$, is a finite ${\goto}_v$-module whose length $\calL(\det\bfG)$ is equal to the length
$\calL(e^*\Omega^1_{\bfG/\bfS})$ (the elementary divisor theorem). The orders of the groups
$\det\bfG$ and $e^*\Omega^1_{\bfG/\bfS}$ are also the same, and this proves property 3). We note
that there is exactly one simple ${\goto}_v$-module (up to isomorphism), namely, the residue
field $k$; this gives us the equality
$$
\#\det\bfG=\#e^*\Omega^1_{\bfG/\bfS}=\bfN v^{\calL(e^*\Omega^1_{\bfG/\bfS})}.
$$
Later (\S\S4 and 5) we shall denote the function $v(\omega)$ by $\calL(A,\omega)$ explicitly
indicating the dependence on the abelian variety $A$ (and we shall call it \textit{the local term
in the logarithmic height} of the abelian variety $A$ over $K_v$). We have
$$
h_{K_v}(A,\omega)=\calL(A,\omega)\log\bfN v.
$$
\oripage The local terms $\calL(A,\omega)$ have the following properties.
\begin{enumerate}[\upshape1)]
\item
If $\lambda\in K^*_v$, then
$$
\calL(A,\lambda\omega)=\calL(A,\omega)+v(\lambda).
$$
\item
Let $L/K_v$ be a finite extension. Then $\calL(A,\omega)=\calL(A\otimes L,\omega\otimes L)$.
If the abelian variety $A$ has \textit{semistable} reduction over $K_v$, and if $L/K_v$ is a
finite extension with ramification index $e$, then
$$
\calL(A\otimes L,\omega\otimes L)=e\calL(A,\omega).
$$
\item
Suppose that the abelian variety $A$ has semistable reduction over $K_v$, and $f\colon
A\rightarrow B$ is an isogeny. Then, in the notation of Proposition 1,
$$
\calL(A,f^*\omega)=\calL(B,\omega)+\calL(e^*\Omega^1_{\bfG/\bfS}),
$$
where $\omega=\Omega^g_{B/K_v}$, $\bfG=\Ker\tilde f$, $e\colon \bfS\rightarrow\bfA$ is the
identity section of the N\'eron model.
\end{enumerate}

It should be noted that our definition of the local term $\calL(A,\omega)$ can be used for
abelian varieties over any discrete valuation field $F$ of characteristic zero having a perfect
residue field of characteristic $p$ (under these assumptions the N\'eron model always exists).
Here the expression $\calL(A,\omega)$ will satisfy conditions 1)--3). In particular, the value of
this expression does not change under unramified extensions of the ground field. For example, one
can consider the local term of the logarithmic height $\calL(A,\omega)$ for abelian varieties
over the maximal unramified extension $K^{\unr}_v$ of a $v$-adic field $K_v$. From the properties
of the N\'eron model (see Lang's book, Chapter 11, \S5) it follows that the local term
$\calL(A,\omega)$ does not change if one passes to the completion $\hat K^{\unr}_v$ of the field
$K^{\unr}_v$, which is a complete discrete valuation field with algebraically closed residue
field of characteristic $p$. We note that the Galois group $I$ of the field $\hat K^{\unr}_v$
coincides with the inertia group $I(v)$ of the place $v$.

In the case of an abelian variety $A$ over the field of functions on a curve $B$, we have
$$
d(A)=\sum_{v\in B}\calL(A\otimes K_v,\omega\otimes K_v),
$$
where $d(A)$ is the invariant \eqref{eq1.6} in \S1.

We now proceed to consider the archimedean components. If $v\in{P(K)}_\infty$, then let $K_v$ be
the (topological) closure of $K$ in $\bfC$, and set
$$
\epsilon_v=\begin{cases}1,\quad\text{if}\ K_v=\bfR,\\2,\quad\text{if}\ K_v=\bfC.\end{cases}
$$
\begin{xxxr}[Definition 2] Suppose that $A$ is an abelian variety over the field $K_v$, $g=\dim A$,
$\omega\in\Omega^g_{A/K_v}$, $\omega\neq0$, is a regular form. We set
$$
h_{K_v}(A,\omega)=-\frac{\epsilon_v}{2}\log\biggl[{\Bigl(\frac i2\Bigr)}^g\int_{A(\bfC)}\omega\wedge\ol\omega\biggr].
$$
\end{xxxr}
\oripage
\begin{xxxr}[Remark 1] The expression in square brackets is connected with the period matrix. It is equal to
$$
\biggl|\det\biggl(\int_{\gamma_i}\omega_i,\int_{\gamma_{i+g}}{\ol\omega}_i\biggr)\biggr|,
$$
where $\omega_i\in\Gamma(A,\Omega^1_A)$, $\omega_1\wedge\dots\wedge\omega_g=\omega$ and
$\gamma_1,\dots,\gamma_{2g}$ is a basis of $H_1(A(\bfC),\bfZ)$.
\end{xxxr}

\begin{xxxi}[Proposition 2] The local term $h_{K_v}(A,\omega)$ has the following properties.
\begin{enumerate}[\upshape1)]
\item If $\lambda\in K^*_v$, then
$$
h_{K_v}(A,\lambda\omega)=h_{K_v}(A,\omega)-\epsilon_v\log|\lambda|.
$$
\item If $L/K_v$ is a finite extension, then
$$
[L:K_v]h_{K_v}(A,\omega)=h_L(A\otimes L,\omega\otimes L).
$$
\item If $f\colon A\rightarrow B$ is an isogeny, then
$$
h_{K_v}(A,f^*\omega)=h_{K_v}(B,\omega)-\frac{\epsilon_v}{2}\log(\deg f).
$$
\end{enumerate}
\end{xxxi}
It is an elementary exercise to verify these properties. We now consider the global situation of
an abelian variety $A$ over a number field $K$. We shall say that $A$ has the reduction
properties given above over the field $K$ if those properties are fulfilled for all
nonarchimedean places $v\in{P(K)}_f$. We let $S(A)$ denote the set of (nonarchimedean) places $v$
of $K$ at which $A$ has bad reduction. $S(A)$ is a finite set. If $B$ is an abelian variety over
$K$, then $S(A\times B)=S(A)\cup S(B)$. If $A$ and $B$ are isogenous to one another over $K$,
then $S(A)=S(B)$. In particular, $S(A)=S(A')$, where $A'$ is the Picard variety of $A$. Let $L$
be a number field containing $K$, and let $A_L$ be the abelian variety $A$ considered over $L$,
i.\,e., $A_L=A{\otimes}_K L$. Then we can consider the finite set $S(A_L)$ of places of $L$—the
set of places of bad reduction of the abelian variety $A_L$. If $v$ is a place of $K$ at which
$A$ has good reduction, i.\,e., if $v\not\in S(A)$, and if $v'$ is a place of $L$ lying over $v$,
then $v'\not\in S(A_L)$. If the field extension $L/K$ is unramified at a place $v$ of $L$ and if
$v\in S(A)$, then all of the places $v'$ of $L$ lying over $v$ belong to the set $S(A)$.

\begin{xxxr}[Definition 3] Suppose that $A$ is an abelian variety over a number field $K$, $g=\dim A$, and
$\omega\in\Omega^g_{A/K}$ is a nonzero invariant form. Then
$$
h(A)=\frac{1}{[K:\bfQ]}\sum_{v\in P(K)}h_{K_v}(A,\omega)
$$
is called the \textit{canonical height} of the abelian variety $A$.
\end{xxxr}

Using the product formula (Chapter 2 in Lang's book) and Propositions 1.1 and 2.1, we see that
$h(A)$ does not depend on the choice of the form $\omega$.

\begin{xxxi}[Proposition 3] The canonical height $h(A)$ has the following properties\r:
\begin{enumerate}[\upshape1)]
\item if $L/K$ is a finite extension and if $A$ has semistable reduction over $K$, then $h(A{\otimes}_KL)=h(A)$;
\oripage
\item if $f\colon A\rightarrow B$ is an isogeny and if $A$ has semistable reduction over $K$, then
$$
h(A)=h(B)-\frac12\log(\deg f)+\frac{1}{[K:\bfQ]}\log(\#e^*\Omega^1_{\bfG/\bfS}),
$$
where $\bfG=\Ker\tilde f$, and $\tilde f$ is the homomorphism of N\'eron models determined by
$f$;
\item $h(A\times B)=h(A)+h(B)$.
\end{enumerate}
\end{xxxi}
Properties 1) and 2) follow from Propositions 1 and 2. Property 3) can be verified directly. We
note that from 2) it follows that $2[K:\bfQ](h(B)-h(A))$ belongs to $\log\bfQ$ and is equal to a
sum of numbers of the form $\log\ell$ with $\ell|\deg f$ (the morphism $\tilde f$ is \'etale over
the places $v$ of $K$ which are prime to the divisors of $\deg f$; hence $\Omega^1_{\bfG/\bfS}$
is equal to zero over such places).

We further note that the number $\exp h(A)$ occurs in an expression for the residue of the
$\zeta$-function of an elliptic curve at the point $s=1$, according to the
Birch--Swinnerton--Dyer conjecture (see~\cite{Man2}). It is not surprising that the corresponding
expression in the geometric analog of this conjecture for a curve over a function field with
finite field of constants ${\bfF}_q$ contains the factor $q^{d(A)}$ where $d(A)$ is the invariant
\eqref{eq1.6} of \S1.

\begin{xxxi}[Finiteness theorem for the height] Suppose that $K$ is a number field,
$c\in \bfR$, and $g$ is an integer.
The set of abelian varieties $A$ over $K$ \r(considered up to $K$-isomorphism\r) for which
\begin{enumerate}[\upshape1)]
\item $\dim A=g$,
\item $h(A)\leq c$, and
\item there exists a principal \r(i.\,e., degree 1\r) polarization, is finite.
\end{enumerate}
\end{xxxi}
(For applications one can replace the inequality in 2) by equality.)

For the proof of this theorem, see~\cite{71}, Expos\'e IV (a more condensed
exposition is in Deligne's report~\cite{37}). We also have another definition
of the height of an abelian variety (discussed in~\cite{17} and~\cite{Par2}).
This modular height $h_{\mod{}}(A)$ has all of the properties of $h(A)$
(Proposition 3 and the finiteness theorem), and differs from $h(A)$ only by a
constant. We limit ourselves to the case of abelian varieties with potentially
good reduction.

\begin{xxxr}[Remark 2] In \S5 we shall also use the following variant of the height $h(A)$. If $A$ is an abelian
variety over a local field $K_v$, then by the \textit{local stable term} ${\calL}_S(A,\omega)$ we
mean the number $e^{-1}\calL(A\otimes L,\omega\otimes L)$, where $L$ is an extension of $K_v$
over which $A$ has semistable reduction and $e$ is the ramification index of the extension
$L/K_v$. This term does not depend on the choice of $L$. We set
$h_s(A)={[K:\bfQ]}^{-1}\sum{\calL}_s(A\otimes K_v,\omega\otimes K_v)$. Then, from the finiteness
theorem for the height $h(A)$, the finiteness theorem for forms, and Hermite's theorem (see \S1),
we conclude that for a given number field $K$ and dimension $g$ there exist only finitely many
(up to $K$-isomorphism) abelian varieties $A$ over $K$ of dimension $g$ having fixed $h_s(A)$ and
$\set S(A)$ (see the argument below in Remark 1).
\end{xxxr}
\oripage

\ssect[The moduli space] Here we give a brief presentation based upon \S3 of~\cite{17}. For
details of the proofs see~\cite{62}, \cite{44}, \cite{50}, and~\cite{48}.

Let ${\scrG}_g\subset\bfC^{g(g+1)/2}$ be the set of complex matrices $Z=X+iY$ which are symmetric
$(Z^t=Z)$ and have positive definite imaginary part $Y(Y>0)$.

The symplectic group $G=\Sp(2g,\bfR)$ acts on ${\scrG}_g$ by the formula
\begin{equation}\label{eq2.1}
g(Z)=(AZ+B){(CZ+D)}^{-1},\quad g=\begin{pmatrix}A&B\\C&D\end{pmatrix}.
\end{equation}

If $\Gamma\subset\Sp(2g,\bfZ)\subset G$ is a subgroup of finite index, then a complex function
$f$ on ${\scrG}_g$ is called a \textit{modular form} of weight $r$ if

\begin{enumerate}[\upshape1)]
\item $f$ is holomorphic on ${\scrG}_g$,
\item $f(\gamma(Z))=\det{(CZ+D)}^rf(Z)$, $\gamma\in\Gamma$, and
\item $f$ is bounded in the region $Y\geq Y_0$, where $Y_0>0$.
\end{enumerate}
This implies that the real function $\tilde f(Z)={(\det Y)}^{r/2}|f(Z)|$ is invariant relative to
$\Gamma$ (\cite{44}, Chapter 1). A modular form $f$ has a Fourier expansion
$$
f(Z)=\sum_{T\geq r}a(T)e^{i\pi\Tr(TZ)}
$$
and is called a \textit{cusp form} if $a(T)\neq0$ implies $T>0$. In this case $\tilde f(Z)$ has a
maximum on ${\scrG}_g$ (\cite{44}, Chapter 1).

In ${\scrG}_g$ we have a volume which is invariant relative to the group $G$; it is given by
\begin{equation}\label{eq2.2}
\Omega=\frac{dXdY}{{(\det Y)}^{g+1}},
\end{equation}
and so for a cusp form $f(Z)$ of weight $N(g+1)$ we have
\begin{equation}\label{eq2.3}
\tilde f{(Z)}^2=f(Z)\ol f(Z){(\det Y)}^{N(g+1)}\leq C.
\end{equation}
We set $\gotM=\gotM(\Gamma)={\scrG}_g/\Gamma$. If $\Gamma$ has a sufficiently large index in
$\Sp(2g,\bfZ)$, then this is a nonsingular complex manifold. One can define a corresponding
affine group $D$, which is an extension of $\Gamma$ by ${\bfZ}^{2g}$ (see~\cite{17}, \S2). Then
$\bfA={\scrG}_g\times{\bfC}^g/D$ determines an abelian scheme $q\colon \bfA\rightarrow\gotM$. For
any $N$ we have the isomorphism
\begin{equation}\label{eq2.4}
i\colon \omega^{\otimes N}_{\gotM/\bfC}={\biggl(\bigwedge^{g(g+1)/2}\Omega^1_{\gotM/\bfC}\biggr)}^{\otimes N}
\tilde\rightarrow\;\omega^{\otimes N(g+1)}_{\bfA/\gotM},
\end{equation}
which commutes with the action of $\Sp(2g,\bfZ)/\Gamma$ on $\gotM(\Gamma)$, if $\Gamma$ is a
normal divisor. The map \eqref{eq2.4} is induced by the Kodaira--Spencer mapping~\cite{50} ($T$
is the tangent bundle)
$$
T_{\gotM}\rightarrow R^1q_*T_{\bfA/\gotM}
$$
followed by the map from $R^1q_*T_{\bfA/\gotM}$ to $R^1q_*{\goto}_{\bfA}\otimes
R^1q_*{\goto}_{\bfA}$ (since for any abelian scheme $T_{\bfA/\gotM}\cong q^*e^*T_{\bfA/\gotM}$,
and then passage to the dual maps of the determinants of these sheaves.

The scheme $\gotM=\gotM(\Gamma)$ or $\gotM^{ab}_g(\Gamma)$ (the \textit{moduli space}) has the
following \textit{universal property} (if the index of $\Gamma$ is sufficiently large): for any
principally \oripage polarized abelian variety $A$ over $\bfC$ there exists a point
$a\in\gotM(\bfC)$ such that the fiber ${\bfA}_a$ of the family $\bfA/\gotM$ is isomorphic to $A$.
The sheaf $\omega_{\gotM}$ (and hence $\omega_{\bfA/\gotM}$) is ample on $\gotM$, and some power
$\omega^{\otimes N}_{\gotM}$ gives a projective imbedding of the manifold $\gotM$.

The regular sections of the sheaf $\omega_{\gotM/\bfC}$ on $\gotM$ may be regarded as modular
forms of weight $g+1$. Hence, if $\eta\in\Gamma(\gotM,\omega^{\otimes N}_{\gotM})$ is a cusp
form, then from \eqref{eq2.3} we have
\begin{equation}\label{eq2.5}
-\log\Bigl|\frac{\eta\wedge\ol\eta}{\Omega^{\otimes N}}\Bigr|\geq C.
\end{equation}
We now construct the tower of \textit{Hecke correspondences} over $\gotM(\Gamma)$. Let
$\Gamma_0(m)\subset\Sp(2g,\bfZ)$ consist of the matrices \eqref{eq2.1} with $C\equiv0\mod m$. We
set
$$
H_{g,m}=H_{g,m}(\Gamma)={\scrG}_g/\Gamma_0(m)\cap\Gamma.
$$
The map $\tau(Z)={-m}^{-1}{Z}^{-1}$ induces an involution on $H_{g,m}$ (under certain conditions
on $\Gamma$ see~\cite{17}, \S2), and hence we have the diagram
\begin{equation}\label{eq2.6}
H_{g,m}\mathop{\rightrightarrows}\limits^{p_1}_{p_2}\;\gotM,
\end{equation}
in which $p_1$ is the projection coming from the imbedding
$\Gamma_0(m)\cap\Gamma\hookrightarrow\Gamma$, and $p_2=p_1\circ\tau$. Lifting the abelian scheme
$\bfA$ over $H_{g,m}$, we obtain the diagram
\begin{equation}\label{eq2.7}
\begin{array}{rcl}
p^{-1}_1(\bfA)&\stackrel{\psi}{\longrightarrow}&p^{-1}_2(\bfA)\\
\searrow&&\swarrow\\
&H_{g,m}&
\end{array}
\end{equation}
in which $\psi$ is the isogeny of degree $m$ which is induced by the map $\delta$ of
${\scrG}_g\times{\bfC}^g$ to itself given by
$\delta(\bfZ,\omega)=(\tau(\bfZ),\omega{\bfZ}^{-1})$. The diagram \eqref{eq2.7} has the following
\textit{universal property}: for any degree $m$ isogeny $f\colon A\rightarrow B$, $A={\bfA}_a$,
$B={\bfA}_b$ (fibers of the family $\bfA/\gotM$), there exists a point $c\in H_{g,m}$ such that
$p_1(c)=a$, $p_2(c)=b$, and $\psi|_c=f$.

If the form $\eta\in\Gamma(\gotM,\omega^{\otimes N}_{\gotM})$ is chosen so that
$$
i(\eta_a)=f^*i(\eta_b)\quad\text{in}\quad\omega^{\otimes N(g+1)}_{\bfA/\gotM,a},
$$
which is possible for $N$ large (because the sheaf $\omega_{\gotM}$ is ample), then a direct
verification shows that
\begin{equation}\label{eq2.8}
{(p^*_1\eta)}_c={(\deg f)}^N{(p^*_2\eta)}_c\quad\text{in}\quad\omega^{\otimes N}_{H,c}.
\end{equation}

\ssect[The modular height] This height is obtained by combining Definition 1 in \S2.1 for the
nonarchimedean places $v$ with a definition of the archimedean components as functions on the
moduli space $\gotM(\Gamma)$ or $\gotM^{ab}_g$. In order to obtain the necessary estimates, we
must have a moduli scheme $\bfM$ over $\Spec\bfZ$ with universal abelian scheme $\bfA/\bfM$. They
must be models over $\bfZ$ of the moduli spaces $\gotM(\Gamma)$ which were introduced for $\bfC$
above. Such a scheme exists, if we consider rigidified abelian varieties (and then
$\bfM\otimes\bfC=\gotM(\Gamma)$ for some subgroup $\Gamma$) \cite{53}. Unfortunately, such a
theorem \oripage can be proved not over $\Spec\bfZ$ but over an open subset of it. There are
several ways of getting around this difficulty. Since all of them involve some technical
complications, we shall suppose in what follows that the moduli scheme exists over $\bfZ$.

\begin{xxxr}[Remark added in translation] Gabber's lemma (item 21 in the diagram in \S 1, and Expos\'e 5 of
\cite{71}) makes it possible to avoid both this difficulty and also the difficulties at points of
bad reduction. Namely, from the moduli scheme $M$ over $\Spec\bfZ$ --- (finite number of points)
and the abelian scheme $A$ over $M$ (and the sheaf $\omega$), using that lemma, one can construct
a proper ``scheme'' $M'$ over all of $\Spec\bfZ$ and a semiabelian ``scheme'' $A'\rightarrow M'$
(and a sheaf $\omega'$ on $M'$) which extend the original objects on $M$. See the formulation of
Gabber's lemma on p.\,376. Then the one-dimensional subscheme $C\in M$ in Definition 2 extends to
a subscheme $C'\in M'$ which is proper over $\Spec\bfZ$. Applying the estimates of this section
to the scheme $M'$ and the sheaf $\omega'$, we obtain a finiteness theorem for the height in the
general case (see also~\cite{71}, Expos\'e 4, no. 2).
\end{xxxr}
\begin{xxxr}
[Definition 1] Suppose that $a\in\gotM^{ab}_g(\bfC)$, $A={\bfA}_a$ is the corresponding abelian
variety over $\bfC$, $\eta\in{\Omega^g_{\gotM/\bfC}}^{\otimes N}$ is a rational form on
$\gotM^{ab}_g$, $a\not\in\Supp(\eta)$, and $N$ is such that the sheaf $\omega^{\otimes
N}_{\gotM}$ is very ample. Then
$$
h_{\mod,v}(A,\eta)=-\frac{\epsilon_v}{2}\log\Bigl|\frac{\eta\wedge\ol\eta}{\Omega^{\otimes N}}\Bigr|(a),
$$
where $\Omega$ is the volume form on $\gotM^{ab}_g(\bfC)$ (see \eqref{eq2.2} in \S2.2).
\end{xxxr}

Since the volume form $\Omega$ is invariant relative to the entire group $\Sp(2g,\bfR)$, it
follows that the functions $h_{\mod,v}$ for moduli spaces $\gotM(\Gamma)$ with different groups
$\Gamma$ are in a natural way compatible with one another. Like the functions in \S2.1, they have
the following properties.

\begin{xxxi}[Properties of the function $h_{\mod,v}$]
\begin{enumerate}[\upshape1)]
\item If $f\in K_v{(\gotM)}^*$, $K_v=\bfR$, $\bfC$, then
$$
h_{\mod,v}(A,f\eta)=h_{\mod,v}(A,\eta)-N\epsilon_v\log|f|.
$$
\item Suppose that $f\colon A\rightarrow B$ is an isogeny, $\omega\in{(\Omega^g_{B/K_v})}^{\otimes N}$,
and the form $\eta\in\Gamma(\gotM,\omega^{\otimes N}_{\gotM})$ is such that
$$
i(\eta)|_a=f^*\omega,\quad i(\eta)|_b=\omega,
$$
where $a$ and $b$ are the points of $\gotM(K_v)$ corresponding to $A$ and $B$. Then
$$
h_{\mod,v}(A,\eta)=h_{\mod,v}(B,\eta)-\frac{N\epsilon_v}{2}\log(\deg f).
$$
\end{enumerate}
\end{xxxi}
\begin{xxxr}[Proof] Property 1) is obvious. Property 2) follows from diagram \eqref{eq2.7} and formula \eqref{eq2.8}
in \S2.2 (we note that in~\cite{17} it was mistakenly concluded from Proposition 6 in \S4 that
the number $h_{\mod,v}(A,\eta)$ is \textit{invariant} relative to isogenies, and the question was
posed of finding nonarchimedean components with the same property).
\end{xxxr}

\oripage We shall need a slight generalization of the number $h_{K_v}(A,\omega)$ in Definition 1
of \S2.1. Namely, we can immediately generalize to the case of a tensor form
$\omega\in{(\Omega^g_{A/K})}^{\otimes N}$ (and a sheaf $\omega^{\otimes N}_{\bfA/\bfS}$), with
Proposition 1 still holding (in 3) one must multiply the second term on the right by $N$).

\begin{xxxr}[Definition 2] Suppose that $K$ is a number field, $A$ is a $g$-dimensional abelian variety with
good reduction over $K$, $\eta\in\Gamma(M,\omega^{\otimes N}_{\gotM/K})$ is a holomorphic form,
and $\gotM=\gotM^{ab}_g$. To these data correspond:
\begin{enumerate}[\upshape1)]
\item the point $a\in\gotM(K)$ for which ${\bfA}_a=A$, and the corresponding points $a_v\in\gotM(\bfC)$,
$v\in{P(K)}_\infty$ for $A\otimes K_v$;

\item the form $\omega=i(\eta)\in\omega^{\otimes N(g+1)}_{\bfA/\gotM}$ obtained using the isomorphism $i$
(see \eqref{eq2.4} of \S2.2); and

\item a one-dimensional subscheme $\bfC\subset\bfM$, $\bfC\simeq\Spec\Lambda_K$, such that $\bfA|_\bfC$ is the N\'eron model of
$A$ over $K$ and $a$ is its generic point.
\end{enumerate}
Then the quantity
$$
h_{\mod}(A)=\frac{1}{N[K:\bfQ]}\sum_{v\in{P(K)}_f}h_{K_v}(A\otimes K_v,\omega\otimes K_v)+
\frac{1}{N[K:\bfQ]}\sum_{v\in{P(K)}_\infty}h_{\mod,v}(A,\eta)(a_v)
$$
is called the \textit{modular height} of the abelian variety $A$.
\end{xxxr}

Using the product formula, we see that $h_{\mod}(A)$, like $h(A)$, does not depend upon the
choice of $\eta$ and makes sense for any $A$. From Proposition 1 of \S2.1 and the properties of
the archimedean components $h_{\mod,v}$, we find that all of the properties in Proposition 3
still hold for the modular height. Thus, it remains only to consider the finiteness properties.

\begin{xxxi}[Proposition] The set of principally polarized $g$-dimensional abelian varieties $A$ with good
reduction over $K$ having a given height $h_{\mod}(A)$ is finite.
\end{xxxi}

\begin{xxxr}[Remark 1] The proposition remains true if we suppose only that the abelian varieties $A$ have
potentially good reduction but that the set $S(A)$ (of places of $K$ of bad reduction) is fixed.
In fact, according to Raynaud's criterion (\S2.1), if we adjoin the coordinates of the points of
order 3 to the field $K$, we obtain an extension $L$ over which $A$ has semistable and hence good
reduction. Since we have fixed the set $S(A)$, the extension $L/K$ has bounded ramification (and
degree). By Hermite's theorem, there are only finitely many extensions $L/K$, and it then remains
to apply the finiteness theorem for forms in \S3.
\end{xxxr}

We further note that, since $S(A)$ is invariant under isogenies, in all of the applications (in
\S5) the fixing of this set has no significance.

In order to prove the proposition, we consider a model $\bfM$ of the moduli scheme over $\bfZ$.
Let $\calL=\omega^{\otimes N(g+1)}_{\bfA/\bfM}$ and $\eta\in\Gamma(\gotM,\omega^{\otimes
N}_{\bfM/K})$. We consider the section $i(\eta)$ (see \eqref{eq2.4} in \S2.2) of the sheaf
$\omega^{\otimes N(g+1)}_{\bfA/\gotM}$. It determines a rational \oripage section $\tau=i(\eta)$
of the sheaf $\calL$ over the scheme $\bfM$. By property 3) of Definition 2 and by Definition 1
of \S2.1, we have
\begin{equation}\label{eq1}
h_{K_v}(A\otimes K_v,i(\eta)\otimes K_v)=v(\tau|_{\bfC})\log\bfN v
\end{equation}
for any place $v\in{P(K)}_f$ of $K$.

\begin{xxxi}[Lemma] Suppose that $X$ is an irreducible algebraic variety over a local field $K_v$, $\bfY/{\goto}_v$ is a smooth
normal scheme with generic fiber $X$, $\calL\in\Pic\bfY$, $s\in\Gamma(X,\calL|X)$. Define the
function
$$
f(P)=v(\tilde s|_{\bfC}),\quad P\in X(K_v),
$$
where $\bfC$ is the section of the scheme $\bfY$ over $\goto$ which corresponds to $P$, and
$\tilde s$ is a rational section $\calL$ which extends $s$. Then the following assertions are
true\r:
\begin{enumerate}[\upshape1)]
\item the function $f$ is bounded from below;
\item if the section $\tilde s$ is regular, then $f\geq0$.
\end{enumerate}
\end{xxxi}
In fact, $f$ is the quasifunction which N\'eron [Ne3] associates to the divisor $(s)$ of the
section $s$ (or the N\'eron divisor, in Lang's terminology, Chapter 10, \S2). Properties 1) and
2) then follow from the fact that $(s)$ is an effective divisor. We shall omit the verification
of the lemma.

Now suppose that $\eta_0,\dots,\eta_n\in\Gamma(\gotM,\omega^{\otimes N}_{\gotM})$ is a basis for
the space of cusp forms of some sufficiently large weight $N$. Then the $(n+1)$-tuple
$\eta_0,\dots,\eta_n$ determines a projective imbedding $\varphi\colon
\gotM\hookrightarrow{\bfP}_n$ \cite{48}. We estimate the height $h_\varphi(a)$ of the point
$a\in\gotM(K)$ in this imbedding (see the definition in Chapter 3). Let $s_j=i(\eta_j)$ and
$s_0(a)\neq0$. By \eqref{eq2.1}, Proposition 1.1 of \S2.1, and property 1) of the function
$h_{\mod,v}$, we have
\begin{align*}
[K:\bfQ]h_\varphi(a)&=\sum_{v\in P(K)}\sup\limits_j\log\biggl|\frac{s_j}{s_0}(a)\biggr|=-\sum_{v\in P(K)}\inf\limits_jv\biggl(\frac{s_j}{s_0}(a)\biggr)={}\\
&{}=-\sum_{v\in{P(K)}_f}\inf\limits_jv({\tilde s}_j|_{\bfC})\log\bfN v+\sum_{v\in{P(K)}_f}v({\tilde s}_0|_{\bfC})\log\bfN v-{}\\
&{}\quad -\sum_{v\in{P(K)}_\infty}\inf\limits_jh_{\mod,v}(A,\eta_j)+\sum_{v\in{P(K)}_\infty}h_{\mod,v}(A,\eta_0)\leq{}\\
&{}\leq\const+\sum_{v\in{P(K)}_f}v({\tilde s}_0|_{\bfC})\log\bfN v+{}\\
&{}\quad +\sum_{v\in{P(K)}_\infty}h_{\mod,v}(A,\eta_0)\leq\const+N[K:\bfQ]h_{\mod}(A),
\end{align*}
where the inequalities follow from the lemma (the divisors $({\tilde s}_j)$ are effective on the
scheme $\bfM$, if we remove from it a finite number of fibers) and the bound \eqref{eq2.5} in
\S2.2. Applying the finiteness property for heights (Lang, Chapter 3, Theorem 2.6), we obtain the
required result.

\begin{xxxr}[Remark 2] It follows from the proposition that the number of polarizations of degree 1 on a
fixed abelian variety over $K$ (considered up to $K$-isomorphism) is finite.
\end{xxxr}
\oripage
\begin{xxxr}[Remark 3] We now show that $h_{\mod}(A)=h(A)+\const$. Let $\omega\in\Gamma(\gotM,\omega_{\bfA/\gotM})$.
Then the form $\omega\wedge\ol\omega$ can be integrated over each fiber
${\bfA}_a$ of the family $\bfA/\gotM$. By the same token there is a scalar
product defined on the fibers of $\omega_{\bfA/\gotM}$ (and also
$\omega^{\otimes N}_{\bfA/\gotM}$). The isomorphism \eqref{eq2.4} takes this
scalar product to a scalar product on $\omega_{\gotM}$. Thus, one obtains a
smooth volume form $\Omega'$ on the tangent bundle $T_{\gotM}$. For any point
$a\in\gotM$ we have the imbeddings $T_{\gotM,a}\hookrightarrow{\scrG}_g$ and
$T_{\bfA/\gotM,e(a)}\hookrightarrow{\bfC}^g$, and the explicit formula for the
Kodaira--Spencer map in terms of these imbeddings (\cite{50}, p.\,408) shows
that $\Omega'=\Omega$ up to a constant. This implies the required equality for
the heights. One can also observe that the formula \eqref{eq2.8} giving the
behavior of the archimedean component of the modular height under isogenies
implies that the form $\Omega'$ is invariant under transformations of
${\scrG}_g$, in $\Sp(2g,\bfQ)$; and since the latter group is dense in
$\Sp(2g,\bfR)$, it follows that the invariance holds relative to the entire
group of automorphisms of ${\scrG}_g$. Consequently, $\Omega'$ is an invariant
volume in ${\scrG}_g$, and so $\Omega'=\Omega$ up to a constant (this proof is
due to Zarkhin).
\end{xxxr}

\begin{xxxr}[Remark 4] An explicit form for the heights $h(A)$ and $h_{\mod}(A)$ on a modular curve with
$g=1$ is given in Deligne's report~\cite{37}.\footletter{\textit{Added in
translation}. A table of the heights of elliptic curves of CM-type for
imaginary quadratic fields with discriminant $-p\geq-2000$, where $p$ is a
prime and $p\equiv3\mod4$, has been constructed by Zarkhin~\cite{2*}.}
\end{xxxr}

\section{$\ell$-divisible groups, Tate modules and abelian varieties}

In this section $\ell$ is a prime number, $F$ is a field of characteristic $\neq\ell$, $\ol F$ is
its separable algebraic closure, and $\Gamma$ is its Galois group.

\begin{xxxi}[Lemma] Let $W_n$ be a finite abelian group which is annihilated by multiplication by ${\ell}^n$. Let
$W_1$ denote the subgroup of $W_n$ consisting of elements annihilated by multiplication by
$\ell$, and let ${\ell}^n$ be the order of this subgroup, i.\,e., $h$ is the dimension of $W_1$
as a vector space over $\bfZ/\ell\bfZ$. Then the following conditions are equivalent\r:
\begin{enumerate}[\upshape\r(a\r)]
\item the group $W_n$ has order equal to ${\ell}^{nh}$;
\item $W_n$ is a free $\bfZ/{\ell}^n\bfZ$-module of rank $h$.
\end{enumerate}
\end{xxxi}
The proof follows immediately from the elementary divisor theorem.

\begin{xxxr}[Remark] Let $m$ denote the greatest natural number such that $W_m=\Ker(W_n\stackrel{{\ell}^m}{\rightarrow}W_n)$ is a free
$\bfZ/{\ell}^m\bfZ$-module, and set $\tau(W_n)=W_m$. Then $W_m$ is a free
$\bfZ/{\ell}^m\bfZ$-module of rank $h$. Conditions (a) and (b) are equivalent to the equality
$W_n=\tau(W_n)$.

If $W_n$ is a direct sum of $r$ cyclic groups, then $W_n/\tau(W_n)$ is a direct sum of at most
$r-1$ cyclic groups.
\end{xxxr}

\begin{xxxr}[Example] Let $G$ be an (infinite) abelian group which satisfies the following conditions: the map
multiplication by $\ell$, $G\stackrel{\ell}{\rightarrow}G$, is surjective, and its kernel $G_1$
is a finite abelian group. Then, for any natural number $n$, the kernel $G_n$ \oripage of
multiplication by ${\ell}^n$, $G\stackrel{{\ell}^n}{\rightarrow}G$, is a free
$\bfZ/{\ell}^n\bfZ$-module of rank $h$, where $h$ is the dimension of $G_1$ as a vector space
over $\bfZ/\ell\bfZ$.
\end{xxxr}

Let $h$ be a nonnegative integer.

\begin{xxxr}[Definition (Tate~\cite{28})] An \textit{$\ell$-divisible group $W$ of height h}
\footnote{Not to be confused with the notion of height in Diophantine geometry!} is an infinite
sequence of pairs $(W_n,i_n)$, consisting of a finite abelian group $W_n$ and an imbedding
$i_n\colon W_n\hookrightarrow W_{n+1}$, which satisfies the following conditions.

\begin{enumerate}[\upshape(a)]
\item The group $W_n$ has order ${\ell}^{nh}$.
\item The imbedding $i_n$ identifies the group $W_n$ with the kernel of multiplication by ${\ell}^n$ in
$W_{n+1}$ i.\,e., one has the exact sequence
$$
0\rightarrow W_n\stackrel{i_n}{\rightarrow}W_{n+1}\stackrel{{\ell}^n}{\rightarrow}W_{n+1}.
$$
\end{enumerate}
In particular, the group $W_n$ is annihilated by multiplication by ${\ell}^n$.
\end{xxxr}

It follows from the lemma that $W_n$ is a free $\bfZ/{\ell}^n\bfZ$-module of rank $h$. The
composition of the imbeddings $i_r$ $(r=n,\dots,n+m-1)$ gives an imbedding
$$
i_{n,m}\colon W_n\hookrightarrow W_{n+m},
$$
which identifies the group $W_n$ with the kernel of multiplication by ${\ell}^n$ in $W_{n+m}$,
i.\,e., one has the exact sequence
$$
0\rightarrow W_n\stackrel{i_{n,m}}{\rightarrow}W_{n+m}\stackrel{{\ell}^n}{\rightarrow}W_{n+m}.
$$
Note that
$$
i_nW_n=\ell W_{n+1},\quad i_{n,m}W_n={\ell}^mW_{n+m}.
$$

By the same token, we have the exact sequence
$$
0\rightarrow W_n\stackrel{i_{n,m}}{\rightarrow}W_{n+m}\rightarrow W_m\rightarrow0.
$$
One can define the group
$$
W_\infty=\lim_{\substack{\longrightarrow\\n}}W_n,
$$
which is the inductive limit relative to the imbeddings $i_n$, and the Tate module
$$
T(W)=\lim_{\substack{\longleftarrow\\n}}W_n,
$$
which is the projective limit relative to the surjections
$$
W_{n+1}\stackrel{\ell}{\rightarrow}\ell W_{n+1}\stackrel{i^{-1}_n}{\rightarrow}W_n.
$$
The group $W_\infty$ is isomorphic to ${({\bfQ}_{\ell}/{\bfZ}_{\ell})}^h$ and the Tate module
$T(W)$ is a free ${\bfZ}_{\ell}$-module of rank $h$. Multiplication by $\ell$ in $W_\infty$ is
surjective, and its kernel is canonically isomorphic to $W_1$. For any $n$ the group $W_n$ is
canonically isomorphic to the kernel of multiplication by ${\ell}^n$ in $W_\infty$. These groups
can also be recovered from the Tate module: there exist canonical isomorphisms
$$
W_n\simeq T(W)/{\ell}^nT(W).
$$
\oripage

An important example of an $\ell$-divisible group is the group
$$
G_n=\Ker(G)\stackrel{{\ell}^n}{\rightarrow}G,\quad G_n\subset G_{n+1},
$$
which is the kernel of multiplication by ${\ell}^n$ in an abelian group $G$ for which
multiplication by $\ell$ is surjective and the kernel $G_1$ of multiplication by $\ell$ is a
finite group; it is the dimension $h$ of the latter group over the field $\bfZ/\ell\bfZ$ which is
the height of the $\ell$-divisible group $G(\ell)=(G_n,i_n)$. Here $i_n\colon G_n\subset G_{n+1}$
is the inclusion map. The group $G{(\ell)}_\infty$ can be identified in a natural way with the
$\ell$-torsion subgroup of $G$, and the Tate module $T(G(\ell))$ coincides with the Tate
$\ell$-module (or ${\bfZ}_{\ell}$-module) $T_{\ell}(G)$ of $G$, which is defined as the
projective limit
$$
\lim_{\substack{\longleftarrow\\n}} G_n
$$
of the groups $G_n$ relative to the transfer maps $G_{n+1}\stackrel{\ell}{\rightarrow}G_n$.

Along with the Tate ${\bfZ}_{\ell}$-module $T_{\ell}(G)$ one often considers the Tate
${\bfQ}_{\ell}$-module
$$
V_{\ell}(G)=T_{\ell}(G)\otimes_{{\bfZ}_{\ell}}{\bfQ}_{\ell},
$$
which is an $h$-dimensional vector space over ${\bfQ}_{\ell}$.

\begin{xxxr}[Examples]
\begin{enumerate}[\upshape1)]
\item Let $G=\bfQ/\bfZ$. Then $G_n={\ell}^{-n}\bfZ/\bfZ$, the $\ell$-divisible group $\bfQ/\bfZ(\ell)$ has height 1, and
$$
G_\infty={\bfQ}_{\ell}/{\bfZ}_{\ell},\quad T_{\ell}(G)={\bfZ}_{\ell},\quad V_{\ell}(G)={\bfZ}_{\ell}.
$$
\item Let $G={\ol F}^*$. Then $G_n=\mu_{{\ell}^n}$ and the $\ell$-divisible group ${\bfG}_m(\ell)=G(\ell)$ also has height 1.
Here $\mu_{{\ell}^n}$ is the group of ${\ell}^n$th roots of 1. We note that the Galois group
$\Gamma$ of $F$ acts compatibly on the groups $G_n$.

\item Let $X$ be an abelian variety over a field $F$. We set $G=X(\ol F)$. Then one can consider the $\ell$-divisible group $X(\ell)=G(\ell)$.
Its height is $2g$, where $g=\dim X$, and $X{(\ell)}_n=X_{{\ell}^n}$, where $X_{{\ell}^n}$ is the
group of points of the abelian variety $X$ which are annihilated by multiplication by ${\ell}^n$.
The Tate module $T(X(\ell))=T_{\ell}(G)$ coincides with the Tate $\ell$-module $T_{\ell}(X)$ of
the abelian variety $X$. As in the previous example, the Galois group acts compatibly on the
groups $G_n=X_{{\ell}^n}$, and these actions glue together to form a continuous homomorphism
$\Gamma\rightarrow\Aut T_{\ell}(X)$ which extends by ${\bfQ}_{\ell}$-linearity to an
\textit{$\ell$-adic representation}
$$
\Gamma\rightarrow\Aut V_{\ell}(X),
$$
where
$$
V_{\ell}(X)=V_{\ell}(X(\ol F))=T_{\ell}(X)\otimes_{{\bfZ}_{\ell}}{\bfQ}_{\ell}.
$$
We let
$$
\chi(X,\ell)\colon \Gamma\rightarrow{\bfZ}_{\ell}\subset{\bfQ}_{\ell}
$$
denote the character of this representation. The ring End $X$ of $F$-endomorphisms of the
abelian variety $X$ acts compatibly on the groups $G_n=X_{{\ell}^n}$ and \oripage commutes
with the action of the Galois group; these actions glue together to form natural imbeddings
\begin{align*}
\End X&\otimes{\bfZ}_{\ell}\hookrightarrow\End\nolimits_{\Gamma}T_{\ell}(X),\\
\End X&\otimes{\bfQ}_{\ell}\hookrightarrow\End\nolimits_{\Gamma}V_{\ell}(X).
\end{align*}
\end{enumerate}
\end{xxxr}
\begin{xxxr}[Definition] An $\ell$-divisible group $W$ is said to be an \textit{$\ell$-divisible group over the field $F$} if the
Galois group $\Gamma$ of the field acts compatibly on all of the $W_n$, i.\,e., if all of the
$W_n$ are Galois modules and the $i_n$ are homomorphisms of Galois modules.
\end{xxxr}

If $W$ is an $\ell$-divisible group over the field $F$, then the group $W_\infty$ and the Tate
module $T(W)$ are Galois modules; in particular, there is a canonical continuous homomorphism
$$
\Gamma\rightarrow\Aut T(W)\subset\Aut V(W),
$$
where $V(W)=T(W)\otimes_{{\bfZ}_{\ell}}{\bfQ}_{\ell}$ is an $h$-dimensional vector space over
${\bfQ}_{\ell}$.

The above examples 1 )-3) of $\ell$-divisible groups are actually examples of $\ell$-divisible
groups over the field $F$ (the Galois group acts trivially in Example 1)).

One defines the notions of a homomorphism of $\ell$-divisible groups (over $F$) and an
$\ell$-divisible subgroup (over $F$) in the obvious way. If $W$ and $H$ are $\ell$-divisible
groups over $F$, then the natural map
$$
\Hom(W,H)\rightarrow{\Hom}_\Gamma(T(W),T(H))
$$
is a bijection. If $W$ is an $\ell$-divisible subgroup of $H$ (over $F$), then its Tate module
$T(W)$ is canonically imbedded in the Tate module $T(H)$, and is a pure (and $\Gamma$-invariant)
${\bfZ}_{\ell}$-submodule. Conversely, any such submodule of the Tate module determines an
$\ell$-divisible subgroup (over $F$). These sub- modules can be simply described using
($\Gamma$-invariant) subspaces $E$ of $V(H)$; the corresponding submodule turns out to be the
intersection of the subspace $E$ with the lattice $T(H)\subset
V(H)=T(H)\otimes_{{\bfZ}_{\ell}}{\bfQ}_{\ell}$.

Suppose that $X$ and $Y$ are abelian varieties over $F$. Any $F$-homomorphism $\pi\colon
X\rightarrow Y$ gives rise to a corresponding homomorphism of $\ell$-divisible groups
$X(\ell)\rightarrow Y(\ell)$ and of Tate modules $T_{\ell}(X)\rightarrow T_{\ell}(Y)$,
$V_{\ell}(X)\rightarrow V_{\ell}(Y)$. If $\pi$ is an isogeny, then the map
$V_{\ell}(X)\rightarrow V_{\ell}(Y)$ is an isomorphism of vector spaces which commutes with the
action of the Galois group. In particular, if $X$ and $Y$ are $F$-isogenous, then the Galois
modules $V_{\ell}(X)$ and $V_{\ell}(Y)$ are isomorphic, and $\chi(X,\ell)=\chi(Y,\ell)$ In the
general case, one has natural imbeddings
\begin{align*}
\Hom(X,Y)&\otimes{{\bfZ}_{\ell}}\hookrightarrow\Hom\nolimits_{\Gamma}(T_{\ell}(X),T_{\ell}(Y)),\\
\Hom(X,Y)&\otimes{{\bfQ}_{\ell}}\hookrightarrow\Hom\nolimits_{\Gamma}(V_{\ell}(X),V_{\ell}(Y))
\end{align*}
(which coincide with the ones given above when $X=Y$). We note that \textit{if
$\chi(X,\ell)=\chi(Y,\ell)$ and if the Galois modules $V_{\ell}(X)$ and $V_{\ell}(Y)$ are
semisimple, then they are isomorphic.}

\oripage

There is a one-to-one correspondence between $\ell$-divisible subgroups $W$ of $X(\ell)$ over $F$
and $\Gamma$-invariant ${\bfQ}_{\ell}$-subspaces $E$ of $V_{\ell}(X)$: the Tate module $T(W)$ is
the intersection of $E$ with $T_{\ell}(X)$. Since the Galois modules $V_{\ell}(X)$ and
$V_{\ell}(Y)$ are isomorphic whenever $X$ and $Y$ are isogenous abelian varieties, there is a
one-to-one correspondence between the $\ell$-divisible subgroups over $F$ in $X(\ell)$ and
$Y(\ell)$ which can be described as follows. Suppose that $\pi\colon X\rightarrow Y$ is an
isogeny. Then an $\ell$-divisible subgroup $W$ of $X(\ell)$ corresponds to the $\ell$-divisible
subgroup $\pi_*W$ in $Y(\ell)$, which is determined by the condition
$$
{(\pi_*W)}_\infty=\pi(W_\infty)\subset Y(\ol F).
$$
\begin{xxxr}[Definition] Suppose that $X$ is an abelian variety over the field $F$. It is said to satisfy the
\textit{$\ell$-finiteness condition} if the following holds for any $\ell$-divisible subgroup $W$
of $X(\ell)$ over $F$.
\end{xxxr}

For any natural number $n$, we define an abelian variety $Y_{(n)}=X/W_n$ over $F$ by taking the
quotient of $X$ by the finite group $W_n$. Then the sequence $Y_{(1)}$,
$Y_{(2)},\dots,Y_{(n)},\dots$ contains only finitely many pairwise nonisomorphic (over $F$)
abelian varieties, i.\,e., there exists a finite set $Y_{(i_1)},\dots,Y_{(i_r)}$ of abelian
varieties such that any of the $Y_{(n)}$ is isomorphic to one of $Y_{(i_1)},\dots,Y_{(i_r)}$.

\begin{xxxr}[Remark] If $X$ and $Y$ are abelian varieties which are $F$-isogenous to one another, and $X$ satisfies
the $\ell$-finiteness condition, then so does $Y$.
\end{xxxr}

\begin{xxxi}[Finiteness theorem for $\ell$-isogenies] Suppose that $X$ is an abelian variety over $F$ which satisfies
the $\ell$-finiteness condition \r(see Corollary 5.4\r). Then the set of abelian varieties $Y$
over $F$ \r(considered up to $F$-isomorphism\r) for which there exists an $\ell$-isogeny
$X\rightarrow Y$ \r(over $F$\r) is finite.
\end{xxxi}

One can give a proof by contradiction using the fact that a projective limit of nonempty finite
sets is nonempty, along with the following concept.

\begin{xxxr}[Definition] An $\ell$-isogeny of abelian varieties $\varphi\colon A\rightarrow B$ which is defined over $F$ is said to be
\textit{minimal} if for any other $\ell$-isogeny $\psi\colon A\rightarrow B$ defined over $F$ one
has $\deg\varphi\leq\deg\psi$.
\end{xxxr}

If two abelian varieties $A$ and $B$ are $\ell$-isogenous over $F$, then there always exists a
minimal isogeny $\varphi\colon A\rightarrow B$. If the minimal isogeny $\varphi\colon
A\rightarrow B$ factors into a composition of $\ell$-isogenies $\alpha\colon A\rightarrow C$ and
$\beta\colon C\rightarrow B$ defined over $F$, i.\,e., $\varphi=\beta\alpha$, then $\alpha$ and
$\beta$ are also minimal isogenies. In particular, if $H$ is a Galois submodule inside the kernel
$\Ker\varphi$ of a minimal isogeny $\varphi\colon A\rightarrow B$, then the natural isogenies
$$
A\rightarrow A/H,\quad A/H\rightarrow(A/H)/(\Ker\varphi/H)\simeq B
$$
are minimal. For example, for $H$ one can take the group
$$
{(\Ker\varphi)}_{{\ell}^m}=\Ker\varphi\cap A_{{\ell}^m}
$$
---the kernel of multiplication by ${\ell}^m$ in $\Ker\varphi$. We choose the largest natural number $m$ such that
${(\Ker\varphi)}_{{\ell}^m}$ is a free $\bfZ/{\ell}^m\bfZ$-module, and we set
$$
\tau(\Ker\varphi)={(\Ker\varphi)}_{{\ell}^m}.
$$
\oripage If $\Ker\varphi$ is a direct sum of $h$ cyclic groups, then $\tau(\Ker\varphi)$ is a
free $\bfZ/{\ell}^m\bfZ$-module of rank $h$, and the quotient group
$\Ker\varphi/\tau(\Ker\varphi)$ is a direct sum of at most $h-1$ cyclic groups (see the remark at
the beginning of this section). Thus, the natural isogenies
$$
A\rightarrow A/\tau(\Ker\varphi),\quad A/\tau(\Ker\varphi)\rightarrow A/\tau(\Ker\varphi)/(\Ker\varphi/\tau(\Ker\varphi))\simeq B
$$
are minimal.

If $A\rightarrow B_{(n)}$ is an infinite sequence of minimal isogenies whose degrees approach
infinity, then the sequence $B_{(1)}$, $B_{(1)},\dots$ contains an infinite number of pairwise
nonisomorphic (over $F$) abelian varieties.

To prove the theorem it suffices to obtain a contradiction from the supposition that there exists
an infinite sequence of minimal isogenies $\varphi_n\colon X\rightarrow Y_{(n)}$ whose degrees
approach infinity. Passing to a subsequence, we may suppose that the kernels $\Ker{\varphi_n}$
are direct sums of the same number $h$ of cyclic groups. We shall use induction on $h$.

Suppose that $h=1$, i.\,e., the $\Ker{\varphi_n}$ are cyclic groups. For any natural number $m$
we let $\calU_{m,1}$ denote the set of all cyclic subgroups of order ${\ell}^m$ in $X$ which have
the form ${\ell}^i\Ker_{\varphi_j}$ (for some $i$ and $j$). It is easy to see that $\calU_{m,1}$
is always a nonempty and finite set (since the set of all cyclic subgroups of $X$ of a given
order is finite). We have natural maps ${\calU}_{m+1,1}\rightarrow{\calU}_{m,1}$, which take the
group $H\in{\calU}_{m+1,1}$, to $\ell H\in{\calU}_{m,1}$. Since the projective limit
$\lim\limits_{\longleftarrow}{\calU}_{m,1}$ is nonempty, it follows that there exists an
$\ell$-divisible subgroup $W$ of $X(\ell)$ defined over $F$ which has height 1 and has the
property that for any $n$ one has $W_n={\ell}^i\Ker\varphi_j$ for some $i$ and $j$. In
particular, $W_n$ is in the kernel $\Ker\varphi_j$ of a minimal isogeny, and so the natural
isogeny $X\rightarrow X/W_n$ is minimal. We obtain an infinite sequence of minimal isogenies
$X\rightarrow X/W_n$ whose degrees, which are ${\ell}^n$, approach infinity; hence, the sequence
$X/W_1$, $X/W_2,\dots$ contains an infinite number of pairwise nonisomorphic abelian varieties.
But this contradicts the $\ell$-finiteness condition. (We essentially already proved the theorem
for elliptic curves.)

Now suppose that $h>1$. We consider the sequence of groups
$$
H_i=\tau(\Ker{\varphi_i})\subset\Ker{\varphi_i}.
$$
Two cases are possible.

\begin{enumerate}[\upshape a)]
\item The orders of the groups $H_i$ are bounded. Then, passing to a subsequence, we may suppose that
there is a finite group $H$ in $X$ such that $H=H_i=\tau(\Ker{\varphi_i})$ for all $i$.
Setting $\tilde X=X/H$, we obtain the sequence of minimal isogenies
$$
\psi_i\colon \tilde X=X/H\rightarrow X/H/(\Ker{\varphi_i}/H)=Y_{(i)}
$$
with kernels
$\Ker\psi_i=\Ker{\varphi_i}/H=\Ker{\varphi_i}/\tau(\Ker{\varphi_i})$ each of
which is a direct sum of at most $h-1$ cyclic groups. Passing to a subsequence,
we may suppose that the $\Ker\psi_i$ are direct sums of a certain number $h'$
of cyclic groups, where $h'<h$.It remains to apply the induction assumption,
since the abelian variety $\tilde X$, which is isogenous to $X$, also satisfies
the $\ell$-finiteness condition. \oripage \item The orders of the groups $H_i$
are unbounded. Then, passing to a subsequence, we may suppose that the orders
of the $H_i$ approach infinity. Let ${\calU}_{m,h}$ denote the set of finite
groups in $X$ which are free $\bfZ/{\ell}^m\bfZ$-modules of rank $h$ and have
the form ${\ell}^iH_j$ (for some $i$ and $j$). The set ${\calU}_{m,h}$ is
nonempty and finite for any natural number $m$. We have natural maps
${\calU}_{m+1,h}\rightarrow{\calU}_{m,h}$, which take the group
$H\in{\calU}_{m+1,h}$ to $\ell H\in{\calU}_{m,h}$ Since the projective limit
$\lim\limits_{\longleftarrow}{\calU}_{m,h}$ is nonempty, it follows that there
exists an $\ell$-divisible subgroup $W$ of $X(\ell)$ defined over $F$ which has
height $h$ and has the property that for any $n$ one has
$W_n={\ell}^iH_j={\ell}^i\tau(\Ker\varphi_j)$ for some $i$ and $j$. In
particular, $W_n$ is in the kernel $\Ker\varphi_j$ of a minimal isogeny, and so
the natural isogeny $X\rightarrow X/W_n$ is minimal. We obtain an infinite
sequence of minimal isogenies $X\rightarrow X/W_n$ whose degrees, which equal
${\ell}^{nh}$, approach infinity; hence, the sequence $X/W_1$, $X/W_2,\dots$
contains an infinite number of pairwise nonisomorphic abelian varieties. But
this contradicts the $\ell$-finiteness condition.
\end{enumerate}

From the proof of Proposition 1 in Tate's paper~\cite{27}, \S2, it is easy to extract the
following result:

\begin{xxxi}[Theorem on semisimplicity of the Tate module (\cite{27}, \cite{11}, \cite{70})] Suppose that $X$ is an abelian
variety over $F$ for which the $\ell$-finiteness condition holds. Then for any $\Gamma$-submodule
$U$ in $V_{\ell}(X)$ there exists an element $u\in\End X\otimes{\bfQ}_{\ell}$, such that $u^2=u$
and $uV_{\ell}(X)=U$. Thus, the subspace $(1-u)V_{\ell}(X)$ is the $\Gamma$-invariant complement
of $U$ in $V_{\ell}(X)$, and the $\Gamma$-module $V_{\ell}(X)$ is semisimple.
\end{xxxi}

Applying this theorem to the graphs of homomorphisms of Tate modules, we obtain the following.

\begin{xxxi}[Homomorphism theorem \textup{(loc.~cit)}] Suppose that $X$ and $Y$ are abelian varieties over $F$, and their
product $X\times Y$ satisfies the $\ell$-finiteness condition. Then the natural maps
\begin{align*}
\Hom(X,Y)&\otimes{{\bfZ}_{\ell}}\hookrightarrow\Hom\nolimits_{\Gamma}(T_{\ell}(X),T_{\ell}(Y)),\\
\Hom(X,Y)&\otimes{{\bfQ}_{\ell}}\hookrightarrow\Hom\nolimits_{\Gamma}(V_{\ell}(X),V_{\ell}(Y))
\end{align*}
are bijections. In particular, the Galois modules $V_{\ell}(X)$ and $V_{\ell}(Y)$ are isomorphic
if and only if $X$ and $Y$ are $F$-isogenous. For example, if the $\ell$-finiteness condition
holds for the square $X^2=X\times X$, then the following natural maps are bijections:
\begin{align*}
\End X&\otimes{{\bfZ}_{\ell}}\hookrightarrow\End\nolimits_{\Gamma} T_{\ell}(X),\\
\End X&\otimes{{\bfQ}_{\ell}}\hookrightarrow\End\nolimits_{\Gamma} V_{\ell}(X)
\end{align*}
\end{xxxi}
\begin{xxxr}[Remark] When verifying the $\ell$-finiteness condition for $X\times Y$, it is sufficient to limit oneself
to $\ell$-divisible groups which are isomorphic to $X(\ell)$; the conclusion of the homomorphism
theorem still remains valid.
\end{xxxr}
\oripage
\begin{xxxi}[Principal polarization theorem (\cite{74}, \cite{52})]\footnote{\cite{52} contains a presentation of Deligne's proof.}
Suppose that $X$ is an abelian variety over $F$, and $X'$ is its Picard variety. Then there
exists a principal polarization on the abelian variety ${(X\times X')}^4$.
\end{xxxi}

\begin{xxxr}[Remark] Suppose that $X\rightarrow Y$ is an $F$-isogeny of abelian varieties (over $F$) whose kernel is
annihilated by multiplication by a natural number $r$. Then there exists an $F$-isogeny
$X'\rightarrow Y'$ whose kernel is also annihilated by multiplication by $r$. From this it easily
follows that there exists an $F$-isogeny of principally polarized abelian varieties ${(X\times
X')}^4\rightarrow{(Y\times Y')}^4$ whose kernel is annihilated by multiplication by $r$.
\end{xxxr}

Suppose that $X$ is an abelian variety over $F$, and $W$ is an $\ell$-divisible sub group of
$X(\ell)$ over $F$. For any natural number $n$ we let $W^\bot_n$ denote the orthogonal complement
of $W_n$ in ${X'}_n$ relative to the Weil pairing. Then the sequence $\{W^\bot_n\}$ gives an
$\ell$-divisible subgroup of $X'(\ell)$ over $F$. Here we have natural $F$-isomorphisms
$$
(X'/W^\bot_n)\simeq(X/W_n)'.
$$
This is easily seen to imply that there exists an $\ell$-divisible subgroup $H$ of ${(X\times
X')}^4$, over $F$ such that one has natural isomorphisms
$$
{(X\times X')}^4/H_n\simeq ((X/W_n)\times(X/W_n)')^4.
$$
\begin{xxxi}[Finiteness theorem for forms] Suppose that $X$ is an abelian variety over $F$, and $E$ is a finite
separable field extension of $F$. Then the set of abelian varieties $Y$ over $F$ \r(considered up
to $F$-isomorphism\r) for which there exists an $E$-isomorphism $X\otimes E\simeq Y\otimes E$ is
finite.
\end{xxxi}

This result follows immediately from a theorem of Borel and Serre on finiteness of the
noncommutative cohomology of groups of arithmetic type (\cite{35}, Proposition 3.8, p.\,136).

For the remainder of the section we shall consider the important concept of the
\textit{determinant of an $\ell$-divisible group}. If $W$ is any $\ell$-divisible group over a
field $F$ of height $h$, one has the free rank 1 ${\bfZ}_{\ell}$-module
$$
\det W=\bigwedge^h_{{\bfZ}_{\ell}}T(W)
$$
which is the maximal exterior power of the Tate module. The Galois group $\Gamma$ acts naturally
on $\det W$. The character $\Gamma\rightarrow\Aut\det W={\bfZ}^*_{\ell}$, which gives this
action, will be denoted by $\chi_W\colon \Gamma\rightarrow{\bfZ}^*_{\ell}$.

\begin{xxxr}[Examples]
\begin{enumerate}[\upshape(a)]
\item $W=\bfQ/\bfZ(\ell)$. Then $\chi_W$ is the trivial character.
\item $W=\bfQ/\bfG_m(\ell)$. Then $\chi_W=\chi_0$ is the \textit{cyclotomic} character, which gives the action of the Galois
group on the ${\ell}^n$th roots of 1 (for all $n$).
\item $W=X(\ell)$, where $X$ is a $g$-dimensional abelian variety over $F$. Then $\chi_W=\chi^g_0$ (this is a consequence of the existence of a
nondegenerate skewsymmetric Weil--Riemann pairing on the Tate module).
\end{enumerate}
\end{xxxr}

\oripage

If $W=W^1\times W^2$ is a product of $\ell$-divisible groups, then $\chi_W=\chi_{W^1}\chi_{W^2}$.

\section*{Digression. Proof of the principal polarization theorem}

This theorem follows from a somewhat more precise statement, given below.

Let $X$ be an abelian variety over a field $F$, and let $\alpha\colon X\rightarrow X'$ be a
polarization which is also defined over $F$. For any natural number $i$ one can determine an
abelian variety $X^i$ and a polarization $\alpha^i\colon X^i\rightarrow{(X')}^i={(X^i)}'$ such
that
$$
\Ker\alpha^i=(\Ker\alpha)^i\mtxt{and}\deg\alpha^i=(\deg\alpha)^i.
$$
\begin{xxxi}[Theorem] Fix a natural number $n$ such that the finite group scheme $\Ker\alpha$ is annihilated by
multiplication by $n$. \r(For example, $n$ can be taken to be $\deg\alpha$.\r) Then the following
assertions are true:
\begin{enumerate}[\upshape 1)]
\item If $-1$ is a square modulo $n$, then there exists a finite group subscheme $W_1$ in $\Ker\alpha^2\subset X^2$ such that\r:
\begin{enumerate}[\upshape a)]
\item $W_1\approx\Ker\alpha$ \r(over $F$\r) \r(in particular, the order of the group scheme $W_1$ is the square root of
$\deg\alpha^2$\r)\r;
\item $W_1$ is a \r(maximal\r) isotropic group subscheme of $\Ker\alpha^2$ relative to the Riemann form
corresponding to the polarization $\alpha^2$\r; and
\item there exists an $F$-isogeny $f_1\colon X^2\rightarrow X\times X'$ with kernel $W_1$. In particular, $\alpha^2$ descends
to a principal polarization on $X\times X'$ which is defined over $F$. By the same token,
$X\times X'$ is principally polarized over $F$.
\end{enumerate}
\item If $-1\mod n$ is a sum of two squares, then there exists a finite group subscheme $W_2$ in
$\Ker\alpha^4\subset X^4$ such that\r:
\begin{enumerate}[\upshape a)]
\item $W_2\approx\Ker\alpha^2$ \r(over $F$\r) \r(in particular, the order of the group scheme $W_2$ is the square root
of $\deg\alpha^4$\r)\r;
\item $W_2$ is a \r(maximal\r) isotropic group subscheme of $\Ker\alpha^4$ relative to the Riemann form
corresponding to the polarization $\alpha^2$\r; and
\item there exists an $F$-isogeny $f_2\colon X^4\rightarrow{(X\times X')}^2$ with kernel $W_2$. In particular,
$\alpha^4$ descends to a principal polarization on ${(X\times X')}^2$ which is defined
over $F$. By the same token, ${(X\times X')}^2$ is principally polarized over $F$.
\end{enumerate}
\item There always exists a finite group subscheme $W_4$ in $\Ker\alpha^8 \subset X^8$ such that\r:
\begin{enumerate}[\upshape a)]
\item $W_4\approx\Ker\alpha^4$ \r(over $F$\r) \r(in particular, the order of the group scheme $W_4$ is the square
root of $\deg\alpha^8$\r)\r;
\item $W_4$ is a \r(maximal\r) isotropic group subscheme of $\Ker\alpha^8$ relative to the Riemann form
corresponding to the polarization $\alpha^8$\r; and
\item there exists an $F$-isogeny $f_4\colon X^8\rightarrow{(X\times X')}^4$ with kernel $W_4$, In particular, $\alpha^8$ descends
to a principal polarization on ${(X\times X')}^4$ which is defined over $F$. By the same
token, ${(X\times X')}^4$ always has a principal polarization defined over $F$ \r(this is
the principal polarization theorem\r).
\end{enumerate}
\end{enumerate}
\end{xxxi}
\oripage
\begin{xxxr}[Proof] \textit{Construction of $W_1$, $W_2$, and $W_4$}.
\begin{enumerate}[\upshape1)]
\item We choose an integer $a\in\bfZ\subset\End X$ such that $a^2\equiv-1\mod n$. Then $W_1$ is the image of the group
scheme $\Ker\alpha\subset X$ under the imbedding $X\rightarrow X^2$, $y\mapsto(y,ay)(y\in
X)$.
\item We choose integers $a$, $b\in\bfZ$ such that $a^2+b^2\equiv-1\mod n$, and we consider the ``complex
number''
$$
I=\begin{pmatrix}a&{-b}\\b&a\end{pmatrix}\in M_2(\bfZ)\in M_2(\End X)=\End X^2.
$$
Then $W_2$ is the image of the group scheme $\Ker\alpha^2\subset X^2$ under the imbedding
$X^2\rightarrow X^2\times X^2=X^4$, $y\mapsto(y,Iy)(y\in X^2)$.
\item We choose a $4$-tuple of integers $a$, $b$, $c$, $d\in\bfZ$ such that $a^2+b^2+c^2+d^2\equiv-1\mod n$,
and we consider the ``quaternion''
$$
I=\begin{pmatrix}a&-b&-c&-d\\b&\phantom{-}a&\phantom{-}d&\phantom{-}c\\
c&-d&\phantom{-}a&\phantom{-}b\\d&\phantom{-}c&-b&\phantom{-}a\end{pmatrix}\in M_4(\bfZ)\in M_4(\End X)=\End X^4.
$$
Then $W_4$ is the image of the group scheme $\Ker\alpha^4\subset X^4$ under the imbedding
$X^4\rightarrow X^4\times X^4=X^8$, $y\mapsto(y,Iy)(y\in X^4)$.
\end{enumerate}
\end{xxxr}

\textit{Construction of the isogenies $f_i\colon X^{2i}\rightarrow{(X\times X')}^i$ with kernel
$W_i$}.
\begin{enumerate}[\upshape1)]
\item In the isomorphism
$$
f\colon X\times X\approx X\times X=X^2,\quad(y,z)\mapsto(y,ay)+(0,z)=(y,ay+z),
$$
we have $W_1=f(\Ker\alpha\times\{0\})\in f(X\times\{0\})$, and this enables us to define an
isogeny
$$
f_1=(\alpha,\id)f^{-1}\colon X^2\rightarrow X^2\rightarrow X'\times X=X\times X'
$$
with kernel $W_1$.
\item In the isomorphism
$$
f\colon X^2\times X^2\approx X^2\times X^2=X^4,\quad(y,z)\mapsto(y,Iy)+(0,z)=(y,Iy+z),
$$
we have $W_1=f(\Ker\alpha^2\times\{0\})\in f(X^2\times\{0\})$, and this enables us to define
an isogeny
$$
f_2=(\alpha^2,\id)f^{-1}\colon X^4\rightarrow X^4\rightarrow(X^2)'\times X^2={(X\times X')}^2
$$
with kernel $W_2$.
\item In the isomorphism
$$
f\colon X^4\times X^4\rightarrow X^4\times X^4=X^8,\quad(y,z)\mapsto(y,Iy)+(0,z)=(y,Iy+z),
$$
we have $W_4=f(\Ker\alpha\times\{0\})\in f(X\times\{0\})$, and this enables us to define an
isogeny
$$
f_4=(\alpha^4,\id)f^{-1}\colon X^8\rightarrow X^8\rightarrow(X^4)'\times X^4={(X\times X')}^4
$$
with kernel $W_4$.
\end{enumerate}

\oripage

\section{Isogenies of abelian varieties over local fields, and Galois modules}

\addtocounter{ssectc}{-1}

\ssect In this section $\hat K$ is a complete discrete valuation field of characteristic zero
with algebraically closed residue field $k$ of finite characteristic $p$, $v\colon {\hat
K}^*\rightarrow\bfZ$ is the discrete valuation of $K$, normalized by the condition $v({\hat
K}^*)=\bfZ$, and $\goto$ is the valuation ring, i.\,e.,
$$
\goto=\{x\in{\hat K}|v(x)\geq0\}.
$$
The number $e=v(p)$ is called the \textit{ramification index} of the field $\hat K$; if $e=1$,
then $\hat K$ is said to be \textit{unramified}. We fix an element $\pi\in\goto\subset\hat K$
with $v(\pi)=1$. (Concerning the concepts in this section, see~\cite{2}, \cite{66},
and~\cite{Ser3}.)

The Galois group $I=\Gal(\ol{\hat K}/\hat K)$ of the field $\hat K$ coincides with the inertia
group. If $n$ is a natural number prime to $p$, then all of the $n$th roots of 1 lie in $\hat K$;
they form a cyclic group of order $n$, which we denote by $\mu_n(\hat K)$. We have canonical
surjective homomorphisms
$$
j_n\colon I\rightarrow{\bfmu}_n(\hat K),
$$
whose kernels correspond to the tamely ramified cyclic extensions ${\hat K}_n=\hat K(\pi^{1/n})$
of degree $n$ over $\hat K$ (which do not depend on the choice of $\pi$). Here
$$
\sigma(\pi^{1/n})=j_n(\sigma)\pi^{1/n}\mtxt{for all}\sigma\in I.
$$
The map $\goto\twoheadrightarrow k$ of taking residues gives a canonical isomorphism
$$
{\bfmu}_n(\hat K)\stackrel\sim{\rightarrow}\bfmu_n(k)
$$
between the groups of $n$th roots of 1 in $\hat K$ and in $k$. Hence, the composition
$$
I\stackrel{j_n}{\rightarrow}{\bfmu}_n(\hat K)\stackrel\sim{\rightarrow}\bfmu_n(k)
$$
gives a surjection $I\twoheadrightarrow{\bfmu}_n(k)$, which we shall also denote by $j_n$.

We consider the case $n=p-1$ in more detail. Clearly ${\bfmu}_{p-1}(k)$ is the group of
invertible elements of the prime subfield $\bfZ/p\bfZ\subset k$. We let $\tau_p$ denote the
corresponding homomorphism
$$
j_{p-1}\colon \Gamma\rightarrow{\bfmu}_{p-1}(k)=(\bfZ/p\bfZ)^*.
$$
Any other (continuous) character $\chi\colon I\rightarrow(\bfZ/p\bfZ)^*$ is equal to a power
$\tau^n_p$, where the exponent $n$ is an integer which can always be taken between 0 and $p-1$.
Thus, for example, the cyclotomic character
$$
{\ol\chi}_0\colon I\rightarrow(\bfZ/p\bfZ)^*
$$
which gives the action of the Galois group on the group ${\bfmu}_p(\ol{\hat K})$ of $p$-th roots
of 1, is equal to $\tau^e_p$. In particular, if the field $\hat K$ is unramified, then
$$
\tau_p=\ol{\chi_0}.
$$

We note that $\ol{\chi_0}\equiv\chi_0\mod p$, where $\chi_0\colon \Gamma\rightarrow{\bfZ}^*_p$ is
the cyclotomic character, which gives the action of the Galois group on the $p^m$th roots of 1
(for all $m$).

\oripage

Suppose that $r$ is a natural number, $q=p^r$ and $F={\bfF}_q$, is a finite field of $q$
elements. Any imbedding of fields $i\colon F\hookrightarrow k$ gives an isomorphism
$F^*\stackrel\sim{\rightarrow}{\bfmu}_{q-1}(k)$, which induces a surjection
$$
\chi_i\colon I\stackrel{j_{q-1}}{\rightarrow}{\bfmu}_{q-1}(k)\simeq F^*,
$$
from which the imbedding can be uniquely recovered. We let $M$ denote this set of surjections
$\chi_i\colon I\rightarrow F^*$. This set (like the set of all imbeddings $F\hookrightarrow k$)
is a principal homogeneous space over the finite cyclic group $\bfZ/r\bfZ$. Here an arbitrary
element $\ell\mod r\bfZ$ of the group $\bfZ/r\bfZ$ takes the character $\chi_i$ to the character
${\chi}_{i+\ell}=\chi_i^{p^{\ell}}$. In particular, if we fix an imbedding $i\colon
F\hookrightarrow k$, then the entire set $M$ consists of the $r$ characters
$$
\chi_i,\chi_{i+1},\dots,\chi_{i+r}(=\chi_i)
$$
of the group $I$ with values in the multiplicative group $F^*$ of the field $F$. Any (continuous)
character $\chi\colon I\rightarrow F^*$ can be uniquely represented in the form of a product
$$
\prod_{1\leq{\ell}\leq r}\chi^{n_{\ell}}_{i+\ell},\quad0\leq n_{\ell}\leq p-1,
$$
with the exception of the identity character, which has two representations
$$
\prod\chi^0_{i+{\ell}}=\prod\chi^{p-1}_{i+{\ell}}.
$$
Conversely, any $r$-tuple $n_1,\dots,n_r$ of nonnegative integers not exceeding $p-1$ gives a
character
$$
\chi_{n_1,\dots,n_r}=\prod\chi^{n_{\ell}}_{i+{\ell}}\colon I\rightarrow F^*.
$$
\begin{xxxr}[Example] If $r=1$, then $F=\bfZ/p\bfZ$, and there is exactly one imbedding $i\colon \bfZ/p\bfZ\hookrightarrow k$. Here
$\chi_i=\tau_p$, ${\ol\chi}_0=\chi^e_i$.
\end{xxxr}

\ssect Suppose that $F$ is a finite field of $q=p^r$ elements, and $V$ is a
\textit{one-dimensional} vector space over $F$. Any character $\chi\colon I\rightarrow F^*$ gives
an $I$-module structure, i.\,e., a Galois module structure, on $V$. Regarding $V$ as a vector
space over the prime field $\bfZ/p\bfZ$, we see that it is $r$-dimensional and is equipped with a
Galois module structure, which we denote $V(\chi)$. It is easy to verify that the Galois module
$V(\chi)$ is simple if and only if one of the following conditions is fulfilled: (a) $F$ is a
prime field, i.\,e., $r=1$; or (b) the sequence of numbers $(n_1,\dots,n_r)$ giving the character
$\chi=\chi_{n_1,\dots,n_r}$ is not invariant under any power of the cyclic permutation.

\begin{xxxi}[Serre's lemma \cite{Ser3}] Suppose that $V$ is a finite abelian group which is annihilated by
multiplication by $p$, i.\,e., it is a vector space over $\bfZ/p\bfZ$ of finite dimension $r$.
Further suppose that a simple $I$-module structure is given on $V$. Then the centralizer
$F={\End}_IV$ is a finite field of $q=p^r$ elements, and $V$ is a one-dimensional vector space
over $F$. The $I$-module structure on $V$ is given by a character
$$
\chi\colon I\rightarrow F^*\subset\Aut(V),
$$

\oripage i.\,e., $V\simeq V(\chi)$.
\end{xxxi}

\begin{xxxr}[Definition] A simple Galois module $V\simeq V(\chi)$ will be said to be \textit{admissible} if the defining
character $\chi$ can be represented in the form
$$
\chi=\prod_{1\leq{\ell}\leq r}\chi^{n_{\ell}}_{i+\ell},\quad0\leq n_{\ell}\leq e\mtxt{for all}\ell.
$$
\end{xxxr}

\begin{xxxr}[Examples] The group $\bfZ/p\bfZ$ (with trivial action of the Galois group) and the Galois module
${\bfmu}_p(\ol{\hat K})$ are always admissible. Conversely, if $\hat K$ is unramified, then any
1-dimensional admissible module (over $\bfZ/p\bfZ$) is isomorphic either to $\bfZ/p\bfZ$ or to
${\bfmu}_p(\ol{\hat K})$. If $e\geq p-1$, then all simple Galois modules are admissible. If
$e=p-1$, then the Galois modules $\bfZ/p\bfZ$ and ${\bfmu}_p(\ol{\hat K})$ are isomorphic.

Suppose that $V$ is a finite abelian group which is annihilated by multiplication by $p$, i.\,e.,
it is a vector space over $\bfZ/p\bfZ$ of some finite dimension $m$. Further suppose that a (not
necessarily simple) Galois module structure is given on $V$. We let $\det V$ denote the maximal
exterior power $\bigwedge^mV$ of the vector space $V$. The one-dimensional space $\det V$ is a
Galois module, and its structure is given by a character $I\rightarrow\Aut(\det
V)={(\bfZ/p\bfZ)}^*$, which we denote $\chi_v\colon I\rightarrow{(\bfZ/p\bfZ)}^*$. For example,
if $V\simeq V(\chi)$, then $\chi_v=\tau^{n_1+\dots+n_r}_p$ (here $\chi=\prod_{1\leq{\ell}\leq
r}\chi^{n_{\ell}}_{i+{\ell}}$); in particular, if $\hat K$ is unramified, then
$$
\chi V=\ol X^{n_1+\dots+n_r}_0.
$$
\end{xxxr}
\ssect If $M$ is a finite $\goto$-module, we let $\calL(M)$ denote its length. If $\bfG$ is a
(quasi)finite commutative flat group scheme over $\Spec\goto$, we set
$$
\dim(\bfG)=\calL(e^*\Omega^1_{\bfG/\Spec\goto}),
$$
where $e\colon \Spec\goto\hookrightarrow\bfG$ is the identity section. If the group scheme $\bfG$
is annihilated by multiplication by some natural number $m$ prime to $p$, then it is \'etale, and
$\Omega^1_{\bfG/\Spec\goto}=\{0\}$; in particular, $\dim(\bfG)=0$.

\sssect Suppose that $\bfG$ is annihilated by multiplication by $p$. Then the group
$\bfG(\ol{\hat K})$ of its $\ol{\hat K}$-points is a vector space over $\bfZ/p\bfZ$ of some
finite dimension $h$, and it is equipped with a natural Galois module structure. This Galois
module uniquely determines (up to isomorphism) the generic fiber ${\bfG}_{\hat K}$ of the group
scheme $\bfG$, which is a group scheme over $\Spec\hat K$. We let $\det\bfG\colon
I\rightarrow{(\bfZ/p\bfZ)}^*$ denote the character $\chi_{\bfG(\ol{\hat K})}$.

\begin{xxxi}[Raynaud's lemma \cite{60}]
\begin{enumerate}[\upshape(a)]
\item Suppose that $\bfG$ is a finite flat commutative group scheme over $\Spec\goto$ which is annihilated by
multiplication by $p$. If the Galois module $\bfG(\hat K)$ is simple, then it is admissible.
\item Suppose that $V(\chi)$ is an admissible Galois module, corresponding to the character $\chi=\chi_{n_1,\dots,n_r}$.
Then there exists a finite flat commutative group scheme $\bfG(\chi)$ over $\Spec\goto$ which
satisfies the following conditions
\end{enumerate}
\oripage
\begin{enumerate}[\upshape1)]
\item The Galois module $\bfG(\ol{\hat K})$ is isomorphic to $V(\chi)$; in particular,
$\det\bfG=\tau^{n_1+\dots+n_r}_p$.
\item As an affine scheme, $\bfG(\chi)$ is isomorphic to the spectrum of the ring
$$
\goto[t_1,\dots,t_r]/\langle t^p_1-\pi^{n_1}t_2,t^p_2-\pi^{n_2}t_3,\dots,t^p_r-\pi^{n_r}t_1\rangle;
$$
for the identity section e one can take the point with zero values of the parameters
$t_1,\dots,t_r$. Thus,
$$
e^*\Omega^1_{\bfG(\chi)/\goto}\simeq\goto/\pi^{n_1}\goto\oplus\dots\oplus\goto/\pi^{n_r}\goto,
$$
$$
\dim\bfG(\chi)=n_1+\dots+n_r\mtxt{and}\det\bfG(\chi)=\tau^{\dim\bfG(\chi)}_p.
$$
In particular, if $\hat K$ is unramified, then $\det\bfG(\chi)={\ol\chi}^{\dim\bfG(\chi)}_0$.
\item If $e<p-1$, then for any admissible module $V(\chi)$ there exists exactly one finite flat
commutative group scheme $\bfG$ over $\Spec\goto$ which is annihilated by multiplication by
$p$ and is such that the Galois modules $\bfG(\ol{\hat K})$ and $V(\chi)$ are isomorphic, in
particular, $\bfG\cong\bfG(\chi)$ and $\det\bfG=\tau^{\dim\bfG}_p$.
\end{enumerate}
\end{xxxi}

\begin{xxxr}[Examples]
(a) If $\chi$ is the trivial character, then $\bfG(\chi)$ is the constant group scheme
$\bfZ/p\bfZ$. (b) If $\chi=\ol\chi_0$ is the cyclotomic character, then $\bfG(\chi)$ is the
scheme ${\bfmu}_p$ of $p$th roots of 1. (c) If the character $\chi=\tau^n_p$ corresponds to a
prime field $F$, i.\,e., if $r=1$ (see \S3.1), then, in the notation of Tate and Oort~\cite{28a},
$\bfG(\chi)$ is the group scheme $G_{a,b}$ of order $p$ with $a=\pi^n$ and $b=\pi^{e-n}$.
\end{xxxr}

\begin{xxxr}[Remark \cite{60}] Raynaud's lemma is also valid when $e=p-1$, except that there are two group
schemes corresponding to the trivial character $(\chi=\ol{\chi_0})$, namely, $\bfZ/p\bfZ$ and
${\bfmu}_p$.

From this it easily follows that, if $\hat K$ is unramified, and if the Galois module
$\bfG(\ol{\hat K})$ is simple, then $\det(\bfG)={\ol\chi}^{\dim(\bfG)}_0$.
\end{xxxr}

\ssect Suppose that $f\colon X\rightarrow Y$ is an isogeny of abelian varieties over $\hat K$
having semistable reduction, $\tilde f\colon \bfX\rightarrow\bfY$ is the corresponding morphism
of their N\'eron models, and $\bfG=\Ker\tilde f$ is its kernel, which is a quasifinite
commutative flat group scheme over $\Spec\goto$. The kernel $\Ker f$ of the isogeny $f$ is a
finite Galois module which coincides with $\bfG(\ol{\hat K})$ and has order equal to the degree
$\deg f$ of the isogeny $f$. If $\Ker f$ is annihilated by multiplication by some natural number
$n$, then $\Ker f\subset X_n$, and we have the commutative diagram
$$
\begin{array}{rcl}
f\colon X\ &\!\longrightarrow\!\!&\ \ \;Y\\
{}_n\!\!\!\searrow&&\swarrow\\
&\!X&
\end{array}
$$
where all the arrows are isogenies of abelian varieties. Passing to morphisms of N\'eron models
$$
\begin{array}{rcl}
\tilde f\colon \bfX\ &\!\longrightarrow\!\!&\ \ \;\bfY\\
{}_n\!\!\!\searrow&&\swarrow\\
&\!\bfX&
\end{array}
$$
\oripage we find that $\bfG$ is annihilated by multiplication by $n$; in particular, $\bfG$ is
annihilated by multiplication by $\deg f$. Thus, if $\deg f$ is prime to $p$, it follows that
$\dim(\bfG)=0$.

\sssect Suppose that $f\colon X\rightarrow Y$ is a $p$-isogeny, i.\,e., $\deg f$ is a power of
some prime $p$. We note that any such isogeny factors into a composition of isogenies over $\hat
K$ whose kernels are annihilated by multiplication by $p$. We next suppose that $\Ker f$, and
hence $\bfG$, are annihilated by multiplication by $p$. Then the isogeny $f$ factors into a
product of isogenies over $\hat K$ whose kernels are finite simple Galois modules. The next
assertion follows easily from results of Raynaud (\S4.2.1).

\sssect \begin{xxxi}[Lemma] Let $f\colon X\rightarrow Y$ be an isogeny of abelian varieties over
$K$ whose kernel $\Ker f$ is annihilated by multiplication by $p$. If $X$ has good reduction over
$K$, then $\bfG$ is a finite flat group scheme which is annihilated by multiplication by $p$. If,
in addition, $\Ker f$ is a simple Galois module, then it is admissible, and if $\hat K$ is
unramified, then $\det\bfG={\ol\chi}^{\dim(\bfG)}_0$.

Conversely, given any admissible Galois module $V$, there exists an isogeny of abelian varieties
with good reduction over $\hat K$ whose kernel is isomorphic to $V$ \r(a special case of
Raynaud's theorem~\cite{33}\r).
\end{xxxi}

\begin{xxxr}[Remark] The statement that the kernel of the isogeny is admissible remains valid if we replace
good reduction by semistable reduction: if $\bfG$ is not finite, and $\Ker f$ is a simple Galois
module, then it is trivial (\cite{71}, p.\,220, Lemma 4.2.4).
\end{xxxr}
\begin{xxxr}[Remark] Given any finite Galois module $V$, there exists an isogeny of abelian varieties over
$\hat K$ whose kernel is isomorphic to $V$ (see~\cite{33}).

\sssect
\begin{xxxi}[Corollary] If $e<p-1$, and if $X$ has good reduction, then the dimension $\dim(\bfG)$ can be uniquely
recovered from the Galois module $\Ker f$. Here $f\colon X\rightarrow Y$ is an arbitrary isogeny
over $\hat K$, and $\bfG$ is the kernel of the corresponding morphism of N\'eron models.
\end{xxxi}

\begin{xxxr}[Proof] It is sufficient to consider the case when $f$ is a $p$-isogeny, and even an isogeny whose
kernel is annihilated by multiplication by $p$. We factor $f$ into a composition of isogenies
$$
X=Y_{(0)}\stackrel{f_1}{\rightarrow}Y_{(1)}\stackrel{f_2}{\rightarrow}\dots\stackrel{f_m}{\rightarrow}Y_{(m)}=Y,
$$
where all of the kernels $V_i=\Ker f_i$ are simple Galois modules. Then the set $V_1,\dots,V_m$
coincides with the set of simple Jordan--H\"older quotients of the Galois module $\Ker f$, and it
can be uniquely (up to rearrangement) recovered from $\Ker f$. We know that each group scheme
${\bfG}_i/\Spec\goto$, which is the kernel of the homomorphism of N\'eron models corresponding to
$f_i$ can be uniquely (up to isomorphism) recovered from $V_i$.
\end{xxxr}
We choose a nonzero global differential form $\omega$ on $Y$ of highest dimension. Its inverse
image on $Y_{(i)}$ under the natural isogeny $Y_{(i)}\rightarrow Y$ will be denoted
\end{xxxr}
\oripage $\omega_i$. We have
$$
\omega_m=\omega,\quad\omega_0=f^*\omega,\quad\omega_{i-1}=f^*_i\omega_i.
$$
From the properties of the local term in the logarithmic height (\S2.1) it follows that
\begin{align*}
\dim(\bfG)&=\calL(X_{(0)},\omega_0)-\calL(Y_{(m)},\omega_m),\\
\dim(\bfG_i)&=\calL(Y_{(i-1)},\omega_i)-\calL(Y_{(i)},\omega_i).
\end{align*}
Hence,
$$
\dim(\bfG)=\calL(Y_{(0)},\omega_0)-\calL(Y_{(m)},\omega_m)=\sum^m_{i=1}
[\calL(Y_{(i-1)},\omega_i)-\calL(Y_{(i)},\omega_i)]=\sum^m_{i=1}\dim({\bfG}_i).
$$
By the same token, the number
$$
\dim(\bfG)=\calL(X,f^*\omega)-\calL(Y,\omega)
$$
can be uniquely recovered from the numbers $\dim({\bfG}_i)$ depending only on the ${\bfG}_i$,
which, in turn, depend only on the $V_i$; this proves Corollary 4.3.3.

\sssect \begin{xxxi}[Corollary] Under the assumptions of Corollary \r{4.3.3}, if $\hat K$ is
unramified, then
$$
\chi_{\Ker f}={\ol\chi}^{\dim(\bfG)}_0={\ol\chi}^{\calL(X,f^*\omega)-\calL(Y,\omega)}_0.
$$
\end{xxxi}

The proof is obtained by a straightforward application of Lemma 3.2 and the following equality
(in the notation of the proof of Corollary 4.3.3):
$$
\chi_{\Ker f}=\chi_{V_1}\chi_{V_2}\dots\chi_{V_m}\colon I\rightarrow{\bfZ/p\bfZ}^*.
$$

\sssect \begin{xxxr}[Example] Let $n\colon X\rightarrow X$ be the isogeny multiplication by $n$
on an arbitrary abelian variety $X$ over $\hat K$. If $\omega$ is any nonzero form of highest
dimension on $X$, then $n^*\omega=n^{\dim X}\omega$. Hence,
$$
\calL(X,n^*\omega)-\calL(X,\omega)=v(n^{\dim X})=(\dim X)v(n).
$$
For example, if $n=p^m$, then $\calL(X,p^{m^*}\omega)-\calL(X,\omega)=\text{me}$.
\end{xxxr}

\sssect We now suppose that $X$ and $Y$ are abelian varieties with good reduction over $\hat K$,
and $e<p-1$. Let $\varphi\colon X_n\rightarrow Y_n$ be a homomorphism of Galois modules, and let
$$
\Gamma_\varphi=\{(x,\varphi x)|x\in X_n\}\subset X_n\times Y_n={(X\times Y)}_n\subset X\times Y
$$
be its graph, which is isomorphic as a Galois module to the group $X_n$. Then $\Gamma_\varphi$ is
the kernel of the natural isogeny
$$
\psi\colon X\times Y\rightarrow Z=(X\times Y)/\Gamma_\varphi,
$$
which is isomorphic to the kernel $X_n$ of the isogeny $n\colon X\rightarrow X$. Hence, if
$\omega'$ is any form of highest dimension on the abelian variety $Z$, we have
$$
\calL(X\times Y,\psi^*\omega')-\calL(Z,\omega')=\calL(X,n^*\omega)-\calL(X,\omega)=v(n).
$$
\oripage

\sssect \begin{xxxi}[Lemma] Let $f\colon X\rightarrow Y$ be an isogeny of abelian varieties over
$\hat K$ with semistable reduction whose kernel $\Ker f$ is annihilated by multiplication by $n$.
Then
$$
0\leq\calL(X,f^*\omega)-\calL(X,\omega)\leq(\dim X)v(n).
$$
\end{xxxi}

\begin{xxxr}[Proof] We recall that
$$
\calL(X,f^*\omega)-\calL(X,\omega)=\dim(\bfG)\geq0,
$$
where $\bfG$ is the kernel of the homomorphism of N\'eron models corresponding to $f$. Let
$g\colon Y\rightarrow X$ be the isogeny such that the composition $gf\colon X\rightarrow X$ is
multiplication by $n$. Let $\omega'$ be the form of highest weight on $X$ such that
$g^*\omega'=\omega$. Then
\begin{align*}
\calL(X,f^*\omega)&-\calL(Y,\omega)\\
&\leq\calL(X,f^*\omega)-\calL(Y,\omega)+[\calL(Y,g^*\omega')-\calL(X,\omega')]\\
&=\calL(X,{(gf)}^*\omega')-\calL(Y,\omega')=\calL(X,n^*\omega')-\calL(X,\omega')\\
&=(\dim X)v(n).
\end{align*}

\ssect Let $X$ be an abelian variety over $\hat K$, and let $W$ be a $p$-divisible group over
$\hat K$ of height $h$ which is isomorphic to some $p$-divisible subgroup of $X(p)$ over $\hat
K$. The theory of Hodge--Tate moduli (see~\cite{28}, \cite{21}, \cite{63}, and~\cite{65}) enables
one to determine a nonnegative number $d=\dim W$, which is called the \textit{dimension} of the
$p$-divisible group $W$ (over $\hat K$) and possesses the following properties.

\begin{enumerate}[\upshape(a)]
\item $\chi_W=\epsilon\chi^d_0$ where $\epsilon\colon I\rightarrow{\bfZ}^*_p$ is a character of finite order. This property may be regarded as
the definition of the dimension $d$. (If $X$ has semistable reduction over $\hat K$, then
$\epsilon=1$ and $\chi_W=\chi^d_0$.) In particular, $p$-divisible groups which are isomorphic
over $\hat K$ have the same dimension.
\item If $W\simeq\bfQ/\bfZ(p)$, then $d=0$; if $W\simeq{\bfG}_m(p)$, then $d=1$; and if $W\simeq Y(p)$ for some abelian
variety $Y$ over $\hat K$, then $d=\dim Y$. In particular,
$$
\dim X(p)=\dim X.
$$
\item $\dim W\leq h$.
\item If $W'$ is a $p$-divisible subgroup (over $\hat K$) in $W$, then $\dim W'\leq \dim W$. In particular, $\dim W\leq\dim X$.
\item If $W$ is the generic fiber of a $p$-divisible group $H=\{H_n\}$ over $\goto$ which is made up of
finite group schemes $H_n$ over $\Spec\goto$ (Tate~\cite{28}), then $W$ and $H$ have the same
dimension. Here $\dim(H_n)=ne\dim(W)$.
\item Consider the natural $p$-adic representation $\rho\colon I\rightarrow\Aut T_p(W)$ of the Galois group in the
Tate module $T(W)$. Its image $\Im\rho$ is a compact $p$-adic Lie group. This group is finite
if and only if $d=0$. If $d=1$, and if $W$ has no nontrivial (i.\,e., different from $\{0\}$
and $W$) $p$-divisible subgroups over $\hat K$ then the centralizer $E={\End}_IV(W)$ is a
$p$-adic field, and $\Im\rho$ is an open subgroup of the group of $E$-linear automorphisms of
the Tate ${\bfQ}_p$-module $T(W)\otimes_{{\bfZ}_p}{\bfQ}_p=V(W)$. \oripage
\item If $W$ is isomorphic to a product of $p$-divisible groups $W'\times W''$ (over $\hat
K$), then $\dim W=\dim W'+\dim W''$.
\end{enumerate}

Suppose that $X$ has semistable reduction over $\hat K$. We define the sequence of abelian
varieties $Y_n=X/W_n$ over $\hat K$, and we consider the isogenies
\begin{align*}
\varphi_n\colon Y_{(n)}&\rightarrow Y_{(n+1)}=Y_{(n)}/(W_{n+1}/W_n),\\
\varphi_{n,m}\colon Y_{(n)}&\rightarrow Y_{(n+m)}=Y_{(n)}/(W_{n+m}/W_n).
\end{align*}

We note that $\Ker{\varphi_n}$ is isomorphic as a Galois module to $W_1$ for all $n$; and
$\Ker{\varphi_{(n,m)}}$ is isomorphic to $W_m$ for any $n$ and $m$. The following fact
(\cite{41}, pp.\,134--136, \S3.4) follows from results of Tate and Raynaud~\cite{28},
\cite{60} concerning $p$-divisible groups and results of Grothendieck and Raynaud~\cite{46}
concerning N\'eron models.
\end{xxxr}

\textit{There exists a natural number $N$ which satisfies the following conditions. For any
$n\geq N$ let ${\bfG}_n$ denote the quasifinite flat group scheme over $\Spec\goto$ which is the
kernel of the homomorphism of N\'eron models corresponding to $\varphi_n$, and let ${\bfG}_n$
denote the group scheme corresponding to the isogeny $\varphi_{n,m}$. Then
$$
e^*\Omega^1_{\bfG_n/\Spec\goto}\simeq{(\goto/p\goto)}^{\dim W},\quad\dim({\bfG}_n)=v(p)\dim W=e\dim W.
$$
Thus, for any nonzero invariant forms $\omega_n$ on $Y_n$ of highest dimension such that
$\varphi^*_n\omega_{n+1}=\omega_n,\ \varphi^*_{n,m}\omega_{n+m+1}=\omega_n$, we have
$$
\calL(Y_{(n)},\omega_n)-\calL(Y_{(n+1)},\omega_n)=\dim({\bfG}_n)=e\dim W,
$$
$$
\dim({\bfG}_{n,m})=\calL(Y_{(n)},\omega_n)-\calL(Y_{(n+m)},\omega_{n+m})=ne\dim W.
$$
}

The proof is based on the construction of $p$-divisible group $H=\{H_n\}$ over $\goto$ such that
$\dim H=\dim W$ and $\dim(H_m)=\dim({\bfG}_{N+m})$.

\ssect Let $X$ be $g$-dimensional abelian variety over $\hat K$, and let $\ell$ be a prime $\neq
p$. If $X$ has semistable reduction, then the inertia group $I$ acts on the Tate module
$T_{\ell}(X)$ by unipotent matrices~\cite{24}.

\sssect Now suppose that $X$ does not necessarily have semistable reduction over $\hat K$. If
$K'$ is a finite Galois extension of $\hat K$ over which all of the points of order $n$ on $X$
are rational, then $X_{K'}=X\otimes K'$ has semistable reduction over $K'$ (Raynaud's criterion
\cite{45}). Here $n$ is a natural number $\geq3$ and prime to $p$. For $K'$ we shall take the
field of definition of all of the points of order $n$, i.\,e., $K'$ corresponds to the kernel of
the homomorphism
$$
\rho_n\colon I\rightarrow\Aut(X_n)\simeq\GL(2g,\bfZ/n\bfZ).
$$
We let $\ol c(n)$ denote the least common multiple of the orders of the cyclic subgroups of
$\GL(2g,\bfZ/n\bfZ)$. For any element $\sigma\in I$ the power $\sigma^{\ol c(n)}$ lies in the
Galois group $I'$ of $K'$, and hence it acts unipotently on the Tate module
$T_{\ell}(X_{K'})=T_{\ell}(X)$, and also on the space $X_{\ell}=T_{\ell}(X)/\ell T_{\ell}(X)$
(over the field $\bfZ/\ell\bfZ$).

\oripage

We should say a few words about the choice of $n$. If $p\neq3$, then we may set $n=3$; if $p=3$,
then $n=5$. Thus, always either $\sigma^{\ol c(3)}$ or $\sigma^{\ol c(5)}$ is a unipotent
automorphism of the Tate module. Since $\ol c(15)=\ol c(3)\ol c(5)$, we conclude that the
automorphism $\sigma^{\ol c(15)}$ is always unipotent, i.\,e., all of the eigenvalues of $\sigma$
acting on $T_{\ell}(X)$, and hence on the $\bfZ/\ell\bfZ$-space $X_{\ell}=T_{\ell}(X)/(\ell)$,
are $\ol c(15)$th roots of unity. Let $c(g)=\ol c(15)$ be the least common multiple of the orders
of the cyclic subgroups of $\GL(2g,\bfZ/15\bfZ)$.

\sssect Let $E$ be a free Galois ${\bfZ}_{\ell}$-submodule of rank $h$ in $T_{\ell}(X)$. Then the
eigenvalues of $\sigma$ acting on $E\subset T_{\ell}(X)$ are $c(g)$th roots of 1. Hence, the
image of the character
$$
\chi_E\colon I\rightarrow\Aut(\bigwedge^hE)={\bfZ}^*_{\ell}
$$
(acting on $\det E=\bigwedge^hE$) consists of $c(g)$th roots of 1, i.\,e., $\chi^{c(g)}_E=1$.

\sssect Let $V$ be a Galois submodule of $X_{\ell}$. Then the eigenvalues of $\sigma$ (for any
$\sigma\in I$) acting on $X_{\ell}$, and hence on $V$, are $c(g)$th roots of 1. Just as above,
the character
$$
\chi_V\colon I\rightarrow\Aut(\det V)={\bfZ/\ell\bfZ}^*
$$
satisfies $\chi^{c(g)}_V=1$.

\section{Behavior of the height under isogenies of abelian varieties over an algebraic number field}

\addtocounter{ssectc}{-1}

\ssect \begin{xxxi}[Theorem] Let $X$ be a $g$-dimensional abelian variety over the rational
number field $\bfQ$ which has good reduction at the prime $\ell$. Let $p$ be a prime $\neq\ell$
which is greater than
$$
C(\ell,g)={({2\ell}^{gc(g)})}^{c(g)\bigl(\begin{smallmatrix}2g\\g\end{smallmatrix}\bigr)},
$$
where $c(g)$ is the least common multiple of the orders of all cyclic subgroups of
$\GL(2g,\bfZ/15\bfZ)$. Let $f\colon X\rightarrow Y$ be an isogeny of abelian varieties over
$\bfQ$ whose kernel $\Ker f$ is annihilated by multiplication by $p$. If $X$ has good reduction
at $p$, then for any nonzero form $\omega$ of highest dimension on $Y$
$$
\calL(X\otimes{\bfQ}_p,f^*\omega)-\calL(Y\otimes{\bfQ}_p,\omega)=\frac12{\dim}_{\bfZ/p\bfZ}\Ker f=\frac12\ord_p(\deg f).
$$
Here ${\ord}_p\colon {\bfQ}^*\rightarrow\bfZ$ is the standard $p$-adic valuation of $\bfQ$.
\end{xxxi}
\sssect \begin{xxxi}[Corollary] Under the assumptions of Theorem 5.0, for any $p$-isogeny
$\pi\colon X\rightarrow Y$ of abelian varieties over $\bfQ$,
$$
\calL(X\otimes{\bfQ}_p,\pi^*\omega)-\calL(Y\otimes{\bfQ}_p,\omega)=\frac12{\ord}_p(\deg\pi).
$$
\end{xxxi}

In fact, any abelian variety isogenous to $X$ satisfies the conditions of Theorem 5.0. It remains
to recall that any $p$-isogeny factors into a composition of $p$-isogenies whose kernels are
annihilated by multiplication by $p$.

\oripage \sssect \begin{xxxi}[Corollary] Let $\varphi\colon X\rightarrow Y$ be an isogeny of
$g$-dimensional abelian varieties over $\bfQ$ which have good reduction at $\ell$. Suppose that
the following holds for all prime divisors $p$ of $\deg\varphi\colon p>C(\ell,g)$ and $X$ has
good reduction at $p$. Then the abelian varieties $X$ and $Y$ have the same canonical height.
\end{xxxi}

In fact, for all primes $t$ not dividing $\deg\varphi$,
$$
\calL(X\otimes{\bfQ}_t,\varphi^*\omega)=\calL(Y\otimes{\bfQ}_t,\omega).
$$
We factor $\varphi$ into a composition of isogenies
$$
X=Y_{(0)}\stackrel{\varphi_1}{\rightarrow}Y_{(1)}\stackrel{\varphi_2}{\rightarrow}\dots Y_{(r-1)}\rightarrow Y_{(r)}=Y,
$$
where $\deg\varphi_i$, is the highest power of the prime $p_i$ dividing the degree of $\varphi$;
here $p_1,\dots,p_r$ are all of the prime divisors of $\deg\varphi$. Then, for any of the primes
$p=p_i$,
$$
\calL(X\otimes{\bfQ}_p,\varphi^*\omega)-\calL(Y\otimes{\bfQ}_p,\omega)=
\calL(Y_{(i-1)}\otimes{\bfQ}_p,\varphi^*_i\omega_{i-1})-\calL(Y_{(i)}\otimes{\bfQ}_p,\omega_i)
$$
(here $\omega_i$ is the form on $Y_{(i)}$ which is the inverse image of the form $\omega$ on $Y$
relative to the isogeny $Y_{(i)}\rightarrow Y$, and, by Corollary 5.0.1 applied to the isogeny
$\varphi_i$,
$$
\calL(X\otimes{\bfQ}_p,\varphi^*\omega)-\calL(Y\otimes{\bfQ}_p,\omega)=
\frac12{\ord}_{p_i}(\deg\varphi_i)=\frac12{\ord}_{p_i}(\deg\varphi).
$$
It remains to make use of the formula which expresses the canonical height in terms of local
terms, and also the equality
$$
\deg\varphi=\prod\deg\varphi_i=\prod p_i^{{\ord}_{p_i}(\deg\varphi_i)}.
$$

\sssect \begin{xxxi}[Corollary] Let $X$ be an abelian variety over $\bfQ$. There exists a finite
set $M$ of prime numbers which depends only on the dimension of $X$ and its set of places of bad
reduction and which satisfies the following condition. The set of abelian varieties $Y$
\r(considered up to isomorphism\r) over $\bfQ$ for which there exists an isogeny $X\rightarrow Y$
of degree not divisible by any of the primes in $M$ is finite.
\end{xxxi}

\begin{xxxr}[Proof] We choose a prime $\ell$ at which $X$ has good reduction, and we let $M$ denote the union of the
set of all primes $\leq C(\ell,\dim X)$, the prime $\ell$, and the set of primes of bad reduction
of $X$. According to Corollary 5.0.2, the canonical height of $X$ is equal to the canonical
height of $Y$. For principally polarized $Y$ the finiteness assertion follows immediately from
the finiteness theorem for heights (\S2.3). In the general case, one must somewhat enlarge the
set $M$, adding to it the primes $p\leq C(\ell,4g)$, and consider the isogeny with principal
polarization of the abelian varieties ${(X\times X')}^4\rightarrow{(Y\times Y')}^4$, whose degree
is also not divisible by any of the primes in $M$ (see the principal polarization theorem and the
remark following it, in \S3). According to Corollary 5.0.2, the abelian varieties ${(X\times
X')}^4$ and ${(Y\times Y')}^4$ have the same height; and, by the finiteness theorem for heights,
the set of abelian varieties \oripage of the form ${(Y\times Y')}^4$ is finite (up to
isomorphism). The required assertion about finiteness of the $Y$ themselves follows from the
well-known theorem of Zassenhaus~\cite{16} about finiteness of the number of classes of exact
left ideals in orders in semisimple $\bfQ$-algebras (in our case the order $\End({(Y\times
Y')}^4)$ in the semisimple $\bfQ$-algebra $\End({(Y\times Y')}^4)\otimes\bfQ$). One could also
use the well-known results of Borel and Harish--Chandra on finiteness of the number of integral
orbits of reductive algebraic groups.
\end{xxxr}

\sssect\begin{xxxi}[Corollary] Let $A$ be an abelian variety over a number field $K$. There
exists a finite set $M$ of prime numbers which satisfies the following condition. The set of
abelian varieties $B$ over $K$ \r(considered up to isomorphism\r) for which there exists an
isogeny $A\rightarrow B$ of degree not divisible by any of the primes in $M$ is finite.
\end{xxxi}

\begin{xxxr}[Proof] Replacing $K$ by a field extension if necessary, we may suppose that it is a Galois
extension of $\bfQ$ (see the finiteness theorem for forms in \S3). One can define the Weil
restriction functor $\Res_{K/\bfQ}$, which takes abelian varieties over $K$ to abelian varieties
over $\bfQ$ and multiplies the dimension of the abelian variety by the degree of $K$. Here for
any abelian variety $C$ over $K$
$$
\Res_{K/\bfQ}C\otimes K\simeq\prod_{\sigma\in\Gal(K/\bfQ)}\sigma C,
$$
where $\sigma C$ is the abelian variety over $K$ whose defining equations are obtained from those
for $C$ by applying the Galois automorphism $\sigma$ to the coefficients. (For more details, see
\cite{23} and~\cite{35}.) The functor $\Res_{K/\bfQ}$ takes the endomorphism multiplication by
$n$ to the endomorphism multiplication by $n$, since it is a functor between categories of
algebraic groups. This implies that it takes isogenies to isogenies, and it takes $p$-isogenies
to $p$-isogenies for any prime $p$.

We set $X=\Res_{K/\bfQ}A$, and for this $X$ we choose the finite set $M$ of primes in accordance
with Corollary 5.0.3, which gives us finiteness of the number of abelian varieties over $\bfQ$ of
the form $\Res_{K/\bfQ}B$, where $B$ is an abelian variety over $K$ for which there exists an
isogeny $A\rightarrow B$ of degree not divisible by any of the primes in $M$. Recalling that the
abelian variety
$$
D=\Res_{K/\bfQ}A\otimes K=\prod_{\sigma\in\Gal(K/\bfQ)}\sigma A
$$
over $K$ contains the abelian variety $A$ as a direct factor, we obtain the required conclusion
about finiteness of the set of all $B$ using Zassenhaus' theorem (see the proof of Corollary
5.0.3) applied to the order $\End(D)$.
\end{xxxr}

\ssect \begin{xxxr}[Proof of theorem 5.0] We set $V=\Ker f$. This is a vector space over the
prime field $\bfZ/p\bfZ$ of some finite dimension $h\leq2g$. We set
$d=\calL(X\otimes{\bfQ}_p,f^*\omega)-\calL(Y\otimes{\bfQ}_p,\omega)$. We must prove that $d=h/2$.
We know that $d$ is an integer between 0 and $g$ (\S4.3.7), and it is connected in the following
way with the action of the inertia group corresponding to $p$ (\S4.3.4). The \oripage Galois
group $\Gamma$ of $\ol{\bfQ}$ acts on $V$. It also acts on the maximal exterior power $\det
V=\bigwedge^hV$ by means of some character
$$
\chi_V\colon \Gamma\rightarrow\Aut(\det V)={(\bfZ/p\bfZ)}^*.
$$
\end{xxxr}
We choose any extension of the $p$-adic valuation ${\ord}_p$ to $\ol{\bfQ}$, and we consider the
corresponding inertia group, which we denote $I(p)$. All of the inertia groups obtained in this
way are subgroups of $\Gamma$ which are conjugate to one another. It follows from local results
(\S4.3.4) that the action of the group $I(p)$ on $\det V$ is given by a power ${\ol\chi}^d_0$ of
the cyclotomic character ${\ol\chi}_0$, i.\,e., the restriction of the character $\chi_V$ to
(any) inertial subgroup $I(p)$ is the same as the restriction of the character
$$
{\ol\chi}^d_0\colon \Gamma\rightarrow{(\bfZ/p\bfZ)}^*,
$$
which is a power of the cyclotomic character ${\ol\chi}_0\colon
\Gamma\rightarrow{(\bfZ/p\bfZ)}^*$ (giving the action of the Galois group of $\ol{\bfQ}$ on the
$p$th roots of 1). Consequently, the characters $\chi_V$ and ${\ol\chi}^d_0$ coincide on all of
the inertia subgroups $I(p)$. If $t$ is a prime $\neq p$, then local results (\S4.5.3) imply that
the restriction of the character $\chi^{c(g)}_V$ to (any) inertia subgroup $I(t)$ is trivial.
Recall that the cyclotomic character ${\ol\chi}_0$ (corresponding to the prime $p$) is also
trivial on $I(t)$. Hence, the character $\chi^{c(g)}_V{\ol\chi}^{-dc(g)}_0$ is everywhere
unramified, and so is trivial (since $\bfQ$ has no nontrivial unramified extensions). Thus,
$$
\chi^{c(g)}_V={\ol\chi}^{dc(g)}_0.
$$

Since $X$ has good reduction at $\ell$, the inertia group $I(\ell)$ acts trivially on the Tate
module $T_p(X)$ (the N\'eron--Ogg--Shafarevich criterion~\cite{24}). In the decomposition group
$D(\ell)$ we choose any element $\sigma$ which acts as the Frobenius automorphism (i.\,e.,
raising to the $\ell$th power) on the algebraic closure of the residue field $\bfZ/\ell\bfZ$.
Then
$$
\chi_0(\sigma)=\ell,\quad{\ol\chi}_0(\sigma)\equiv\ell\mod p,
$$
where $\chi_0\colon \Gamma\rightarrow{\bfZ}^*_p$ is the cyclotomic character which gives the
action of the Galois group of $\ol{\bfQ}$ on the $p^m$th roots of 1 (for all $m$), and
${\ol\chi}_0\colon \Gamma\rightarrow{({\bfZ/p\bfZ})}^*$ is the character which gives the action
on $p$th roots of 1. The image of $\sigma$ in $\Aut T_p(X)$ under the natural homomorphism
$\Gamma\rightarrow\Aut T_p(X)$ is denoted by $F_{\ell}$ and called the \textit{Frobenius element}
\cite{21}, \cite{24}. Its characteristic polynomial
$$
{\calP}_{\ell}(T)=\det(T\id-F_{\ell},T_p(X))
$$
has integer coefficients, and all of the (complex) roots of the polynomial have absolute value
$\sqrt{\ell}$ (Weil's theorem). This implies that the characteristic polynomial $\calP(T)$ of the
automorphism a acting on the ${\bfZ}_{\ell}$-module ${[\bigwedge^hT_p(X)]}^{\otimes c(g)}$, has
integer coefficients, and all of its (complex) roots have absolute value equal to
${\ell}^{hc(g)/2}$. On the other hand, since the action of the Galois group $\Gamma$ on $\det
V^{\otimes c(g)}$ is given by the character
$$
\chi_V^{c(d)}={\ol\chi}_0^{dc(g)}\colon \Gamma\rightarrow\Aut(\det V^{\otimes c(g)})={({\bfZ/p\bfZ})}^*,
$$
\oripage the element $\sigma$ acts on $\det V^{\otimes c(g)}$ by multiplication by
${\ell}^{dc(g)}$. Since
$$
\det V^{\otimes c(g)}\subset{[\smash{\bigwedge^h}X_p]}^{\otimes c(g)}={[\smash{\bigwedge^h}T_p(X)]}^{\otimes c(g)}/(p)
$$
is a Galois submodule of
$$
{[\smash{\bigwedge^h}T_p(X)]}^{\otimes c(g)}/(p),
$$
the integer $a=\calP({\ell}^{dc(g)})$ must be divisible by $p$. If $d\neq h/2$, then the number
${\ell}^{dc(g)}$ could not be a root of $\calP$, i.\,e., $a\neq0$. Hence, if $d\neq h/2$, then
$p$ cannot exceed
\begin{align*}
|a|&=|\calP({\ell}^{dc(g)})|\leq({\ell}^{dc(g)}+{\ell}^{hc(g)/2})\deg\calP\\
&\leq{(2{\ell}^{dc(g)})}^{c(g)\bigl(\begin{smallmatrix}2g\\h\end{smallmatrix}\bigr)}\\
&\leq{(2{\ell}^{dc(g)})}^{c(g)\bigl(\begin{smallmatrix}2g\\g\end{smallmatrix}\bigr)}=C(\ell,g).
\end{align*}

\ssect Suppose that $X$ is an abelian variety over $\bfQ$, $W$ is a $p$-divisible subgroup of
$X(p)$ defined over $\bfQ$, $T(W)$ is its Tate module, and $V(W)=T(W)\otimes{\bfQ}_p$ is the Tate
${\bfQ}_{\ell}$-module. The Galois group $\Gamma$ of $\ol\bfQ$ acts on $T_p(X)$, and $T(W)$ is an
invariant $\Gamma$-submodule. The rank of the ${\bfZ}_p$-module $T(W)$ is equal to the height $h$
of the $p$-divisible group $W$. We consider the action of $\Gamma$ on the maximal nonzero
exterior power of the Tate module
$$
\det W=\bigwedge^h_{\bfZ_p}T(W),
$$
which is given by a character
$$
\chi_W\colon \Gamma\rightarrow\Aut(\det W)={\bfZ}^*_p.
$$

From local results (\S4.4) it follows that the restriction of $\chi_W$ to (any) inertia subgroup
$I(p)$ is equal to $\epsilon\chi^d_0$, where $\chi_0\colon \Gamma\rightarrow{\bfZ}^*_p$ is the
cyclotomic character corresponding to the prime $p$, and $\epsilon\colon I\rightarrow{\bfZ}^*_p$
is a character of finite order. The exponent $d$ is a nonnegative integer.

\sssect \begin{xxxr}[Theorem] $d=h/2$.\end{xxxr} The proof will be given in \S5.5.

\sssect\begin{xxxi}[Corollary] Under the assumptions of \S5.2, the sequence $Y_{(1)}=X/W_1$,
$Y_{(2)}=X/W_2,\dots$ contains only a finite number of pairwise nonisomorphic \r(over $\bfQ$\r)
abelian varieties, i.\,e., there exists a finite set $Y_{(i_1)},\dots,Y_{(i_r)}$ of abelian
varieties such that any $Y_{(n)}$ is isomorphic to one of the $Y_{(i_1)},\dots,Y_{(i_r)}$.
\end{xxxi}

\ssect \begin{xxxr}[Proof of corollary 5.2.2 (assuming Theorem
5.2.1)]\end{xxxr}

\begin{xxxi}[Lemma] Let

$$
f_n\colon Y_{(n)}\rightarrow Y_{(n+1)}=Y_{(n)}/(W_{n+1}/W_n)
$$

\oripage be the natural isogenies of degree $\#W_1=p^h$. There exists a natural number $N$ such
that for all $n\geq N$ and for $\omega$ a nonzero form of highest weight on $Y_{(n+1)}$
$$
\calL_s(Y_{(n)}\otimes{\bfQ}_p,f^*_n\omega)-\calL_s(Y_{(n+1)}\otimes{\bfQ}_p,\omega)=d=h/2={\ord}_p(\deg f_n)/2.
$$

In particular, all of the abelian varieties $Y_{(N)}$, $Y_{(N+1)},\dots$ have the same canonical
height.
\end{xxxi}

\begin{xxxr}[Proof of the lemma] If $X$ has semistable reduction at $p$, then our assertion follows immediately
from Theorem 5.2.1 by applying local results (\S4.4).
\end{xxxr}

Now suppose that $X$ does not necessarily have semistable reduction at $p$. Let $K$ denote the
finite Galois extension of $\bfQ$ which is the field of definition of all of the points of order
fifteen on $X$. We fix a place $\gotp$ of $K$ lying over $p$, and we let $\ord_{\gotp}\colon
K^*\rightarrow\bfZ$ be the standard discrete valuation corresponding to $\gotp$. Let $e$ denote
the ramification index of the place $\gotp$ over $p$; by definition, $e={\ord}_{\gotp}(p)$. For
any natural number $m$, the local (stable) term for $(Y_{(m)}\otimes{\bfQ}_p,\omega)$
corresponding to $p$ is equal to (see Remark 2 of \S2.1)
$$
\calL(Y_{(m)}\otimes K_{\gotp},\omega)/e,
$$
since $X\otimes K_{\gotp}$, and hence also any of the $Y_{(m)}$, have semistable reduction at
$\gotp$. It is not hard to see that the restriction of the character $\chi_W$ to the inertia
group $I(\gotp)$ is also equal to the product of a character of finite order and $\chi^d_0$ (here
$I(\gotp)$, the inertia group of the place $\gotp$, is a subgroup of finite index in $I(p)$). As
before, local results ensure us that there exists a natural number $N$ such that for all $n\geq
N$
$$
\calL(Y_{(n)}\otimes K_{\gotp},f^*_n\omega)-\calL(Y_{(n+1)}\otimes K_{\gotp},\omega)=d{\ord}_{\gotp}(p)=de=he/2.
$$

But this gives us the required equality for the local (stable terms) correspond ing to
$Y_{(n)}\otimes{\bfQ}_p$ and $Y_{(n+1)}\otimes{\bfQ}_p$ (see Remark 2 in \S2.1).

\textit{End of the proof of Corollary} 5.2.2. If all of the $Y_{(n)}$ are principally polarized,
then our assertion follows from Lemma 5.3.1 and the finiteness theorem for the height (see Remark
2 in \S2.1). In the general case, there exists (\S3) a $p$-divisible subgroup $H$ of ${(X\times
X')}^4$ such that
$$
{(X\times X')}^4/H_n\simeq{(Y_{(n)}\times {Y'}_{(n)})}^4.
$$

Since all of the ${(Y_{(n)}\times {Y'}_{(n)})}^4$ are principally polarized, we find (applying
the lemma in \S5.3 to ${(X\times X')}^4$) that the sequence ${\{{(Y_{(n)}\times
{Y'}_{(n)})}^4\}}^\infty_{n=1}$ contains only a finite number of nonisomorphic abelian varieties.
But this gives us the required finiteness assertion for the $Y_{(n)}$ (see the end of the proof
of Corollary 5.0.3).

\ssect \begin{xxxi}[Corollary] Let $A$ be an abelian variety over an arbitrary number field $K$,
and let $H$ be an arbitrary $p$-divisible subgroup $A(p)$ which is defined over $K$. Then the
sequence $B_{(1)}=A/H_1$, $B_{(2)}=A/H_2,\dots$ contains only a finite number of pairwise
nonisomorphic \r(over $K$\r) abelian varieties.
\end{xxxi}
\oripage
\begin{xxxr}[Proof] We may suppose that $K$ is a Galois extension of $\bfQ$. Let $\varphi_n\colon
X\rightarrow Y_{(n)}$ be the canonical isogenies. We set $X={\Res}_{K/\bfQ}A$,
$Y_{(n)}={\Res}_{K/\bfQ}B_{(n)}$. Then the kernels $V_n$ of the isogenies
$f_n={\Res}_{K/\bfQ}\varphi_n\colon A\rightarrow B_{(n)}$ form a $p$-divisible group over $\bfQ$
whose height is equal to the height of $H$ multiplied by the degree of the field $K$. In order to
see this, it is sufficient to note that under the natural isomorphisms (see the proof of
Corollary 5.0.4)
$$
X\otimes K\simeq\prod\sigma A,\quad Y_{(n)}\otimes K\simeq\prod\sigma B_{(n)}\quad(\sigma\in\Gal(K/\bfQ))
$$
the isogeny $f_n\otimes K$ goes to $\prod\sigma\varphi_n$. Now Corollary 5.2.2 gives us
finiteness of the number of pairwise nonisomorphic abelian varieties $Y_n$ over $\bfQ$, which
implies finiteness for the abelian varieties $\prod_{\sigma\in\Gal(K/\bfQ)}\sigma B_{(n)}$ over
$K$. From this one deduces the required finiteness assertion for the varieties $B_{(n)}$ (see the
end of the proof of Corollary 5.0.4).
\end{xxxr}
\ssect \begin{xxxr}[Proof of theorem 5.2.1] This proof is in many ways similar to that of Theorem
5.0. We fix a natural number $r$ such that the character $\epsilon'$ is trivial, and we consider
the action of the Galois group $\Gamma$ on the free rank 1 ${\bfZ}_{\ell}$-module $\det
W^{\otimes rc(g)}$, which is a Galois submodule of $\bigwedge^hT_p{(X)}^{\otimes rc(g)}$ (here
$g$ is the dimension of $X$). The action of $\Gamma$ on $\det W^{\otimes rc(g)}$ is given by the
character $\chi^{rc(g)}_W\colon \Gamma\rightarrow\Aut(\det W^{\otimes rc(g)})={\bfZ}^*_p$. The
restriction of $\chi^{rc(g)}_W$ to the inertia group $I(p)$ is equal to
$$
{(\epsilon\chi^d_0)}^{rc(g)}=\chi^{drc(g)}.
$$
For any prime $t\neq p$, the action of the inertia group $I(t)$ on $\det W^{\otimes rc(g)}$ turns
out to be trivial, by \S4.5.2. By the same token, the character
$$
\chi^{rc(g)}_W\chi^{-drc(g)}_0\colon \Gamma\rightarrow{\bfZ}^*_p
$$
of the Galois group of $\ol\bfQ$ is unramified, and hence is trivial, i.\,e.,
$$
\chi^{rc(r)}_W=\chi^{drc(g)}_0.
$$
We now choose a prime $\ell\neq p$ at which $X$ has good reduction, and we fix an element
$\sigma$ in the decomposition group $D(\ell)$ which acts as the Frobenius automorphism on the
algebraic closure of the residue field $\bfZ/\ell\bfZ$. We recall (\S5.1) that
$\chi_0(\sigma)=\ell$, and that the characteristic polynomial of the action of $\sigma$ on
$T_p(X)$ has integer coefficients and has (\textit{complex}) roots all of absolute value
$\sqrt{\ell}$. This means that all of the eigenvalues of $\sigma$ acting on
$\bigwedge^hT_p{(X)}^{\otimes rc(g)}$ are algebraic numbers of absolute value
${\ell}^{hrc(g)/2}$. On the other hand, $\sigma$ acts on $\det W^{\otimes
rc(g)}\subset\bigwedge^pT_p{(X)}^{\otimes rc(g)}$ by multiplication by the number
$$
\chi^{drc(g)}_0(\sigma)=\chi_0{(\sigma)}^{drc(g)}={\ell}^{drc(g)}.
$$
Consequently, $d=h/2$.
\end{xxxr}

\ssect As we noted in the Introduction and in \S3, Corollary 5.4 implies Tate's conjecture on
homomorphisms of Tate modules and the semisimplicity of the Tate module. The purpose of this
subsection is to propose a simpler approach to proving Tate's conjecture in certain important
special cases.

\oripage
\begin{xxxi}[Theorem] Let $X$ be an abelian variety over $\bfQ$ which has good reduction at the odd prime $p$. Let
$f\colon X\rightarrow Y$ be a $p$-isogeny for which there exist an abelian variety $Z$ over
$\bfQ$ having good reduction at $p$ and a power $q=p^n$ of $p$ such that the Galois modules $\Ker
f$ and $Z_q$ are isomorphic. Then the abelian varieties $X$ and $Y$ have the same canonical
height. The set of abelian varieties $Y$ over $\bfQ$ \r(considered up to isomorphism\r) for which
there exists an isogeny $X\rightarrow Y$ with the above properties is finite.
\end{xxxi}

\begin{xxxr}[Proof] The claim that the heights are equal follows from local results (\S4.3.5). The finiteness
theorem for heights tells us that there are only finitely many $Y$ with principal polarization.
In order to reduce everything to the principally polarized case, one observes that the kernel of
the dual isogeny $f'\colon Y'\rightarrow X'$ is isomorphic (as a Galois module) to $Z'_q$, and
hence the abelian varieties $X'$ and $Y'$ have the same canonical height. This means that the
principally polarized abelian varieties ${(X\times X')}^4$ and ${(Y\times Y')}^4$ have the same
canonical height, which implies that there are only finitely many abelian varieties of the form
${(Y\times Y')}^4$. But then, as above, we conclude that the set of abelian varieties $Y$ is
finite.

We note that, to prove Tate's homomorphism conjecture for an abelian variety $X$ and a prime $p$
(see~\cite{27}, \cite{11}, and~\cite{12}), we must establish finiteness of the set of abelian
varieties $Y$ for which there exists an isogeny $X\times X\rightarrow Y$ whose kernel is the
graph of some endomorphism of the Galois module $X_{p^n}$, which, being a graph, is isomorphic to
$X_{p^n}$. By the same token, Theorem 5.6.1 implies Tate's conjecture on homomorphisms of Tate
$p$-modules for abelian varieties over $\bfQ$ having good reduction at the odd prime $p$.
\end{xxxr}

\section{The finiteness theorem for characters of Galois representations. Proof of the Mordell
conjecture}

In order to complete the proof of the Shafarevich conjecture according to the outline in \S1, it
now remains for us to obtain the following result.

\begin{xxxi}[Finiteness theorem for characters of Galois representations] Let $K$ be an algebraic number field,
and let $G=\Gal(\ol K/K)$ be the Galois group of its algebraic closure $\ol K$. Fix a finite set
$S$ of finite places of $K$, a prime $\ell$ and an integer $n\geq1$. Then there exists a finite
set $Q$ of finite places of $K$ which is disjoint from $S$ and has the property that the
character $\chi_\rho(g)=\Tr(\rho(g))$ of any $n$-dimensional continuous $\ell$-adic
representation $\rho\colon G\rightarrow\GL(n,{\bfQ}_{\ell})$ unramified outside $S$ is uniquely
determined by its values $\chi_\rho(\Fr_v)$ on the Frobenius automorphisms $\Fr_v$, $v\in Q$.
\end{xxxi}

We recall that if $L/K$ is a Galois extension (perhaps infinite) which is unramified at the place
$v$ of $K$, then the Frobenius element ${\Fr}_v$ in $\Gal(L/K)$ is defined uniquely up to
conjugation~\cite{2}, \cite{21}. In the situation of the theorem, the representation $\rho$
factors a representation of the group $G_S=\Gal(K_S/K)$, where $K_S$ is the maximal extension of
$K$ which is unramified \oripage outside $S$. Thus, for any place $v\not\in S$, the elements
${\Fr}_v$ of $G$ are (nonuniquely) defined (first in $G_S$, and then they are lifted to $G$).

We now construct the set $Q$. Let $L/K$ be the compositum of all Galois extensions of $K$ which
are unramified outside $S$ and have degree $\leq{\ell}^{2n^2}$. According to Hermite's theorem
\cite{L5}, this is a finite extension. Applying the Chebotarev density theorem (see~\cite{21}) to
the extension $L/K$, we see that there exists a finite set $Q$ of places of $K$ which is disjoint
from $S$ and has the property that the automorphisms $\Fr_v$, $v\in Q$, exhaust the Galois group
$\Gal(L/K)$.

We now verify that the set $Q$ satisfies the conclusion of the theorem. Let $\rho_1$ and $\rho_2$
be two representations whose characters coincide on the elements $\Fr_v$, $v\in Q$. We consider
the representation
$$
\rho_1\times\rho_2\colon {\bfZ}_{\ell}[G]\rightarrow\Mat(n,{\bfQ}_{\ell})\times\Mat(n,{\bfQ}_{\ell})
$$
of the group algebra ${\bfZ}_{\ell}[G]$. Its image $M$ is a ${\bfZ}_{\ell}$-subalgebra of rank at
most $2n^2$. The image of the Galois group $G$ in the group of invertible elements of the algebra
$M/\ell M$, which is a vector space over ${\bfF}_{\ell}$, determines a \textit{finite} Galois
extension which is contained in $L$. By the construction of $Q$, the elements
$\rho_1\times\rho_2(\Fr_v)$, $v\in T$, span $M/\ell M$ over ${\bfF}_{\ell}$, and hence they span
the module $M$ over ${\bfZ}_{\ell}$ (by Nakayama's lemma). We now consider the linear form
$$
f(a_1,a_2)=\Tr(a_1)-\Tr(a_2),\quad a_{1,2}\in\Mat(n,{\bfQ}_{\ell}),
$$
on $M$. By assumption, $\chi_{\rho_1}(\Fr_v)=\chi_{\rho_2}(\Fr_v)$, $v\in T$, and hence $f=0$ on
the elements $\rho_1\times\rho_2(\Fr_v)$, $v\in Q$, and so on all of $M$; this gives
$\chi_{\rho_1}=\chi_{\rho_2}$. The theorem is proved.

We shall show that the finiteness theorem for abelian varieties (the fundamental theorem of \S1)
implies the original Shafarevich conjecture for algebraic curves.

Recall that if $\Lambda_K$ is the ring of integers of $K$, then any curve $X$ of genus $g>1$ can
be associated to a two-dimensional regular scheme $f\colon \bfY\rightarrow\Spec\Lambda_K$ with
generic fiber $X$ (the \textit{minimal model}; see~\cite{Sh3} or~\cite{38}, Expos\'e X).

By analogy with the definition of the N\'eron minimal model in \S2, the scheme $\bfY$ is
characterized by the condition that, given any regular projective scheme $g\colon
\bfZ\rightarrow\Spec\Lambda_K$, the isomorphism $\bfZ\otimes K\stackrel\sim{\rightarrow} X$
extends uniquely to an isomorphism between $\bfZ$ and $\bfY$. This is equivalent to requiring
that the fibers of $f$ not have any exceptional curves of the first kind (\cite{32}, Chapter 1)
(see below). For all $v$ except for finitely many $v\in S$---the \textit{places of bad
reduction}---the fibers of $f$ are nonsingular projective curves of genus $g$.

\begin{xxxi}[Theorem] Suppose that $K$ is a number field, $S$ is a finite set of places of $K$, and $g\geq1$ is an
integer. There exist only finitely many genus $g$ algebraic curves $X$ defined over $K$
\r(considered up to isomorphism over $K$\r) which have good reduction everywhere outside $S$.
\end{xxxi}

\begin{xxxr}[Proof] To every curve $X$ as in the theorem we associate its Jacobian $A=J(X)$ with principal
polarization $\theta$. The variety $A$ and the polarization \oripage $\theta$ are defined over
$K$, and $A$ has good reduction outside of $S$. According to Remark 2 in \S2, the degree 1
polarizations over $K$ on a given abelian variety are finite in number up to the action of
$K$-automorphisms (this fact holds for any field $K$; see~\cite{54}). Thus, from the fundamental
theorem of \S1 it follows that the number of $K$-isomorphism classes of pairs $(A,\theta)$ is
finite. By the Torelli theorem~\cite{72}, the curve $X$ is determined by the pair $(A,\theta)$ up
to isomorphism over $\ol K$. In order to obtain finiteness over $K$, for each abelian variety $A$
we consider the extension $L_A/K$ obtained by adjoining all of the coordinates of the points of
order $n\geq3$ ($n$ is fixed in the sequel). The extensions $L_A/K$ are unramified outside
$S\cup\{\text{divisors of}\ n\}$ and are of bounded degree. By Hermite's theorem, they are all
contained in a finite extension $K'$ of $K$. We now show that the curves $X$ in the theorem form
a finite set up to $K'$-isomorphism. To do this it suffices to note that, if $\sigma\in\Gal(\ol
K/K')$, and if $\varphi$ is an isomorphism of two curves with the same Jacobian $A$, then
$\varphi^\sigma\circ\varphi^{-1}$ induces the identity map on $A_n$, and hence it is the identity
(Serre's lemma; see~\cite{Mum2}, p.\,207). In order to go from $K'$ to $K$, one uses the fact
that the group $H^1(\Gal(K'/K),{\Aut}_{K'}(X))$ is finite in our situation. The proof is
complete.
\end{xxxr}

In order to obtain the Mordell conjecture from this theorem, we shall need the following
\textit{construction of ramified coverings}. Let $X$ be an algebraic curve of genus $g\geq2$ over
an \textit{arbitrary} field $K$, and let $P\in X(K)$. We fix an integer $n>1$, and we set
\begin{equation}\label{eq5.1}
\begin{CD}
\kern-5em X'=X\times_{J(X)}J(X) @>>> J(X)\\
@V\pi VV @VVnV\\
X @>\varphi >> J(X)
\end{CD}
\end{equation}
where $\varphi$ is the map taking $Q$ to $Q-P$. Let $L/K$ be a field extension of $K$ over which
the cycle $\pi^{-1}(P)$ is a sum of rational points. We choose $P_0\in\pi^{-1}(P)$ and we
consider the generalized Jacobian $J_{\gotm}(X')$ of the curve $X'$ corresponding to the cycle
$\gotm=\pi^{-1}(P)-P_0$ \cite{Ser1}.
We set
\begin{equation}\label{eq5.2}
\begin{CD}
\kern-5em X_P=X'\times_{J_{\gotm}(X')}J_{\gotm}(X') @>>> J_{\gotm}(X')\\
@VVV @VVnV\\
X' @>\psi >> J_{\gotm}(X')
\end{CD}
\end{equation}
where $\psi$ is the rational map taking $Q\not\in\gotm$ to $Q-P_0$. Then $X_P$ is a nonsingular
projective curve, and it is uniquely determined by this diagram of rational maps. It is defined
over the field $L_P$.

\begin{xxxi}[Properties of the construction of $X_P$ and $L_p$] Let $K$ be an algebraic number field, and let $S$ be the
set of places of bad reduction of the curve $X$. Then the following assertions are true.
\begin{enumerate}

\item The extension $L/K$ is unramified outside $S\cup\{\text{divisors of}\ n\}$, and its degree depends only on
$g$ and $n$.

\item The genus of the curve $X_P$ depends only on $g$ and $n$, and the covering $X_P\rightarrow X$ is ramified
at the point $P$ and only at that point.

\oripage

\item The curve $X_P$ has good reduction outside of the places of $L$ lying over $S$ and the divisors
of $n$.
\end{enumerate}
\end{xxxi}

Property 1 follows from the general properties of extensions obtained by adjoining points of
finite order on an abelian variety~\cite{24}. Property 2 is purely geometric in nature, and is
verified immediately.

As for property 3, we shall only give the idea of the proof, referring the reader interested in
the details (and different versions) of the proof to~\cite{Par1} \cite{55}, \cite{69}, pp.
71--76, \cite{41}, pp.\,191--197, and~\cite{71}. We remove from $\bfR=\Spec\Lambda_K$ all places
lying over $S$ or dividing $n$, and we remove the corresponding fibers from $\bfY$. We keep the
same notation $f\colon \bfY\rightarrow\bfR$ for the resulting scheme. The point $P$ determines a
section $s\colon \bfR\rightarrow\bfY$ of the morphism $f$. We have the commutative diagram
\begin{equation}\label{eq5.3}
\begin{array}{clc}
\bfY&\stackrel\varphi{\longrightarrow}&\bfA\\
\!\!{\scriptstyle f}\downarrow\ &\ \swarrow\!\!{\scriptstyle g}&\\
\bfR&&
\end{array}
\end{equation}
in which $g\colon \bfA\rightarrow\bfR$ is the N\'eron model of the abelian variety $A$ (see \S2),
and $\varphi$ is an extension of the morphism $\varphi$ in diagram \eqref{eq5.1}. Multiplication
by $n$ in the group scheme $\bfA$ leads to the surface $\bfY'\rightarrow\bfR$, where
$\bfY'=\bfY\times_{\bfA}\bfA$ and the generic fiber of $\bfY'\otimes K$ coincides with the curve
$X'$ in \eqref{eq5.1}. At every place $v\in\bfR$ the fiber ${\bfY'}_v$ will contain a nonsingular
projective curve whose genus is equal to the genus of the generic fiber $X'$. If we pass from
$\bfY'$ to a minimal model $\bfY''$ of $X'$ over $\bfR$, then this property of the fibers is
preserved (neither resolving the singularities of the scheme $\bfY'$ nor contracting the
exceptional curves of the first kind in the fibers can destroy such a curve). Thus, property 3
for the curve $X'$ follows from the following result.

\begin{xxxi}[Lemma] Let $f\colon \bfY\rightarrow\bfS$, $\bfS=\Spec\goto$, be a regular scheme having smooth irreducible generic fiber
which is a projective curve $X$ of genus $g>1$, and let $\goto$ be a complete discrete valuation
ring. Suppose that the closed fiber ${\bfY}_0$
\begin{enumerate}[\upshape1)]
\item does not contain any exceptional curves of the first kind, and
\item contains an irreducible curve ${\bfY}_1$ whose normalization has genus equal to $g$.
\end{enumerate}
Then the morphism $f$ is smooth, and hence the fiber ${\bfY}_0$ is a nonsingular projective curve
of genus $g$.
\end{xxxi}
\begin{xxxr}[Proof] We shall assume that the residue field of $\goto$ is algebraically closed. It is well known
(see~\cite{Sh3} or~\cite{38}, Expos\'e X) that there is an intersection theory for divisors on
the scheme $\bfY$. If $N=\mathop{\mathrm{NS}}(\bfY)$ is the N\'eron--Severi group, then the space
$N\otimes\bfQ$ is spanned by the components $Y_1,\dots,Y_n$ of the fiber, and the intersection
index $C\cdot D$ determines a quadratic form of signature $(n-1,0)$ on this space. In particular,
this means that if $n>1$, then $C\cdot C<0$ for any irreducible curve $C$ on $\bfY$. In addition
to the intersection theory on $\bfY$, one has the relative canonical class $\omega_{\bfY}$
(\cite{38}, Expos\'e X) and the \oripage adjunction formula
\begin{equation}\label{eq5.4}
C\cdot\omega_{\bfY}+C\cdot C=2\pi-2,
\end{equation}
in which $C$ is any complete curve on $\bfY$ and $\pi$ is its arithmetic genus. If $C$ is
irreducible, then $\pi=g'+\delta$, where $g'$ is the genus of its normalization and $\delta $is
the number of singular points counting multiplicity. The invariance of the intersection index and
the arithmetic genus under specialization implies that
\begin{equation}\label{eq5.5}
{\bfY}_0\cdot{\bfY}_0=0\mtxt{and}{\bfY}_0\cdot\omega_{\bfY}=2g-2.
\end{equation}
\end{xxxr}
Suppose that $n>1$ and $C$ is an irreducible curve. Then $C\cdot C<0$, and from \eqref{eq5.4} we
find that $C\cdot\omega_{\bfY}\geq-1$ and the equality $C\cdot\omega_{\bfY}=-1$ is equivalent to
$C$ being an exceptional curve of the first kind ($g'=\delta=0$, $C\cdot C=-1$). But by 1), this
is impossible, and so always $C\cdot\omega_{\bfY}\geq0$. Now let ${\bfY}_0=\sum_{i\geq1}a_iY_i$
be the expansion of the cycle ${\bfY}_0$ in terms of irreducible components. In this case the
equality ${\bfY}_0\cdot\omega_{\bfY}=\sum a_iY_i\cdot\omega_{\bfY}$ takes the following form, by
\eqref{eq5.5}:
$$
2g-2=a_1(2g-2+\delta)+\sum a_iY_i\cdot\omega_{\bfY},
$$
where all of the terms on the right are nonnegative. This implies that $a_1=1$, $\delta=0$,
$Y_1\cdot\omega_{\bfY}=2g-2$, and hence $Y_1\cdot Y_1=0$, which contradicts the supposition that
$n>1$. Thus, $n=1$, and $a_1=1$, $\delta=0$, which gives us the required conclusion.

In order to complete the verification of property 3, we must now construct a smooth group scheme
over $\bfR'$ ($=\Spec\Lambda_L$ with the places lying over $S$ and the ones dividing $n$
removed)---a model of the generalized Jacobian variety $J_{\gotm}(X')$ in \eqref{eq5.2}---and
determine the analog of diagram \eqref{eq5.3} in this situation. As before, the required property
is obtained from the lemma just proved.

We proceed to the proof of the central result of this survey.

\begin{xxxi}[The Mordell conjecture] Let $K$ be an algebraic number field, and let $X$ be a nonsingular projective
algebraic curve over $K$ of genus $g>1$. Then the set of rational points $X(L)$ is finite for any
finite extension $L/K$ of $K$.
\end{xxxi}

In this form, this result was not known for a single curve before Faltings.

The proof proceeds in exactly the same way for any field $L$. Let $P\in X(K)$.

We consider the set of algebraic curves of the form $X_P$ in the construction in \eqref{eq5.2}.
By the properties of that construction, the fields of definition $L_P$ of the curves $X_P$ form a
finite set (Hermite's theorem). Furthermore, the curves $X_P$ with a fixed field of definition
satisfy all of the hypotheses of the Shafarevich conjecture, and so they form a finite set. We
consider a covering $X_P\rightarrow X$. Since $P$ is the only ramification point, it is uniquely
determined by the covering, and one must use the classical fact that, given curves $X_P$ and $X$
of genus greater than 1, there exist only finitely many coverings $X_P\rightarrow X$ (de
Franchis' lemma, which itself is a special case of the Mordell conjecture for algebraic function
fields, see the commentary section of Chapter 8 in Lang's book).

\oripage

We conclude by noting some unsolved problems which arise in connection with the
proof of the Mordell conjecture. Problems 1--5 are due to Parshin; 6 and 7 to
Zarkhin.

\medskip{\bf1.} \textbf{The problem of effectiveness.} How can one construct an algorithm to find the
rational points on a curve $X$? Or, more precisely, how can one write an explicit bound for the
heights of the rational points in terms of the invariants of the field $K$, the set $S$, and the
genus $g$? It is natural to take for the height the canonical height on $X(K)$ relative to the
invertible sheaf $\Omega^1_X$ (see \S2). One possible approach to this problem is described in
\S1. We note that a method which in principle enables one to bound the \textit{number} of
rational points was proposed by Parshin (see~\cite{71}, the Introduction and the reports by
Raynaud (Expos\'e VII) and by Szpiro (Expos\'e XI), which contain a detailed elaboration of the
method). However, this route can lead only to very rough estimates (thus, in the case of the
Fermat curve $x^{\ell}+y^{\ell}=1$ with $\ell$ prime, the constant cannot be better than
$2^{2^{\ell}}$ ).

\medskip{\bf2.}  As we have seen, the proof of the Mordell conjecture does not at all make use of
the Mordell--Weil theorem on the group of rational points on an abelian
variety---which had always been an ingredient in the earlier attempts to prove
the conjecture. Nevertheless, it is worthwhile to compare the two results.
Barry Mazur has conjectured that the Mordell--Weil theorem remains true if one
replaces the ground field $K$ by its (infinite) cyclotomic $\ell$-extension
$K(\ell)$
\cite{4**,Man2}.

\begin{xxxi}[Conjecture] Let $X$ be an algebraic curve of genus $g>1$ defined over an algebraic number field
$K$. Then the set $X(K(\ell))$ is finite.
\end{xxxi}
In support of this conjecture one can cite the corresponding result in the case
of an algebraic function field with finite field of constants. In addition,
part of the proof of the Mordell conjecture over a number field carries over to
this situation. In particular, using constructions of Zarkhin in~\cite{15}, one
can show that Tate's homomorphism conjecture (and his conjecture on
semisimplicity of the Tate module) remain valid over $K(\ell)$.


\medskip\textbf{Note added in December 2009}. In order to prove the Tate's homomorphism conjecture for
 over $K(\ell)$, one may use Bogomolov's theorem on $\ell$-adic Lie algebras \cite{Bo1,Bo2} mentioned in Subsection
  7 below.

\medskip{\bf3.} For which varieties $X$ of dimension greater than 1 is the Mordell conjecture in the
above form valid? As Raynaud showed (\cite{61} combined with Faltings'
theorem), the assertion holds for closed subvarieties of abelian varieties
which do not contain a translation of an abelian variety. In the case of
algebraic function fields with complex coefficients, another proof of this fact
was found by Parshin (his method uses the properties of the hyperbolic
Kobayashi metric).
\footletter{\textit{Added in translation}. Using this method, one can also
prove Lang's conjecture on integral points of affine subsets of abelian
varieties (in the case when the ground held is a field of algebraic functions
over $C$ and the
hyperplane section does not contain a shift of the abelian variety).
See~\cite{3*} and Chapter 8 of Lang's book.} Other suitable candidates for the
conjecture are varieties with ample cotangent bundle (in the geometrical
situation Grauert's \oripage method~\cite{Gra} is easily applicable to them,
see~\cite{No}) and varieties which are Kobayashi hyperbolic as complex-analytic
manifolds (here
even the function field analog of the conjecture is unknown; see~\cite{L14}).


\medskip\textbf{Note added in December 2009}. M. Raynaud \cite{61} has proved that
the finiteness theorem for closed subvarieties of abelian varieties can be
reduced to the number field case. However, this does not give  the proof in the
number field case. The finiteness theorem  was later proven by  Faltings
\cite{1**,2**}
 in the number field case together with
Lang's conjecture mentioned in footnote (E).


\medskip{\bf4.} \textbf{Construction of ramified coverings.} This construction also leads to several interesting
problems. If we fix an integer $n$, then the set of abelian varieties $A_P=J(X_P)$, $P\in
X(\bfC)$, gives a smooth family of abelian varieties over the curve $X$ which has no
degenerate fibers. By the same token, one has a map of $X$ (or a finite covering of it) to
the moduli space $\gotM$ of curves or abelian varieties (the first examples of maps of
complete curves to a moduli space that is necessarily not complete were constructed by
Kodaira~\cite{49}). Which conditions characterize the images of curves in $\gotM$, for
example, relative to the natural metrics on $\gotM$---the Teichm\"uller or Weil--Petersson
metrics? What is the structure of the monodromy of the families $A_P/X$? From an arithmetic
point of view, it would be interesting to understand the nature of the correspondences which
this construction gives canonically on any curve. To a point $P\in X$ one associates all
points $Q\in X$ for which the variety $J(X_Q)$ is isogenous to $J(X_P)$. If $X$ is a modular
curve or a curve that is uniformized by quaternionic groups, then is there a connection
between these correspondences and the usual Hecke correspondences for such curves?

\medskip{\bf5.} As noted in \S1, Shafarevich's conjecture implies that the number of $S$-integral points in
the moduli space of algebraic curves of genus $g\geq1$ is finite. In the case of curves of
genus 1, we find that there are only finitely many $S$-integral points on the moduli curves
of level $n\geq3$. On the other hand, according to the Mordell conjecture, the set of all
rational points on a modular curve of level $n>6$ (in which case the genus is greater than 1)
is finite. In accordance with the conjectures in 3 above, one would expect that the same is
true for the moduli space ${\gotM}^{ab}_{g,n}$  of abelian varieties of sufficiently large
level $n$ (depending only on the dimension $g$). Thus, it is natural to investigate whether
the corresponding finiteness assertion holds for the $\ell$-adic representations. The exact
statement is as follows. Let $\rho\colon \Gal(\ol K/K)\rightarrow\GL(V)$ be a continuous
$\ell$-adic representation, where $V\cong D{\bfZ}^m_{\ell}$, and suppose that the image of
$\rho$ in $V/{\ell}^aV$ is trivial for some $a$. Does there exist a function $a=a(m)$ such
that the number of semisimple representations satisfying the Riemann hypothesis (see
\cite{21}, Chapter 1) and this condition and having the given dimension $m$, is finite? As
before, this question reduces to finding a suitable finiteness theorem for characters
(i.\,e., constructing a finite set $Q$). The analogous question for complex representations
of discrete groups (see \S1) has also not been studied.

\medskip{\bf6.} \textbf{Variants of the Tate homomorphism conjecture.} Suppose that $K$ is a number field, $\ol K$ is
its algebraic closure, and $G=\Gal(\ol K/K)$ is the Galois group. Faltings (\cite{39}, p.\,205,
Theorem 5) has proved that for any abelian \oripage varieties $X$ and $Y$ over $K$, the natural
map
$$
\Hom(X,Y)\rightarrow{\Hom}_G(X(\ol K), Y(\ol K))
$$
is a bijection.

Consider the infinite-dimensional vector spaces
$$
X{(\ol K)}_{\bfQ}=X(\ol K)\otimes{\bfQ}, \quad Y{(\ol K)}_{\bfQ}=Y(\ol K)\otimes{\bfQ}
$$
over the rational number field $\bfQ$, on which the group $G$ acts naturally. Choosing a Tate
height makes these actions unitary. For any finite algebraic extension $L$ of $K$ lying in
$\ol K$, the invariants of the Galois group $\pi=\Gal(\ol K/L)$ of $L$ are $X(L)\otimes\bfQ$
and $Y(L)\otimes\bfQ$, respectively. In particular, if the Galois modules $X{(\ol K)}_{\bfQ}$
and $Y{(\ol K)}_{\bfQ}$ are isomorphic, then the groups $X(L)$ and $Y(L)$ have the same rank
(for all $L$). Brauer theory for representations of finite groups over $\bfQ$-vector spaces
(\cite{23}, \S12.5, corollary of Proposition 25') enables one to establish the following
criterion for the Galois modules $X{(\ol K)}_{\bfQ}$ and $Y{(\ol K)}_{\bfQ}$ to be
isomorphic.

\begin{xxxi}[Theorem] The representations of $G$ in the $\bfQ$-vector spaces $X{(\ol K)}_{\bfQ}$
and $Y{(\ol K)}_{\bfQ}$ are isomorphic if and only if the groups $X(L)$ and $Y(L)$ have the
same rank for all finite algebraic extensions $L$ of $K$.
\end{xxxi}

It is obvious that the representations $X{(\ol K)}_{\bfQ}$ and $Y{(\ol K)}_{\bfQ}$ are
isomorphic when $X$ and $Y$ are isogenous.

The question arises: Are the abelian varieties $X$ and $Y$ isogenous if the representations
$X{(\ol K)}_{\bfQ}$ and $Y{(\ol K)}_{\bfQ}$ are isomorphic? Perhaps this question should be
given a refined formulation, using the Tate height (for example, one could require that an
isomorphism exist which preserves the scalar product).


\medskip{\bf Note added in December 2009}. For algebraic families of complex
abelian varieties the natural analog of the  Tate homomorphism conjecture deals
with the actions of the fundamental group of the base on the integral homology
groups of the fibers \cite{7*,Fa2,11**}.

 In characteristic $p$ there are
variants of the Tate homomorphism conjecture that deal with $p$-divisible
groups (Barsotti--Tate groups) instead of Tate modules  (see \cite{3**},
\cite{10**}).


\medskip{\bf7.} Let $\rho_{\ell}\colon \Gal(\ol K/K)\rightarrow\Aut T_{\ell}(X)$ be a representation of the Galois group in the Tate
module of an abelian variety $X$. Its image $\Im\rho_{\ell}$ is a compact $\ell$-adic Lie
group. The Lie algebra
$$
{\gotg}_{\ell}=\Lie(\Im\rho_{\ell})\subset\End V_{\ell}(X)
$$
is a reductive (because of the semisimplicity of the Tate module, proved by Faltings)
algebraic (see~\cite{Bo1} and~\cite{Bo2}) ${\bfQ}_{\ell}$-Lie algebra containing homotheties.
It is abelian if and only if $X\otimes\ol K$ is an abelian variety of CM-type. There is a
conjecture (Serre~\cite{64}) to the effect that the Lie algebra ${\gotg}_{\ell}$ is
``independent of $\ell$'', i.\,e., there exists a $\bfQ$-algebra $\gotg_0$ such that
${\gotg}_{\ell}\cong\gotg_0\otimes_{\bfQ}{\bfQ}_{\ell}$ for all $\ell$. Concerning the
question of how the Lie algebra $\gotg$ must look, see \cite{29} and~\cite{64}. (A conjecture
of Mumford and Tate asserts that one can take ${\gotg}_0$ to be the Lie algebra of the
Mumford--Tate group of the abelian variety $X$.)

We note that the results in~\cite{73}, combined with Faltings' theorems on the action of the
Galois group on the Tate module, imply that the rank of the Lie algebra ${\gotg}_{\ell}$ does
not depend on $\ell$. The dimension of the center of this Lie algebra is also independent of
$\ell$. We expand the reductive $\bfQ_{\ell}$-Lie algebra ${\gotg}_{\ell}$ as a direct sum
${\gotg}_{\ell}={\gotg}^s\oplus\gotc$ of its center $\gotc$ and a semisimple
$\bfQ_{\ell}$-Lie algebra ${\gotg}^s$. \oripage We consider the finite-dimensional vector
space ${\ol V}_{\ell}(X)=V_{\ell}(X)\otimes_{\bfQ_{\ell}}\ol{\bfQ}_{\ell}$ over the algebraic
closure $\ol{\bfQ}_{\ell}$ of the field ${\bfQ}_{\ell}$ and the semisimple
${\bfQ}_{\ell}$-Lie algebra
$$
\gotg_X={\gotg}^s\otimes_{\bfQ_{\ell}}\ol{\bfQ}_{\ell}\subset\End V_{\ell}(X)\otimes_{\bfQ_{\ell}}\ol{\bfQ}_{\ell}=\End{\ol V}_{\ell}(X).
$$
The space ${\ol V}_{\ell}(X)$ is
 a faithful finite-dimensional ${\gotg}_X$-module. In the
theory of finite-dimensional representations of semisimple Lie algebras in characteristic
0---this is \'Elie Cartan's well-known theory of moduli with highest weight
(Bourbaki~\cite{4a})---it is natural to ask which simple ${\gotg}_X$-modules occur in ${\ol
V}_{\ell}(X)$.

\medskip\textbf{The microweight conjecture~\cite{15a}.} Let $\gotg$ be a simple
${\ol\bfQ}_{\ell}$-Lie algebra which is an ideal in ${\gotg}_X$. Then $\gotg$ is a classical Lie
algebra, i.\,e., a Lie algebra of type $A_m$, $B_m$, $C_m$, or $D_m$. Let $V$ be a nontrivial
simple $\gotg$-submodule of ${\ol V}_{\ell}(X)$. Then $V$ is a fundamental $\gotg$-module.
Moreover, its highest weight is a microweight in the sense of Bourbaki.

If $\gotg$ is a Lie algebra of type $A_m$, then the $\gotg$-module $V$ is isomorphic to the
exterior power $\Lambda^iW$ for a suitable $i$, $1\leq i\leq m$ where $W$ is a standard
$\gotg$-module of dimension $m+1$ for which the image of the Lie algebra $\gotg$ in the algebra
$\End(W)$ coincides with the Lie algebra $\gotsl(W)$ of zero trace operators. If $\gotg$ is a Lie
algebra of type $B_m$, then $V$ is a spinor representation of dimension $2^m$. If $\gotg$ is a
Lie algebra of type $C_m$, then $\dim V=2m$, and the image of $\gotg$ in the algebra $\End(V)$
coincides with the Lie algebra $\gotsp(V)$ of the symplectic group. If $\gotg$ is a Lie algebra
of type $D_m$, then one of the following conditions holds: a) $\dim V=2m$ and the image of
$\gotg$ in $\End(V)$ is the Lie algebra $\gotso(V)$ of the orthogonal group; or b) $V$ is one of
the two semi-spinor representations of dimension $2^{m-1}$.

\medskip\textbf{Note added in December 2009}. The microweight conjecture is proven by Richard
Pink \cite{7**}.

\medskip\textbf{Note in proof (January 1986).} Gerhard Frey~\cite{75} has proved that the nonexistence of
a nontrivial point on the Fermat curve of degree $p$ would follow from the
nonexistence of stable elliptic curves over $\bfQ$ for which $\Delta\geq
N^{2p}/2$ (here $\Delta$ is the minimal discriminant of the curve and $N$ is
its conductor). It was shown by Parshin that for large $p$ the latter assertion
follows from the arithmetic analog of the Bogomolov--Miyaoka--Yao inequality,
see \S1. Concerning the connection of these problems with the Taniyama--Weil
conjecture, see \cite{75}.\footletter{\textit{Added in translation.} See also
Note 1 and~\cite{4*}.}

\textsc{NOTE 1} (added in translation by S.\,Lang). Readers should have the following facts
available in connection with the conjecture that all elliptic curves over the rationals are
modular. Taniyama, in Problem 12 at the Tokyo--Nikko \oripage conference in 1955 stated:
\begin{quote}
Let $C$ be an elliptic curve defined over an algebraic number field $k$, and $L_C(s)$ denote the
$L$-function of $C$ over $k$. If a conjecture of Hasse is true for $\zeta_C(s)$, then the Fourier
series obtained from $L_C(s)$ by the inverse Mellin-transformation must be an automorphic form of
dimension $-2$, of some special type (cf.~Hecke). If so, it is very plausible that this form is
an elliptic differential of the field of that automorphic function. The problem is to ask if it
is possible to prove Hasse's conjecture for $C$, by going back to these considerations, and by
finding a suitable automorphic form from which $L_C(s)$ can be obtained.
\end{quote}
In a letter to me dated 13 August 1986, Shimura wrote:
\begin{quote}
[Taniyama] doesn't say \textit{modular} form. I am sure he was thinking of Hecke's paper No.\,33
(1936) which involves some Euchsian groups not necessarily commensurable with $\SL_2(\bfZ)$... As
for Weil he was far from the conjecture. (It seems that strictly speaking, Weil has never made
the conjecture...) Indeed, in his lecture titled ``On the breeding of bigger and better zeta
functions'' at the University of Tokyo, sometime in August or September 1955, he mentions
Eichler's result and adds: ``But already in the next simplest case, that is the case of an
elliptic curve which cannot be connected with modular functions in Eichler's fashion, the
properties of its zeta function are completely mysterious...''
\end{quote}
Shimura stated the conjecture that all elliptic curves over the rationals are modular in
conversations with Serre and Weil in 1962--1963. In his letter to me, he goes on to say:
\begin{quote}
Knowing the above passage and Taniyama's problem, and having stated the conjecture in my own way,
I couldn't and wouldn't have attributed the origin of the conjecture to Weil. Besides, there is
one point which almost all people seem to have forgotten. In his paper~[1967a], Weil views the
statement as problematic. In other words, he was not completely for it, and so he didn't have to
attribute it to me. Thus there is nothing for which you can take him to task.
\end{quote}
\oripage In another letter dated 16 September 1986, Shimura has written:
\begin{quote}
In these papers [certain papers of Eichler cited previously] it is shown that the zeta function
of an ``arithmetic quotient (especially a modular) curve'' has analytic continuation. The same
applies to the Jacobian... I was conscious of this fact when I talked with Serre. In fact, I
explained about it to Weil, perhaps in 1965. He mentions it at the end of his paper~[1967a]:
``nach eine Mitteilung von G.\,Shimura...'' I even told him at that time that the zeta function
of the curve $C'$ mentioned there is the Mellin transform of the cusp form in question, but he
spared that statement... Of course Weil made a contribution to this subject on his own, but he is
not responsible for the result on the zeta functions of modular elliptic curves, nor for the
basic idea that such curves will exhaust all elliptic curves over $\bfQ$.
\end{quote}
Finally, in a letter to me dated 3 December 1986, Weil stated:
\begin{quote}
Concerning the controversy which you have found fit to raise, Shimura's letters seem to me to put
an end to it, once and for all.
\end{quote}
Copies of all letters and a full dossier of exchanges concerning this history
have been given to the office of the AMS in Providence, and had also been sent
to Parshin and Zarkhin.

\textsc{NOTE 2} (added in translation by S.\,Lang). In the communication, \textit{Valeur,
asymptotique du nombre des points rationnels de hauteur born\'ee sur une courbe \'elliptique}, at
the International Congress of Mathematicians in Edinburgh (1958) N\'eron stated the conjecture
that the height was quadratic. He elaborated an extensive theory proving his conjecture by
showing how to decompose the height into local factors, in two papers: \textit{Mod\'eles minimaux
des vari\'et\'es ab\'eliennes sur les corps locaux et globaux}, Inst.~Hautes \'Etudes
Sci.~Publ.~Math. No.\,21 (1964), and \textit{Quasi-fonctions et hauteurs sur les vari\'et\'es
ab\'eliennes}, Ann.~of Math. (2) {\bf82} (1965). As the Annals paper was being completed, Tate
gave a simple direct global proof, published by Manin (with Tate's permission) in Izv.~Akad.~Nauk
SSSR Ser.~Mat. {\bf28} (1964), 1363--1390 (English transl., Amer.~Math.~Soc.~Transl. (2) {\bf59}
(1966), 82--110). Tate's proof and N\'eron's proofs were done independently and simultaneously.
Both used the limiting trick, Tate directly on the height and N\'eron to show the existence of
his local symbol, before identifying it with the intersection number on the minimal model. In
1964 I had both manuscripts available, and I reproduced both Tate's proof and that part of
N\'eron's paper in my report to the Bourbaki seminar in May 1964.

\oripage

\def\same{\rule{10mm}{0.4pt}\;}

\newpage

\centerline{\bf Bibliography}

\medskip

\textbf{\normalsize A. Items from Lang's bibliography cited in this Appendix}

\def\refname{}

Translated by N.\,KOBLITZ

\end{document}